\documentclass[12pt, letter]{amsart}

\DeclareRobustCommand{\SkipTocEntry}[5]{}
\usepackage{etoolbox}

\makeatletter
\let\old@tocline\@tocline
\let\section@tocline\@tocline
\newcommand{\subsection@dotsep}{4.5}
\newcommand{\subsubsection@dotsep}{4.5}
\patchcmd{\@tocline}
  {\hfil}
  {\nobreak
     \leaders\hbox{$\m@th
        \mkern \subsection@dotsep mu\hbox{.}\mkern \subsection@dotsep mu$}\hfill
     \nobreak}{}{}
\let\subsection@tocline\@tocline
\let\@tocline\old@tocline

\patchcmd{\@tocline}
  {\hfil}
  {\nobreak
     \leaders\hbox{$\m@th
        \mkern \subsubsection@dotsep mu\hbox{.}\mkern \subsubsection@dotsep mu$}\hfill
     \nobreak}{}{}
\let\subsubsection@tocline\@tocline
\let\@tocline\old@tocline

\let\old@l@subsection\l@subsection
\let\old@l@subsubsection\l@subsubsection

\def\@tocwriteb#1#2#3{%
  \begingroup
    \@xp\def\csname #2@tocline\endcsname##1##2##3##4##5##6{%
      \ifnum##1>\c@tocdepth
      \else \sbox\z@{##5\let\indentlabel\@tochangmeasure##6}\fi}%
    \csname l@#2\endcsname{#1{\csname#2name\endcsname}{\@secnumber}{}}%
  \endgroup
  \addcontentsline{toc}{#2}%
    {\protect#1{\csname#2name\endcsname}{\@secnumber}{#3}}}%

\newlength{\@tocsectionindent}
\newlength{\@tocsubsectionindent}
\newlength{\@tocsubsubsectionindent}
\newlength{\@tocsectionnumwidth}
\newlength{\@tocsubsectionnumwidth}
\newlength{\@tocsubsubsectionnumwidth}
\newcommand{\settocsectionnumwidth}[1]{\setlength{\@tocsectionnumwidth}{#1}}
\newcommand{\settocsubsectionnumwidth}[1]{\setlength{\@tocsubsectionnumwidth}{#1}}
\newcommand{\settocsubsubsectionnumwidth}[1]{\setlength{\@tocsubsubsectionnumwidth}{#1}}
\newcommand{\settocsectionindent}[1]{\setlength{\@tocsectionindent}{#1}}
\newcommand{\settocsubsectionindent}[1]{\setlength{\@tocsubsectionindent}{#1}}
\newcommand{\settocsubsubsectionindent}[1]{\setlength{\@tocsubsubsectionindent}{#1}}

\renewcommand{\l@section}{\section@tocline{1}{\@tocsectionvskip}{\@tocsectionindent}{}{\@tocsectionformat}}%
\renewcommand{\l@subsection}{\subsection@tocline{1}{\@tocsubsectionvskip}{\@tocsubsectionindent}{}{\@tocsubsectionformat}}%
\renewcommand{\l@subsubsection}{\subsubsection@tocline{1}{\@tocsubsubsectionvskip}{\@tocsubsubsectionindent}{}{\@tocsubsubsectionformat}}%
\newcommand{\@tocsectionformat}{}
\newcommand{\@tocsubsectionformat}{}
\newcommand{\@tocsubsubsectionformat}{}
\expandafter\def\csname toc@1format\endcsname{\@tocsectionformat}
\expandafter\def\csname toc@2format\endcsname{\@tocsubsectionformat}
\expandafter\def\csname toc@3format\endcsname{\@tocsubsubsectionformat}
\newcommand{\settocsectionformat}[1]{\renewcommand{\@tocsectionformat}{#1}}
\newcommand{\settocsubsectionformat}[1]{\renewcommand{\@tocsubsectionformat}{#1}}
\newcommand{\settocsubsubsectionformat}[1]{\renewcommand{\@tocsubsubsectionformat}{#1}}
\newlength{\@tocsectionvskip}
\newcommand{\settocsectionvskip}[1]{\setlength{\@tocsectionvskip}{#1}}
\newlength{\@tocsubsectionvskip}
\newcommand{\settocsubsectionvskip}[1]{\setlength{\@tocsubsectionvskip}{#1}}
\newlength{\@tocsubsubsectionvskip}
\newcommand{\settocsubsubsectionvskip}[1]{\setlength{\@tocsubsubsectionvskip}{#1}}

\patchcmd{\tocsection}{\indentlabel}{\makebox[\@tocsectionnumwidth][l]}{}{}
\patchcmd{\tocsubsection}{\indentlabel}{\makebox[\@tocsubsectionnumwidth][l]}{}{}
\patchcmd{\tocsubsubsection}{\indentlabel}{\makebox[\@tocsubsubsectionnumwidth][l]}{}{}

\newcommand{\@sectypepnumformat}{}
\renewcommand{\contentsline}[1]{%
  \expandafter\let\expandafter\@sectypepnumformat\csname @toc#1pnumformat\endcsname%
  \csname l@#1\endcsname}
\newcommand{\@tocsectionpnumformat}{}
\newcommand{\@tocsubsectionpnumformat}{}
\newcommand{\@tocsubsubsectionpnumformat}{}
\newcommand{\setsectionpnumformat}[1]{\renewcommand{\@tocsectionpnumformat}{#1}}
\newcommand{\setsubsectionpnumformat}[1]{\renewcommand{\@tocsubsectionpnumformat}{#1}}
\newcommand{\setsubsubsectionpnumformat}[1]{\renewcommand{\@tocsubsubsectionpnumformat}{#1}}
\renewcommand{\@tocpagenum}[1]{%
  \hfill {\mdseries\@sectypepnumformat #1}}

\let\oldappendix\appendix
\renewcommand{\appendix}{%
  \leavevmode\oldappendix%
  \addtocontents{toc}{%
    \protect\settowidth{\protect\@tocsectionnumwidth}{\protect\@tocsectionformat\sectionname\space}%
    \protect\addtolength{\protect\@tocsectionnumwidth}{2em}}%
}
\makeatother

\makeatletter
\settocsectionnumwidth{2em}
\settocsubsectionnumwidth{2.5em}
\settocsubsubsectionnumwidth{3em}
\settocsectionindent{1pc}%
\settocsubsectionindent{\dimexpr\@tocsectionindent+\@tocsectionnumwidth}%
\settocsubsubsectionindent{\dimexpr\@tocsubsectionindent+\@tocsubsectionnumwidth}%
\makeatother

\settocsectionvskip{5pt}
\settocsubsectionvskip{0pt}
\settocsubsubsectionvskip{0pt}

\settocsectionformat{\scshape}%
\settocsubsectionformat{\mdseries}
\settocsubsubsectionformat{\mdseries}
\setsectionpnumformat{\scshape}%
\setsubsectionpnumformat{\mdseries}
\setsubsubsectionpnumformat{\mdseries}

\let\oldtableofcontents\tableofcontents
\renewcommand{\tableofcontents}{%
  \vspace*{-\linespacing}%
  \oldtableofcontents}

\usepackage{graphicx} %
\usepackage[caption=false]{subfig}
\usepackage[margin=1in,marginparwidth=0.8in, marginparsep=0.1in]{geometry}
\newcommand{\bA}{\mathbb{A}}

\newcommand{\cP}{\mathcal{P}}
\newcommand{\bC}{\mathbb{C}}

\newcommand{\cF}{\mathcal{F}}

\newcommand{\fs}{\mathfrak{s}}

\newcommand{\Bl}{\mathrm{Bl}}

\DeclareMathOperator{\Spec}{Spec}

\newcommand{\cN}{\mathcal{N}}

\newcommand{\cA}{\mathcal{A}}

\newcommand{\sU}{\mathscr{U}}

\newcommand{\bR}{\mathbb{R}}

\newcommand{\bZ}{\mathbb{Z}}

\newcommand{\cO}{\mathcal{O}}

\DeclareMathOperator{\Fuk}{Fuk}

\DeclareMathOperator{\QCoh}{QCoh}

\DeclareMathOperator{\End}{End}
\DeclareMathOperator{\Hom}{Hom}

\newcommand{\cU}{\mathcal{U}}
\newcommand{\cX}{\mathcal{X}}

\newcommand{\Perf}{\mathrm{Perf}}

\newcommand{\red}[1]{\textcolor{red}{#1}}

\let\Im\relax
\DeclareMathOperator{\Im}{Im}

\DeclareMathOperator{\Mod}{Mod}

\newcommand{\R}{\mathbb R}
\newcommand{\C}{\mathbb C}
\newcommand{\id}{\mathrm{id}}

\DeclareMathOperator{\coker}{coker}

\DeclareMathOperator{\Crit}{Crit}

\DeclareMathOperator{\Skel}{Skel}

\DeclareMathOperator{\supp}{supp}
\DeclareMathOperator{\rank}{rank}
\DeclareMathOperator{\diag}{diag}
\DeclareMathOperator{\Mat}{Mat}

\DeclareMathOperator{\fib}{fib}

\newcommand{\laff}{\lambda_{\mathrm{aff}}}

\newcommand{\T}{\mathsf T}

\newcommand{\dd}{\mathrm{d}}

\usepackage[dvipsnames]{xcolor}
\usepackage[draft]{say}
\definecolor{MHcol}{RGB}{100, 180, 255}

\definecolor{egyptianblue}{rgb}{0.06, 0.2, 0.65}

\usepackage{amsthm}
\usepackage{amssymb}
\usepackage{amsmath}
\usepackage{parskip}
\usepackage{tikz}
\usepackage{tikz-cd}
\usepackage{bbm}
\usepackage{enumitem}
\usepackage{xcolor}
\usepackage{import}
\usepackage{yhmath}
\usepackage{aliascnt}
\usepackage{standalone}
\usepackage{caption}
\usepackage{subcaption}
\usepackage{mathtools}

\usepackage{mathrsfs}

\usepackage{amsmath,amssymb,mathtools,amsthm,mathrsfs}
\usepackage{tikz-cd,enumitem,booktabs,longtable,array}
\usepackage[colorlinks=true,linkcolor=blue,citecolor=blue,urlcolor=blue]{hyperref}

\newtheorem{theorem}{Theorem}[section]
\newtheorem{proposition}[theorem]{Proposition}
\newtheorem{lemma}[theorem]{Lemma}
\newtheorem{corollary}[theorem]{Corollary}
\theoremstyle{definition}
\newtheorem{definition}[theorem]{Definition}
\newtheorem{remark}[theorem]{Remark}
\newtheorem{example}[theorem]{Example}
\newtheorem{warning}[theorem]{Warning}

\newcommand{\Z}{\mathbb Z}
\newcommand{\Q}{\mathbb Q}

\newcommand{\OO}{\mathcal O}

\newcommand{\HH}{\operatorname{HH}}
\newcommand{\SH}{\operatorname{SH}}
\newcommand{\D}{\mathcal D}

\newcommand{\Ac}{\mathcal A}
\newcommand{\PP}{\mathcal P}
\newcommand{\U}{\mathscr U}
\newcommand{\A}{\mathscr A}

\newcommand{\Nn}{\mathcal N}
\newcommand{\pr}{\mathrm{prin}}

\newcommand{\fr}{\mathrm{fr}}

\newcommand{\otL}{\mathbin{\otimes^{\mathbb L}}}
\DeclareMathOperator{\Ind}{Ind}

\DeclareMathOperator{\RHom}{RHom}

\DeclareMathOperator{\Supp}{Supp}

\DeclareMathOperator{\Vit}{Vit}

\setlist[enumerate]{label=\textup{(\roman*)}}
\numberwithin{equation}{section}
\title{Ice and mirrors}
\author{Mingyuan Hu}
\author{Vivek Shende}
\begin{document}
\begin{abstract}
Given a mirror pair of smooth affine principal upper cluster varieties admitting reddening sequences, we show that homological mirror symmetry for principal cluster varieties -- assumed compatible with certain cluster structures -- descends under freezing.

This allows us to induct on the following class.  A skew-symmetrizable matrix is peelable if, by mutations and removing sinks or sources, it can be reduced to the empty matrix;  it is sub-peelable if it is a principal submatrix of a peelable matrix.  We establish  homological mirror symmetry for full rank upper cluster varieties whose mutable matrix is sub-peelable.

We show sub-peelability for the
Fock--Goncharov $\mathrm{SL}_n$  quiver for any compact surface with nonempty boundary, all of whose boundaries have markings.  Consequently, we establish  homological mirror symmetry between the moduli of decorated
$\mathrm{SL}_n$ local systems and the corresponding pinned
$\mathrm{PGL}_n$ spaces. 
The results were obtained using ChatGPT Astra.

\end{abstract}

\maketitle

\setlength{\parskip}{0pt}

\renewcommand\contentsname{\vspace*{-40pt}}
\tableofcontents

\thispagestyle{empty}

\setlength{\parskip}{.3\baselineskip plus 2pt}

\newpage

\section{Introduction}
\label{sec:intro}

Cluster algebra is organized around two kinds of variables: the mutable and the frozen.  A basic operation in the theory changes one of the former into one of the latter, which (under appropriate hypotheses) corresponds to passing to the Zariski open subset of the cluster variety characterized by the invertibility of the previously mutable variable.
The main observation of the present article is that, in good circumstances, homological mirror symmetry descends under  freezing.

This new observation allows us to significantly extend the reach of our previous work \cite{HS}.  There we established homological mirror symmetry for the cluster varieties with full-rank sink-recurrent seeds -- this class being inductively defined by asking that removing a sink, and also removing the sink and its neighbors, both yield seeds which are mutation equivalent to sink-recurrent seeds.  Our new freezing result allows us to establish homological mirror symmetry for a larger class.

\begin{definition}
\label{def:subpeelable}
The class of \emph{peelable} matrices is the smallest class of skew symmetrizable matrices closed under mutation, closed under adjoining a source or sink, and containing the empty matrix.
A matrix is \emph{sub-peelable} if it is obtained from a peelable matrix by deleting rows and columns with matching indices.
\end{definition}

The class of peelable matrices was introduced in \cite{Ladkani}, and there denoted $\cP'$.  It  contains the class of Louise seeds \cite{EJ}.  Sink-recurrent seeds are obviously peelable.
Our interest in this class stems from the following (new) example:

\begin{theorem}[\ref{thm:surface-quiver}] \label{thm:fg-sub-peelable}
    Let $S$ be a compact surface with nonempty boundary, for which all boundaries have marked points, and there are at least three marked points if $S$ is a disk.  Then the Fock-Goncharov seed for the moduli of decorated local systems $\mathcal{A}^\times_{\mathrm{SL}_n, S}$ is sub-peelable.
\end{theorem}

When $n=2$ and there are at least two marked points in total, these seeds are sink-recurrent \cite[Proposition~D.1]{HS}, but the once-holed torus with one boundary marked point is sub-peelable but not sink-recurrent (Example~\ref{ex:strict}).  
We note by contrast that the Fock-Goncharov $\mathrm{SL}_2$ quiver for the once punctured torus is not sub-peelable (Example \ref{nonexample}).

In this article, we prove homological mirror symmetry for all upper cluster varieties with full rank sub-peelable seeds.

We recall the notational conventions of cluster algebra.
A square integer matrix \(B\) is said to be {\em skew-symmetrizable} if there is a positive diagonal matrix $D$ with $BD = - DB^{\mathsf T}$.
A column index $s$ is a \emph{sink} if $B_{is}\ge0$ for every $i$, and a
\emph{source} if $B_{is}\le0$ for every $i$.  This specializes to the intuitive meaning of these words when $B$ is the matrix of a quiver.  A {\em seed} $\fs$ is the data of a skew symmetrizable $B$, along with:
\[
 \widetilde B=\binom BC,
 \qquad 
 C\in\operatorname{Mat}_{f\times m}(\Z).
\]
We call
$\widetilde B$ \emph{full rank} if $\rank_{\Q}\widetilde B=m$.

We write $\mathscr{A}(\widetilde B)$ for the cluster algebra (generated by the cluster variables in any seed) and $\mathscr{U}(\widetilde B)$ for the upper cluster algebra (elements of the fraction field  regular in all seeds).  In many situations -- in particular, all locally acyclic cluster algebras, which includes the sink-recurrent class -- one has $\mathscr{A}(\widetilde B) = \mathscr{U}(\widetilde B)$; however, we do not assume this here, and will work with the upper cluster algebra.  We will find it convenient to systematically work with the `principal coefficients' of \cite{FZ-IV}.  We write:
\[
 \mathcal{U}(\widetilde B)=\Spec\U(\widetilde B),\qquad
 \mathcal{U}_{\pr}(B)=\Spec\U\binom B{I_m}.
\]
The group $H=(\C^*)^m$ acts on the coordinates $(X;Z;U) = (X_1,\ldots,X_m;Z_1,\ldots,Z_m;U_1,\ldots,U_f)$ 
in a seed of $\mathcal{U}_\pr(-B^{\mathsf T})\times(\C^*)^f$ by the weights
\begin{equation}\label{eq:weights}
 [\,I_m\mid B^{\mathsf T}\mid C^{\mathsf T}\,],
\end{equation}
and homogeneity of the exchange relations extends this action to the upper
cluster variety.  The expected mirror stack is:
\begin{equation}\label{eq:mirror}
  \mathcal{X}(\widetilde B^\vee)
 :=\left[
   \bigl(\mathcal{U}_{\pr}(-B^{\mathsf T})\times(\C^*)^f\bigr)/H
  \right]
\end{equation}
The coarse moduli space of this stack is the affinization of the $\cX$ variety considered in \cite{FG09, FG, GHK, GHKK}, and agrees with the stack considered in \cite{Goncharov-Shen}. 

The seed $\fs$ with matrix $\widetilde B$ also determines the embedding of a parameterized torus $j_{\fs}: (\C^*)^{m+f} \hookrightarrow \cX(\widetilde{B}^\vee)$ by asking that the coordinates $\chi_i$ on the torus are:
\begin{equation} \label{eq: x variables}
 \chi_i=Z_i\prod_{j=1}^mX_j^{-B_{ij}}\quad(1\le i\le m),\qquad
 \chi_{m+a}=U_a\prod_{j=1}^mX_j^{-C_{aj}}\quad(1\le a\le f).
\end{equation}
Mutate $\fs$ at the $k$'th variable to obtain $\fs'$ with matrix $\widetilde B'$. Then the comparison
$$j_{\fs}: (\C^*)^{m+f} \hookrightarrow \cX(\widetilde{B}^\vee) = \cX((\widetilde{B}')^\vee) \hookleftarrow  (\C^*)^{m+f} :j_{\fs'}$$
is, in coordinates 
\[
 \chi'_k=\chi_k^{-1},\qquad
 \chi'_i=\chi_i\chi_k^{[-\widetilde B_{ik}]_+}
                    (1+\chi_k)^{\widetilde B_{ik}}\quad(i\ne k),
 \qquad [r]_+=\max(r,0).
\]
This formula is  dual to the more usual $\cA$ cluster transformation.
These quotient coordinates and their mutation are given in
\cite[Section~3.6, equations~(25) and~(28)]{HS}. 
The original $\cX$ mutation is \cite[Section~1.2, equation~(13)]{FG09}.
Their convention is
$X'_i=X_i(1+X_k^{-\operatorname{sgn}(\varepsilon_{ik})})^{-\varepsilon_{ik}}$
for $i\ne k$. Setting $X_i=\chi_i$ and
$\varepsilon_{ik}=-\widetilde B_{ik}$ gives the formula above.

The frozen coordinates define a map
$y:\mathcal U(\widetilde B)\to(\C^*)^f$. We write
$
 T_{\fr}(\widetilde B):=\operatorname{Hom}(\Z^f,\C^*)$.
For $\kappa\in T_{\fr}(\widetilde B)$, let $\mathcal L_\kappa$
be the rank-one local system on $(\C^*)^f$ with these monodromies.  
$T_{\fr}(\widetilde B)$ acts on the Fukaya category by
$(L,\mathcal E)\longmapsto
 (L,\mathcal E\otimes(y|_L)^*\mathcal L_\kappa)$.  The same torus acts on $\mathcal{X}(\widetilde B^\vee)$
by $U_a\mapsto\kappa_aU_a$, fixing the principal factor. This
commutes with $H$ and therefore defines an action on the quotient
stack.  The role of these actions in mirror symmetry is explained previously in \cite{GHK, GHKK, GS}.

Positivity of cluster transformations gives a well defined positive real locus $L_+ := \mathcal{U}(\widetilde B)(\bR_+)$.  For a seed $\fs$ obtained from the initial seed by iterated mutations, we write $L_{\fs} \subset (\C^*)^{m+f} \subset \cU(\widetilde B)$ for the unitary locus of the corresponding seed torus.

\begin{definition}
Fix an exact symplectic form $\omega$ on $\cU(\widetilde B)$.
We say that a Liouville form $\lambda$ on $\cU(\widetilde B)$ is
\emph{cluster-admissible} for $\omega$ if
$\dd\lambda=\omega$, and $\lambda$ is of finite type and weakly Weinstein, $L_+$ is a properly
embedded strictly exact Lagrangian, and every unitary seed torus $L_{\fs}$
is an exact Lagrangian, and in addition the following holds. 

Write $q_{\fs}=L_{\fs}\cap L_+$, the point with all seed coordinates equal to $1$.
Then there is a cotangent disk embedding $i_{\fs}:D^*_{\epsilon_{\fs}}L_{\fs}\hookrightarrow
\cU(\widetilde B)$ and a compactly supported function $h_{\fs}$ such that
\[
 i_{\fs}|_{L_{\fs}}=\mathrm{id},\qquad
 i_{\fs}^{-1}(L_+)=D^*_{q_{\fs},\epsilon_{\fs}}L_{\fs},\qquad
 h_{\fs}|_{L_+}=0,\qquad
 i_{\fs}^*(\lambda+\dd h_{\fs})=\lambda_{\mathrm{can}}.
\]
Here $\lambda_{\mathrm{can}}=\sum_i p_i\,\dd q_i$ is the cotangent Liouville form.
\end{definition}

For principal coefficients we take $\omega=\Im\Omega_B$, where
$\Omega_B$ is given by \eqref{eq:principal-omega}.
Its extension to a regular algebraic symplectic form on the smooth
finite-type variety $\cU_{\pr}(B)$ is proved in Lemma~\ref{lem: form}.  In this case we construct cluster-admissible forms in Lemma~\ref{lem:principal-geometry} whenever $B$ has large cluster complex 
in the sense of \cite{GHKK}; this happens in particular when $B$ has a reddening sequence 
\cite[Lemma~8.29]{GHKK}.

Lemma~\ref{lem:change-primitive} below gives an equivalence of Fukaya categories for any
two such forms, preserving the object $L_+$ and compatibility with all seed restrictions.

Fix such neighborhoods and the deformations $\lambda+t\,\dd h_{\fs}$.
They preserve $\lambda|_{L_+}=0$ and give Viterbo functors
$
 \Vit_{\fs}:\Fuk(\cU(\widetilde B))\to\Fuk(T^*L_{\fs})$ and isomorphisms $
 \Vit_{\fs}(L_+)\simeq F_{\fs}:=T^*_{q_{\fs}}L_{\fs}.
$

For any torus $M$ and $q\in M$, a cotangent fiber $F_q$ gives an equivalence
\[
\chi:\Fuk(T^*M)\xrightarrow{\sim}\Perf(\C[H_1(M;\Z)])
\]
and a fixed isomorphism 
 $\chi(F_q)\simeq\C[H_1(M;\Z)]$,
by \cite[Theorem~1.1]{Ab}. (For a torus $M$, one has 
$C_{-*}(\Omega_qM) \cong \C[H_1(M;\Z)]$).

The coordinate circles of $L_{\fs}$ give a basis
$
 A_{\fs}:=\C[H_1(L_{\fs};\Z)]
     =\C[t_1^{\pm1},\ldots,t_{m+f}^{\pm1}]$. 
We write 
$ j_{\fs}:\Spec A_{\fs}\hookrightarrow\cX(\widetilde B^\vee)$ for the cluster chart characterized by 
$ j_{\fs}^*\chi_{\fs,i}=t_i$, where  $\chi_{\fs,i}$ denotes the quotient coordinate \eqref{eq: x variables} for the seed $\fs$. We also write $\chi_{\fs}: \Fuk(T^*L_{\fs}) \xrightarrow{\sim} \Perf( A_{\fs})$, and 
$
 \kappa_{\fs}:\chi_{\fs}\Vit_{\fs}(L_+)\xrightarrow{\sim}A_{\fs}
$ for the aforementioned fixed identification.

\begin{definition}\label{def: cms}
By \emph{cluster homological mirror symmetry} for $\widetilde B$, we mean a
cluster-admissible Liouville form
on $\cU(\widetilde B)$, an equivalence of categories $ \Phi:\Perf(\cX(\widetilde B^\vee))
\xrightarrow{\sim}\Fuk(\cU(\widetilde B))$, a fixed isomorphism $ \eta:\Phi(\cO)\xrightarrow{\sim}L_+$, and, for every seed $\fs$, we require a natural transformation $\alpha_{\fs}:\chi_{\fs}\Vit_{\fs}\Phi\xRightarrow{\sim}j_{\fs}^*$ with $\alpha_{\fs,\cO}=\kappa_{\fs}\circ\chi_{\fs}\Vit_{\fs}(\eta)$, making the following diagram commute:

\begin{equation}\label{eq:cms}
\begin{tikzcd}[column sep=large]
 \Perf(\cX(\widetilde B^\vee))\arrow[r,"\Phi"]\arrow[d,"j_{\fs}^*"']
 &\Fuk(\cU(\widetilde B))\arrow[d,"\Vit_{\fs}"]\\
 \Perf(A_{\fs})\arrow[r,"\chi_{\fs}^{-1}"']&\Fuk(T^*L_{\fs})
\end{tikzcd}
\end{equation}

By \emph{equivariant cluster homological mirror symmetry}, we mean such data with a
$T_{\fr}(\widetilde B)$-equivariant structure on $\Phi$, compatible with
the comparisons $\alpha_{\fs}$ and the unit and multiplication of the torus.
\end{definition}

When  $\cX(\widetilde B^\vee)$ is (nonstacky)  affine, as is the case for principal coefficients, the required coherences are properties, rather than  structures.  Indeed, for commutative rings $A, B$,  a map $f: \Perf\, A \to \Perf\, B$ along with an isomorphism $f(A) \cong B$ induces 
$$A = \Hom_A(A, A) \to \Hom_B(f(A), f(A)) = \Hom_B(B, B) = B;$$
one readily sees that the original $f$ is extension of scalars along this map $\overline{f}: A \to B$.  Thus, given two such maps, $f, g: \Perf\, A \to \Perf\, B$, either $\overline{f} = \overline{g}$, or not.  Thus $\alpha$ is unique (up to isomorphism) if it exists. Similarly, in the affine case, since all the isomorphisms are fixed and the left vertical arrow is extension of scalars along a localization of a domain, an equivalence $\Phi$ and isomorphism $\eta: \Phi(\mathcal{O}) \cong L_+$ compatible with a single $\kappa_{\fs}$ is unique up to unique isomorphism if it exists.  Uniqueness here includes
the isomorphism $\eta$ and the normalization of each $\alpha_{\fs}$.
In particular, if $\widetilde{B}$ satisfies cluster HMS, then the aforementioned equivalence of Fukaya categories for different cluster-admissible forms constructed ad-hoc in  Lemma \ref{lem:change-primitive} is unique up
to unique isomorphism by Lemma~\ref{lem:normalized-hms-unique}.

Finally, while it may seem daunting to verify commutativity for all seeds, the works \cite{GLe, PT} calculate the Floer theory for the unitary Lagrangians in mutation-adjacent seeds.  We will also explain in Lemma \ref{lem:one-seed-hms} below how to deduce this essentially formally from the simplest case $\cN := \cU_{\pr}([0])$.  This calculation implies that, given an equivalence $\Phi$ and fixed isomorphism $\eta: \Phi(\cO) \cong L_+$, and remaining in the aforementioned affine setting, if \eqref{eq:cms} holds for one seed, it holds for them all.

While $\cU(\widetilde{B})$, if smooth of finite type, always carries its K\"ahler form, which has a Weinstein primitive, the tori $L_{\fs}$ are not Lagrangian for this form.  At least for principal coefficients, there is a natural holomorphic symplectic form $\Omega_B$ (see Lemma \ref{lem: form})  with respect to whose imaginary part they are Lagrangian~\cite{GSV05, FG09};  said form however does not obviously admit an appropriate primitive.  To construct the desired form, we will  need a certain hypothesis, which we now recall.

For a skew-symmetrizable $m\times m$ matrix $B$, a
\emph{reddening sequence} (see \cite[Section~3]{Keller-DT} and
\cite[Definition~1.0.4]{BHIT})
is a finite sequence of mutations for which
\[
 \mu_{k_r}\cdots\mu_{k_1}\binom{B}{I_m}
 =\binom{B_r}{C_r},\qquad (C_r)_{ij}\le0\quad\text{for all }i,j.
\]

Every sub-peelable matrix has a reddening sequence by Lemma~\ref{lem:class}.

The key new mirror symmetry result of this article is the following: 

\begin{theorem}[\ref{thm:freezing}] \label{thm: intro freezing}
Assume that $B$ has a reddening sequence.
Suppose 
$\cU_{\pr}(B)$ and $\cU_{\pr}(B^\vee)$ are smooth affine of finite type, and that $\cU_{\pr}(B)$ satisfies the hypotheses of 
Lemma~\ref{lem:principal-geometry}, and, equipped with the resulting symplectic form, satisfies cluster homological mirror symmetry.

Then for any $I\subset \{1,\ldots,m\}$, the varieties
$\cU_{\pr}(B_{I,I})$ and $\cU_{\pr}((B^\vee)_{I,I})$ have the corresponding geometric properties, and satisfy
cluster homological mirror symmetry.
\end{theorem}

The key point in the proof of Theorem \ref{thm: intro freezing} is to understand when Viterbo restriction, a priori a categorical localization, is carried to a classical commutative-algebra localization on the mirror.  When the mirror is affine and smooth of finite type, we give sufficient conditions for this in Theorem \ref{thm:affine-freezing} below.  Certain of these conditions will be verified by the  cluster algebraic structure of the problem at hand; the remaining condition is that the symplectic cohomology is finite over its zero degree component.  This we will verify using the methods of Ganatra and Pomerleano \cite{GP, Pom}, who have described the symplectic cohomology of an affine variety in terms of the ordinary cohomology of strata of an appropriate compactification.  There is a subtle point: cluster mirror symmetry is natural for the imaginary part of the holomorphic symplectic form (for which the compact seed tori are Lagrangian), whereas the Ganatra-Pomerleano machinery applies to the K\"ahler form.  Fortunately, we can show these are symplectomorphic using the criteria of \cite{HSinv}.\footnote{We also constructed such symplectomorphisms in \cite{HS}, but they were not fundamental to that article: we could have just declared there that we would work with $\Im \Omega$, at the cost of proving ad-hoc that it was Weinstein.  For our argument here, the use of both symplectic forms is essential.}

We use a de/equivariantization argument to show that mirror symmetry in the principal case implies equivariant mirror symmetry for any cluster variety with the same mutation matrix:
\begin{theorem}[\ref{thm:coefficients}] \label{thm:coefficients-intro}
Suppose $\mathcal U_{\pr}(B)$ and
$\mathcal U_{\pr}(-B^{\mathsf T})$ are smooth affine of finite type.
Suppose cluster mirror symmetry for $\mathcal U_{\pr}(B)$ is realized by
a Liouville form $\lambda_{\pr}$ satisfying the conclusions of
Lemma~\ref{lem:principal-geometry}.
Then for  $C\in\operatorname{Mat}_{f\times m}(\Z)$ such that $\widetilde{B} = \binom BC$ has full rank, there is  equivariant
cluster mirror symmetry, and in particular,
$\Fuk(\mathcal U(\widetilde B))
 \simeq
 \Perf( \mathcal X(\widetilde B^\vee) )$.
\end{theorem}

Using these results, we generalize our argument from \cite{HS} to prove:

\begin{theorem}[\ref{thm:sub-peelable}]\label{thm:main} 
Let $\widetilde B=(B;C)$ have full rank, where $B$ is
skew-symmetrizable and sub-peelable.
Then $\cU(\widetilde B)$ is smooth affine of finite type and
carries a cluster-admissible Liouville structure symplectomorphic to the Liouville form from the affine K\"ahler structure. $\cX(\widetilde B^\vee)$ is a smooth Deligne-Mumford stack, and is a smooth affine scheme if $\Z^{m+f}/\widetilde B \Z^m$ is torsion free.\footnote{This condition is called `really full rank' in the cluster algebra literature}   These spaces satisfy equivariant cluster mirror symmetry; in particular, $\Fuk(\cU(\widetilde B)) \cong \Perf(\cX(\widetilde B^\vee))$.
\end{theorem}

A series of works have shown that the moduli of decorated local systems introduced in \cite{FG, GS} are in fact the cluster varieties associated to the seeds proposed in the original works;  for the Fock-Goncharov seed for $\mathrm{SL}_n$ on a marked surface $S$, we have 
$\cU(\widetilde B)=\Ac^\times_{\mathrm{SL}_n,S}$ and
$\cX(\widetilde B^\vee)=\PP_{\mathrm{PGL}_n,S}$; see
\cite[Lemma~3.7, Remarks~3.8, 3.10, Theorems~3.21, 4.3 and Lemma~4.4]{IOS},
\cite[Lemma~3.4, Theorem~3.15, Remark~3.17 and Section~4.6]{IO}, and
\cite[Theorem~1.1]{Shen}.
These results, together with Theorems~\ref{thm:fg-sub-peelable} and~\ref{thm:main}, give:

\begin{corollary}\label{cor:FG}
Let $(S,M)$ be a compact connected oriented marked surface with
nonempty boundary, no punctures, at least one mark on each boundary
component, and at least three marked points if the surface is a disk.
Then there is an equivariant cluster mirror symmetry $\Fuk(\Ac^\times_{\mathrm{SL}_n,S})
       \simeq\Perf(\PP_{\mathrm{PGL}_n,S})$. 
\end{corollary}
This establishes the natural generalization to surfaces of \cite[Conjecture 1.5]{Goncharov-Shen} for $\mathrm{SL}_n$.  
(We note in this context also that our Theorem \ref{thm:coefficients-intro} upgrades our proof in \cite{HS} of the nonequivariant version of \cite[Conjecture 1.5]{Goncharov-Shen} to the full version.) 

\begin{remark}
    One can also obtain an equivalence
$\Fuk(\PP_{\mathrm{PGL}_n,S})
 \simeq\Perf(\Ac^\times_{\mathrm{SL}_n,S})$ by applying the finite descent argument
of \cite[Section~5.5]{HS} to the quotient presentation in
Proposition~\ref{prop:upper-finite-atlas}, and similarly such an isomorphism exists in the context of Theorem \ref{thm:main} whenever the $\cX$ space is nonstacky, hence smooth affine. 
\end{remark}

\begin{remark}
    We have written Fukaya categories with $\C$ coefficients and taken algebraic varieties over $\C$.  The need to work over a field stems only from the (essential) use of \cite{Pom}.   In particular, the splitting construction in
\cite[Section~4.3]{Pom} divides by factorials, and its integral extension
is left open in \cite[Remark~4.23]{Pom}.
\end{remark}

We work throughout with stable $k$-linear categories, all functors in this context are understood to respect these structures.

{\bf AI Statement.} 
The results of this article are obtained after weeks of discussion with ChatGPT-Astra, building on and inputting our previous (also AI assisted) work \cite{HSinv, HS} along with many previous discussions of the first-named author with AI. 
The key new results, namely Theorem \ref{thm:fg-sub-peelable} and Theorem \ref{thm:main} were produced by Astra, though the idea to use the methods of \cite{ShendeToric}  to get equivariance and go from the principal to  the general case is essentially due to us.  Aside from this, we asked the questions,  improved the exposition, and formulated the results -- and  the definitions -- in their present forms.

{\bf Acknowledgements.} 
We found \href{https://webusers.imj-prg.fr/~bernhard.keller/quivermutation/}{Keller's quiver mutation applet} very helpful.
We are supported by Villum Fonden Villum Investigator 37814.

\section{Sub-peelability of Fock-Goncharov quivers of type A}
\label{sec:surface-quivers}

Let $S$ be a compact connected oriented surface, and let
$M\subset\partial S$ be finite. Assume
\begin{equation}\label{eq:surface-assumption}
 \partial S\ne\varnothing,\qquad
 M\cap C\ne\varnothing\quad\text{for every boundary component }C,
 \qquad -2\chi(S)+|M|>0.
\end{equation}
There are no punctures. For a disk the last condition is
$|M|\ge3$; otherwise it is automatic. Fix $n\ge2$.
We write $\Ac^\times_{\mathrm{SL}_n,S}$ for the space of
boundary-generic twisted decorated local systems. All boundary
frozen coordinates are inverted.

\begin{theorem}\label{thm:surface-quiver}
The Fock--Goncharov cluster structure on
$\Ac^\times_{\mathrm{SL}_n,S}$ is sub-peelable.
More precisely, the mutable quiver of every seed is a full
induced subquiver of a peelable quiver.
\end{theorem}

\subsection{Sub-peelable matrices}

Recall sub-peelable matrices were defined in Definition~\ref{def:subpeelable}.
\begin{lemma}\label{lem:class}
The class of sub-peelable matrices is closed under mutation, taking principal
submatrices, direct sums, adjoining a source or sink, and
$B\mapsto B^\vee = -B^\T$. Every sub-peelable matrix has a reddening sequence.
\end{lemma}
\begin{proof}
The assertion for mutation, direct sum, taking principal submatrices is obvious.

Assume \(B\) is sub-peelable,
and let $B'$ be obtained from $B$ by adjoining a source or sink,
and fix a positive right symmetrizer $D'$ of $B'$.
For each connected component $B_a$ of $B$, choose a peelable
matrix $R_a$ containing $B_a$ as a principal submatrix.
Since $B_a$ is connected, we can rescale a positive right
symmetrizer of $R_a$ to agree with $D'$ on $B_a$.
Adjoin the same source or sink to $\bigoplus_a R_a$,
keeping its entries with the vertices of $B$ unchanged
and setting its entries with all extra vertices to zero.
This gives a peelable matrix containing $B'$ as a principal submatrix.

Duality preserves sources and sinks and satisfies
$
 \mu_k(B^\vee)=-(\mu_kB)^{\mathsf T}$ and $ (B_{I,I})^\vee=(B^\vee)_{I,I}$.
It therefore preserves peelability and sub-peelability.

The empty matrix has the empty reddening sequence.
Existence of a reddening sequence is preserved by mutation
\cite[Corollary~3.8]{LP} and adjoining a source or sink
\cite[Theorem~4.5 and Remark~4.6]{CL}, so every peelable matrix has a reddening sequence.
Existence  passes to principal submatrices
\cite[Lemma~4.9]{CL}, so every sub-peelable matrix has a
reddening sequence.
\end{proof}
\begin{example} \label{nonexample}
An ideal triangulation of a torus with one puncture
gives the Markov quiver, with exchange matrix
\[
 B=\begin{pmatrix}
 0&2&-2\\-2&0&2\\2&-2&0
 \end{pmatrix}.
\]
This quiver has no reddening sequence  \cite[Remark~3.4]{Ladkani}, so  is not sub-peelable by Lemma~\ref{lem:class}.
\end{example}

Every sink-recurrent matrix is peelable. At each stage of the
sink-recurrent recursion, remove the chosen source or sink and
continue with the remaining matrix.
In the skew-symmetric case, Louise and Banff quivers are peelable
by \cite[Corollary~3.6]{EJ}. However, we are not aware of an example that is peelable but not sink-recurrent. 

The following example is a sub-peelable quiver but it is not sink recurrent.

\begin{example}\label{ex:strict}
Let $W$ be the quiver with matrix
\begin{equation}\label{eq:W}
 B_W=\begin{pmatrix}
 0&2&-1&-1\\-2&0&1&1\\1&-1&0&1\\1&-1&-1&0
 \end{pmatrix}.
\end{equation}
This is the quiver corresponding to an ideal triangulation of a torus with a disk removed, with one boundary marked point.
For $k=1,2,3,4$, respectively, the permutations
\[
 (2,1,3,4),\quad(2,1,3,4),\quad(2,4,1,3),\quad(3,1,4,2)
\]
satisfy $(\mu_k B_W)_{p(i),p(j)}=(B_W)_{ij}$.
Thus every quiver mutation equivalent to $W$ is isomorphic to
$W$ and has no source or sink. In particular, $W$ is not
sink recurrent.

Adjoin a vertex with arrows $4\to5\to3$, obtaining $R$.
Mutations at $(5,3,1,2,5)$ give
\[
 B_{R'}=\begin{pmatrix}
 0&1&-2&1&0\\-1&0&1&0&1\\2&-1&0&-1&0\\
 -1&0&1&0&1\\0&-1&0&-1&0
 \end{pmatrix}.
\]
The vertex $5$ is a sink with neighbors $2,4$. The quivers are shown in Figure~\ref{fig:quiver_nonSR}.
\begin{figure}
\centering
\begin{tikzpicture}[x=.9cm,y=.9cm,>=stealth,line width=.65pt,
  vertex/.style={circle,fill=white,inner sep=2pt,font=\small}]
\begin{scope}
\node[vertex] (w1) at (-1.2,0) {$1$};
\node[vertex] (w2) at (1.8,0) {$2$};
\node[vertex] (w3) at (0,2.8) {$3$};
\node[vertex] (w4) at (.2,.933) {$4$};
\draw[->,bend left=8] (w1) to (w2);
\draw[->,bend right=8] (w1) to (w2);
\draw[->] (w3) -- (w1);
\draw[->] (w4) -- (w1);
\draw[->] (w2) -- (w3);
\draw[->] (w2) -- (w4);
\draw[->] (w3) -- (w4);
\node at (0,-.65) {$W$};
\end{scope}
\begin{scope}[xshift=4.1cm]
\node[vertex] (r1) at (-1.2,0) {$1$};
\node[vertex] (r2) at (1.8,0) {$2$};
\node[vertex] (r3) at (0,2.8) {$3$};
\node[vertex] (r4) at (.2,.933) {$4$};
\node[vertex] (r5) at (2.6,1.7) {$5$};
\draw[->,bend left=8] (r1) to (r2);
\draw[->,bend right=8] (r1) to (r2);
\draw[->] (r3) -- (r1);
\draw[->] (r4) -- (r1);
\draw[->] (r2) -- (r3);
\draw[->] (r2) -- (r4);
\draw[->] (r3) -- (r4);
\draw[->] (r4) -- (r5);
\draw[->] (r5) -- (r3);
\node at (0,-.65) {$R$};
\end{scope}
\begin{scope}[xshift=8.2cm]
\node[vertex] (d1) at (-1.2,0) {$1$};
\node[vertex] (d2) at (1.8,0) {$2$};
\node[vertex] (d3) at (0,2.8) {$3$};
\node[vertex] (d4) at (.2,.933) {$4$};
\node[vertex] (d5) at (2.6,1.7) {$5$};
\draw[->] (d1) -- (d2);
\draw[->,bend left=10] (d3) to (d1);
\draw[->,bend right=10] (d3) to (d1);
\draw[->] (d1) -- (d4);
\draw[->] (d2) -- (d3);
\draw[->] (d2) -- (d5);
\draw[->] (d4) -- (d3);
\draw[->] (d4) -- (d5);
\node at (0,-.65) {$R'$};
\end{scope}
\end{tikzpicture}
\caption{The quivers $W$, $R$, and $R'$.
An entry $B_{ij}>0$ gives $B_{ij}$ arrows $i\to j$.}
\label{fig:quiver_nonSR}
\end{figure}
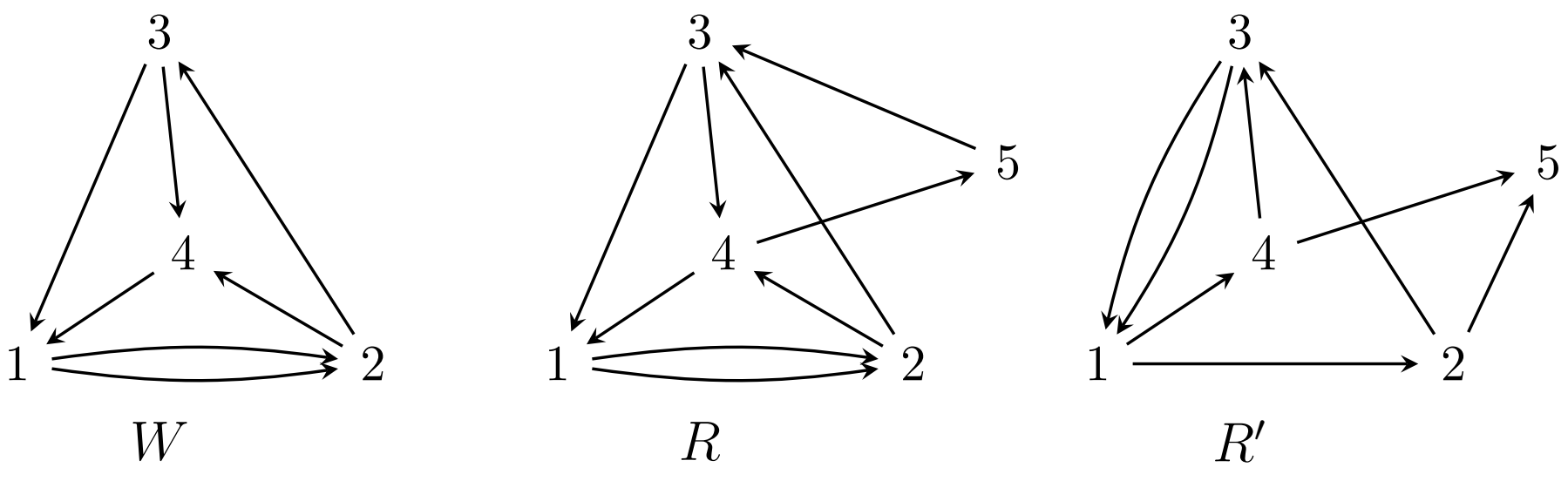

After deleting $5$, mutations at $2,4$ give an acyclic quiver. 
So \(R'\), hence $R$, is peelable.
Since $W=R|_{\{1,2,3,4\}}$, $W$ is sub-peelable.
\end{example}

\subsection{Quivers for Fock--Goncharov spaces}
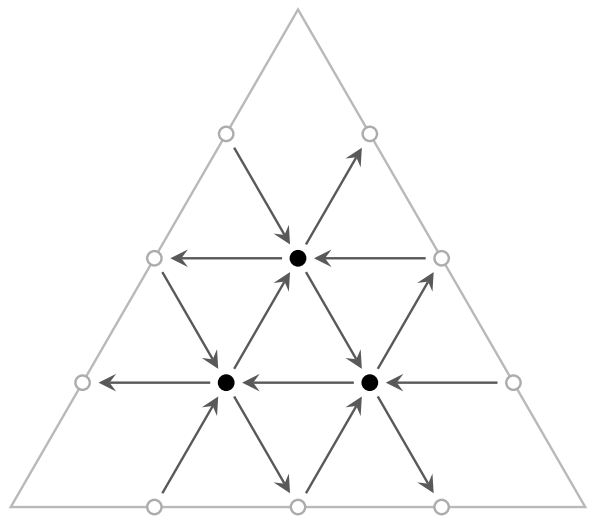
\begin{figure}
\centering
\begin{tikzpicture}[x=0.95cm,y=0.95cm,line width=.45pt,font=\small]
\draw[gray!55] (0,0)--(4,0)--(2,3.4641)--cycle;
\draw[->,>=stealth,black!65,shorten >=3pt,shorten <=3pt] (3.0000,1.7321) -- (2.0000,1.7321);
\draw[->,>=stealth,black!65,shorten >=3pt,shorten <=3pt] (3.5000,0.8660) -- (2.5000,0.8660);
\draw[->,>=stealth,black!65,shorten >=3pt,shorten <=3pt] (1.5000,2.5981) -- (2.0000,1.7321);
\draw[->,>=stealth,black!65,shorten >=3pt,shorten <=3pt] (2.0000,1.7321) -- (2.5000,2.5981);
\draw[->,>=stealth,black!65,shorten >=3pt,shorten <=3pt] (2.0000,1.7321) -- (2.5000,0.8660);
\draw[->,>=stealth,black!65,shorten >=3pt,shorten <=3pt] (2.0000,1.7321) -- (1.0000,1.7321);
\draw[->,>=stealth,black!65,shorten >=3pt,shorten <=3pt] (2.5000,0.8660) -- (3.0000,1.7321);
\draw[->,>=stealth,black!65,shorten >=3pt,shorten <=3pt] (2.5000,0.8660) -- (3.0000,0.0000);
\draw[->,>=stealth,black!65,shorten >=3pt,shorten <=3pt] (2.5000,0.8660) -- (1.5000,0.8660);
\draw[->,>=stealth,black!65,shorten >=3pt,shorten <=3pt] (1.0000,1.7321) -- (1.5000,0.8660);
\draw[->,>=stealth,black!65,shorten >=3pt,shorten <=3pt] (1.5000,0.8660) -- (2.0000,1.7321);
\draw[->,>=stealth,black!65,shorten >=3pt,shorten <=3pt] (1.5000,0.8660) -- (2.0000,0.0000);
\draw[->,>=stealth,black!65,shorten >=3pt,shorten <=3pt] (1.5000,0.8660) -- (0.5000,0.8660);
\draw[->,>=stealth,black!65,shorten >=3pt,shorten <=3pt] (2.0000,0.0000) -- (2.5000,0.8660);
\draw[->,>=stealth,black!65,shorten >=3pt,shorten <=3pt] (1.0000,0.0000) -- (1.5000,0.8660);
\draw[draw=gray!65,fill=white] (2.5000,2.5981) circle (0.05);
\draw[draw=gray!65,fill=white] (3.0000,1.7321) circle (0.05);
\draw[draw=gray!65,fill=white] (3.5000,0.8660) circle (0.05);
\draw[draw=gray!65,fill=white] (1.5000,2.5981) circle (0.05);
\draw[draw=black,fill=black] (2.0000,1.7321) circle (0.05);
\draw[draw=black,fill=black] (2.5000,0.8660) circle (0.05);
\draw[draw=gray!65,fill=white] (3.0000,0.0000) circle (0.05);
\draw[draw=gray!65,fill=white] (1.0000,1.7321) circle (0.05);
\draw[draw=black,fill=black] (1.5000,0.8660) circle (0.05);
\draw[draw=gray!65,fill=white] (2.0000,0.0000) circle (0.05);
\draw[draw=gray!65,fill=white] (0.5000,0.8660) circle (0.05);
\draw[draw=gray!65,fill=white] (1.0000,0.0000) circle (0.05);
\end{tikzpicture}
\caption{The type A triangle quiver for $n=4$. Vertices on $\partial S$ are frozen.}
\label{fig:quiver-triangle}
\end{figure}

An ideal triangulation gives the seed of
\cite[Definition~1.11 and Section~9]{FG}. In each oriented triangle,
take the lattice points
\[
 (a,b,c)\in\Z_{\ge0}^3,\qquad a+b+c=n.
\]
Each point gives a quiver vertex except for the three corners
$(n,0,0)$, $(0,n,0)$ and $(0,0,n)$. Add the oriented cycle
\begin{equation}\label{eq:small-triangle}
 (a+1,b+1,c)\longrightarrow(a,b+1,c+1)
 \longrightarrow(a+1,b,c+1)\longrightarrow(a+1,b+1,c)
\end{equation}
for each $a,b,c\ge0$ with $a+b+c=n-2$, as in Figure~\ref{fig:quiver-triangle}.
The three edges are $a=0$, $b=0$, and $c=0$.
Glue the lattice vertices along common edges and cancel opposite
arrows. Vertices on $\partial S$ are frozen. We omit them
when discussing the mutable quiver.

Write $Q_{S,T}$ for the mutable quiver of a triangulation $T$ of $S$,
or $Q_S$ when the triangulation is understood.

\begin{lemma}\label{lem:disks}
For a disk with $k\ge3$ marked boundary points, the
Fock--Goncharov quiver for $\mathrm{SL}_n$ is peelable.
\end{lemma}
\begin{proof}
Let $\Delta$ be the positive braid lift of the longest Weyl group element,
and let $H\subset\mathrm{SL}_n$ be the diagonal Cartan subgroup.
Then (see, e.g.,~\cite[\S 6.2.2]{IOS})
$
 \Ac^\times_{\mathrm{SL}_n,S}
 \cong X_{\mathrm{SL}_n}(\Delta^{k-1})\times H^{k-1}
$
where $X_{\mathrm{SL}_n}(\Delta^{k-1})$ is the braid variety.
Its mutable quiver is sink-recurrent by \cite{CGGLSS, GLSB, COMP},
hence peelable.
\end{proof}

\subsection{Adding boundary marked points}
\label{subsec:ear}
\begin{figure}
\centering
\begin{tikzpicture}[x=1cm,y=1cm,line width=.55pt,font=\small,>=stealth]
\fill[gray!7] (-1.25,0)--(1.25,0)--(0,2)--cycle;
\draw[gray!65] (-1.25,0)--(0,2)--(1.25,0);
\draw[line width=1pt] (-1.25,0)--(1.25,0);
\coordinate (u) at (-.625,1); \coordinate (v) at (.625,1);
\coordinate (f) at (0,0);
\draw[->,shorten >=4pt,shorten <=4pt] (u)--(v) (v)--(f) (f)--(u);
\fill (u) circle (.055); \fill (v) circle (.055);
\draw[fill=white] (f) circle (.065);
\node[left] at (u) {$u$}; \node[right] at (v) {$v$};
\node[below] at (f) {$f$ frozen};
\node[font=\scriptsize] at (0,2.4) {before adding a mark};
\draw[->] (1.7,.8)--(2.8,.8);
\node[above,font=\scriptsize] at (2.25,.8) {one ear};
\begin{scope}[shift={(4.5,0)}]
\fill[gray!7] (-1.25,0)--(1.25,0)--(0,2)--cycle;
\fill[blue!10] (-1.25,0)--(1.25,0)--(0,-1.25)--cycle;
\draw[gray!65] (-1.25,0)--(0,2)--(1.25,0);
\draw[line width=1pt] (-1.25,0)--(0,-1.25)--(1.25,0);
\draw[blue!65!black,dashed] (-1.25,0)--(1.25,0);
\coordinate (u2) at (-.625,1); \coordinate (v2) at (.625,1);
\coordinate (f2) at (0,0);
\coordinate (g2) at (-.625,-.625); \coordinate (h2) at (.625,-.625);
\draw[->,shorten >=4pt,shorten <=4pt] (u2)--(v2) (v2)--(f2) (f2)--(u2);
\draw[->,blue!65!black,shorten >=4pt,shorten <=4pt]
  (f2)--(h2) (h2)--(g2) (g2)--(f2);
\fill (u2) circle (.055); \fill (v2) circle (.055);
\fill[blue!65!black] (f2) circle (.065);
\draw[fill=white] (g2) circle (.065); \draw[fill=white] (h2) circle (.065);
\fill[blue!65!black] (0,-1.25) circle (.05);
\node[left] at (u2) {$u$}; \node[right] at (v2) {$v$};
\draw[blue!65!black,thin] (.12,.03)--(1.4,.25);
\node[right,text=blue!65!black] at (1.4,.25) {$f$ mutable};
\node[below,font=\scriptsize,text=blue!65!black] at (0,-1.25) {new boundary mark};
\node[font=\scriptsize] at (0,2.4) {old mutable arrow $u\to v$ unchanged};
\end{scope}
\path (0,-1.8);
\end{tikzpicture}
\caption{Attaching an ear for $\mathrm{SL}_2$.}
\label{fig:ear}
\end{figure}
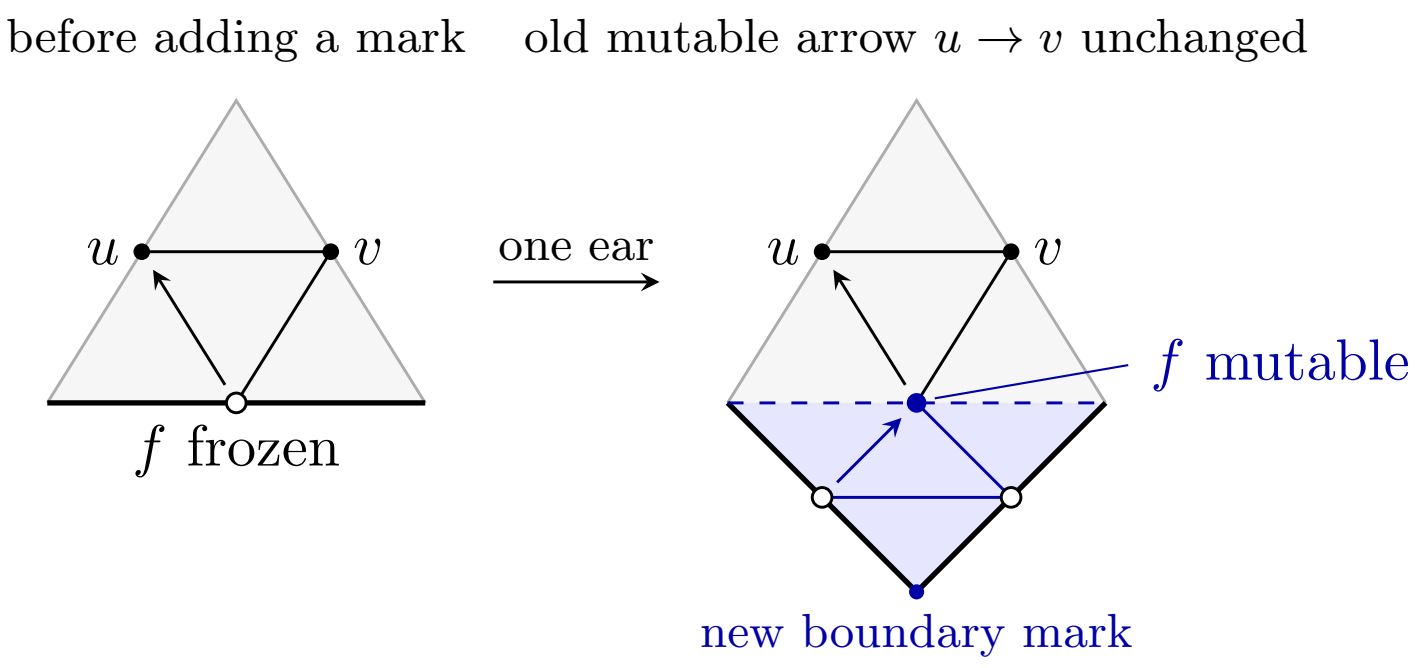

Attaching an ``ear'' along a boundary edge adds one boundary marked point,
as in Figure~\ref{fig:ear}.
The old mutable quiver is a full induced subquiver of the new one,
so it is sub-peelable if the new quiver is.

\subsection{Cutting the surface along a strip} 
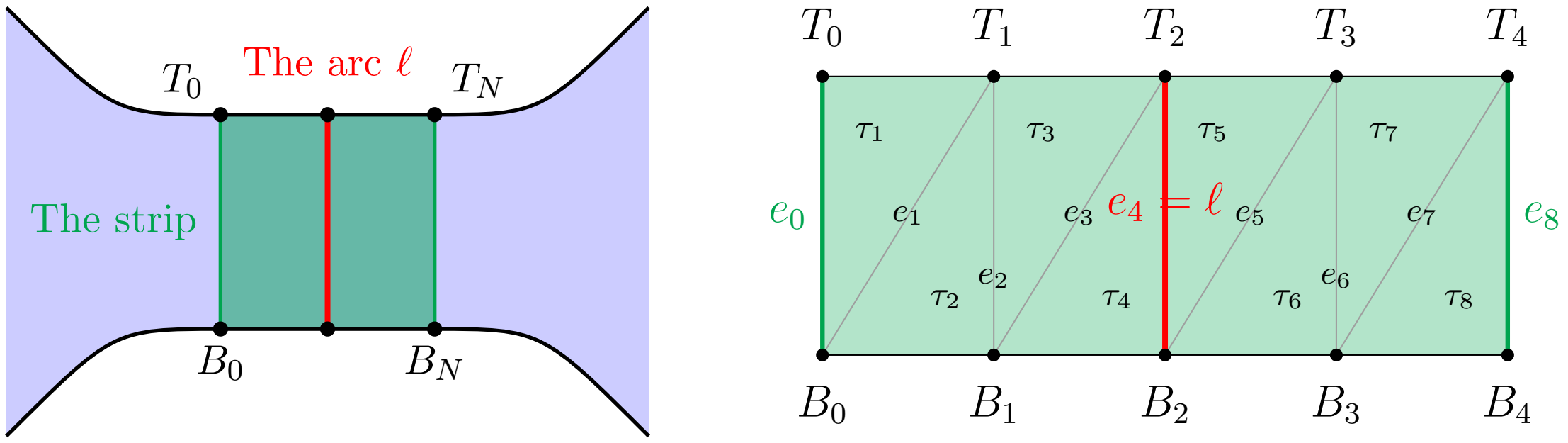
\begin{figure}
    \centering
    \begin{tikzpicture}[line width=1pt, scale = 1]
        \fill[blue, fill opacity=0.2]
            (-2, 1) .. controls (-1, 0) .. (0, 0)
            -- (2, 0) .. controls (3, 0) .. (4, 1)
            -- (4, -3) .. controls (3, -2) .. (2, -2)
            -- (0, -2) .. controls (-1, -2) .. (-2, -3)
            -- cycle;

        \fill[Green, fill opacity=0.5]
            (0, 0) rectangle (2, -2);

        \draw
            (-2, 1) .. controls (-1, 0) .. (0, 0)
            -- (2, 0) .. controls (3, 0) .. (4, 1);
        \draw
            (-2, -3) .. controls (-1, -2) .. (0, -2)
            -- (2, -2) .. controls (3, -2) .. (4, -3);

        \draw[Green] (0, 0) -- (0, -2);
        \draw[Green] (2, 0) -- (2, -2);

        \draw[red, line width=1.5pt]
            (1, 0) -- (1, -2);

        \foreach \x in {0,1,2} {
            \fill[black] (\x, 0) circle (2pt);
            \fill[black] (\x, -2) circle (2pt);
        }
        \node[above left] at (0, 0) {\small $T_0$};
        \node[above right] at (2, 0) {\small $T_{N}$};
        \node[below] at (2, -2) {\small $B_{N}$};
        \node[below] at (0, -2) {\small $B_0$};
        \node[red]   at (1, .5) {\small{The arc \(\ell\)}};
        \node[Green] at (-1, -1) {\small{The strip}};
    \end{tikzpicture}
    \qquad
    \begin{tikzpicture}[x=1.6cm,y=1.3cm]
        \fill[Green,fill opacity=.3] (0,0) rectangle (4,2);
        \draw (0,2)--(4,2) (0,0)--(4,0);
        \foreach \i in {1,3} {
            \draw[gray!75] (\i,0)--(\i,2);
        }
        \foreach \i in {0,...,3} {
            \draw[gray!75] (\i,0)--(\i+1,2);
            \pgfmathtruncatemacro{\j}{2*\i+1}
            \node[ inner sep=1pt,font=\scriptsize]
                at (\i+.5,1) {$e_{\j}$};
            \node[inner sep=1pt,font=\scriptsize]
                at (\i+.28,1.6) {$\tau_{\j}$};
            \pgfmathtruncatemacro{\k}{2*\i+2}
            \node[inner sep=1pt,font=\scriptsize]
                at (\i+.72,.4) {$\tau_{\k}$};
        }
        \draw[Green,line width=1.2pt] (0,0)--(0,2) (4,0)--(4,2);
        \draw[red,line width=1.5pt] (2,0)--(2,2);
        \foreach \i in {0,...,4} {
            \fill (\i,2) circle (1.7pt);
            \fill (\i,0) circle (1.7pt);
            \node[above=4pt] at (\i,2) {$T_{\i}$};
            \node[below=4pt] at (\i,0) {$B_{\i}$};
        }
        \node[left,text=Green] at (0,1) {$e_0$};
        \node[right,text=Green] at (4,1) {$e_8$};
        \node[ inner sep=1pt,font=\scriptsize] at (1,.55) {$e_2$};
        \node[ inner sep=1pt,text=red] at (2,1.1) {$e_4 = \ell$};
        \node[ inner sep=1pt,font=\scriptsize] at (3,.55) {$e_6$};
    \end{tikzpicture}
    \caption{The arc $\ell$ (red), its strip neighbourhood (green), and a triangulation with $2N=8$ triangles.}
    \label{fig:arc-strip}
    \label{fig:strip}
\end{figure}

Let $g$ and $b$ be the genus and the number of boundary components of $S$,
and put $h(S)=2g+b-1$.
Suppose $h(S)>0$. By Section~\ref{subsec:ear}, we may assume that
there are at least two boundary marked points.
Choose an arc $\ell$ joining distinct marked points such that
$S\setminus\ell$ is connected and $h(S\setminus\ell)=h(S)-1$.
Here $S\setminus \ell$ denotes the compact surface obtained by cutting
along $\ell$, with both copies of each endpoint marked.
Take a strip neighbourhood of $\ell$ as in Figure~\ref{fig:arc-strip}.

By Section~\ref{subsec:ear}, we can add marked points on the boundary sides of the strip. Choose an integer
$ N\ge\max\{\lfloor n/2\rfloor,2\}$ and label the points on its
two boundary sides
\[
  T_0,\ldots,T_{N},\qquad B_0,\ldots,B_{N}
\]
as in Figure~\ref{fig:strip}. 
Triangulate the strip by the arcs
\begin{equation}\label{eq:strip-triangulation}
 e_{2i}=T_iB_i\quad(0\le i\le N),\qquad
 e_{2i+1}=T_{i+1}B_i\quad(0\le i<N)
\end{equation}
and let \(\ell= e_N\).
Denote the triangles by
$\tau_1,\ldots,\tau_{2N}$, with $\tau_i$ bounded by
$e_{i-1}$, $e_i$, and a segment of $\partial S$.
Extend the triangulation \eqref{eq:strip-triangulation} to $S$ and consider the quiver \(Q_S\). 

\begin{definition}
A \emph{deletion move} at a vertex $s$ of a quiver \(Q\) consists of applying
a sequence of mutations avoiding $s$ that makes $s$ a sink
or source, deleting $s$, and undoing the mutations in reverse order.
The result is the full induced subquiver $Q\setminus\{s\}$. By Lemma \ref{lem:class}, $Q$ is sub-peelable if and only if
$Q\setminus\{s\}$ is sub-peelable.
\end{definition}

A deletion move at $s$ in $Q$, with mutation sequence $W$, gives
\begin{equation}\label{eq:delete-undo}
 \mu_W^{-1}\bigl((\mu_W Q)\setminus\{s\}\bigr)
       =Q\setminus\{s\}.
\end{equation}

\begin{proposition}\label{prop:strip}
Let $S$, $\ell$, and the strip be as above.
There is a sequence of deletion moves taking $Q_S$ to
the full induced subquiver $Q_S\setminus Z$, where $Z$
contains all vertices on $\ell$.
All deletion moves take place inside the strip, and
no vertex on $e_0$ or $e_{2N}$ is mutated or deleted.
\end{proposition}

\begin{proof}[Proof of Theorem~\ref{thm:surface-quiver}]
We argue by induction on $h(S)$, for all boundary marked sets
satisfying \eqref{eq:surface-assumption}.
If $h(S)=0$, then $S$ is a disk and Lemma~\ref{lem:disks} applies.
Suppose $h(S)>0$ and assume the assertion for smaller $h$.
Take an arc \(\ell\) and a strip neighborhood as above.

Cut \(S\) along $\ell$ and replace it by two boundary edges, denoted
$\ell$ and $\ell'$. The following picture shows the cut strip for $N=4$.
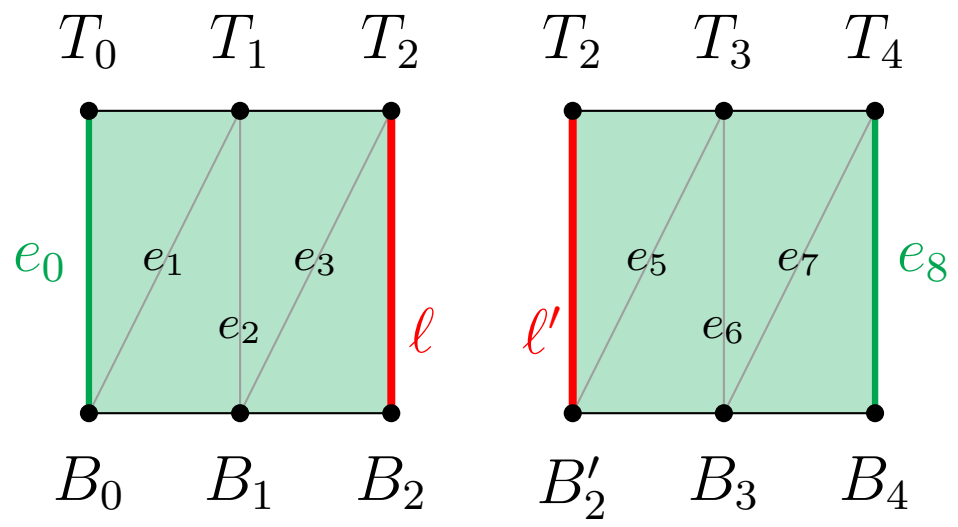
\begin{center}
\begin{tikzpicture}[scale=1]
    \begin{scope}
        \fill[Green,fill opacity=.3] (0,0) rectangle (2,2);
        \draw (0,2)--(2,2) (0,0)--(2,0);
        \draw[gray!75] (1,0)--(1,2);
        \foreach \i in {0,1} {
            \draw[gray!75] (\i,0)--(\i+1,2);
            \pgfmathtruncatemacro{\j}{2*\i+1}
            \node[inner sep=1pt,font=\scriptsize]
                at (\i+.5,1) {$e_{\j}$};
        }
        \draw[Green,line width=1.2pt] (0,0)--(0,2);
        \draw[red,line width=1.5pt] (2,0)--(2,2);
        \foreach \i in {0,1,2} {
            \fill (\i,2) circle (1.7pt);
            \fill (\i,0) circle (1.7pt);
            \node[above=4pt] at (\i,2) {$T_{\i}$};
            \node[below=4pt] at (\i,0) {$B_{\i}$};
        }
        \node[left,text=Green] at (0,1) {$e_0$};
        \node[ inner sep=1pt,font=\scriptsize]
            at (1,.55) {$e_2$};
        \node[ inner sep=1pt,text=red]
            at (2.2,.55) {$\ell$};
    \end{scope}

    \begin{scope}[xshift=1.2cm]
        \fill[Green,fill opacity=.3] (2,0) rectangle (4,2);
        \draw (2,2)--(4,2) (2,0)--(4,0);
        \draw[gray!75] (3,0)--(3,2);
        \foreach \i in {2,3} {
            \draw[gray!75] (\i,0)--(\i+1,2);
            \pgfmathtruncatemacro{\j}{2*\i+1}
            \node[ inner sep=1pt,font=\scriptsize]
                at (\i+.5,1) {$e_{\j}$};
        }
        \draw[Green,line width=1.2pt] (4,0)--(4,2);
        \draw[red,line width=1.5pt] (2,0)--(2,2);
        \foreach \i in {2,3,4} {
            \fill (\i,2) circle (1.7pt);
            \fill (\i,0) circle (1.7pt);
            \node[above=4pt] at (\i,2) {$T_{\i}$};
        }
        \node[below=4pt] at (2,0) {$B_2'$};
        \foreach \i in {3,4} {
            \node[below=4pt] at (\i,0) {$B_{\i}$};
        }
        \node[right,text=Green] at (4,1) {$e_8$};
        \node[ inner sep=1pt,text=red]
            at (1.8,.55) {$\ell'$};
        \node[ inner sep=1pt,font=\scriptsize]
            at (3,.55) {$e_6$};
    \end{scope}
\end{tikzpicture}
\end{center}
By induction, $Q_{S\setminus\ell}$ is sub-peelable.
Proposition~\ref{prop:strip} gives deletion moves from $Q_S$ to $Q_S\setminus Z$.
Since $Z$ contains every vertex on $\ell$, the latter is a full induced
subquiver of $Q_{S\setminus\ell}$.
It is therefore sub-peelable, and so is $Q_S$.
\end{proof}

We prove Proposition~\ref{prop:strip} in the next subsection. Proposition~\ref{prop:strip} is also formalized in \href{https://github.com/Hu-Mingyuan/peeling-strip}{Lean}.

\subsection{Mutations inside the strip}
\label{sec:strip-local}

\begin{figure}
\centering
\begin{tikzpicture}[x=1.65cm,y=1.2cm,font=\small,>=stealth,
    every node/.style={fill=white,inner sep=2pt}]
\draw[gray!55] (0,0)--(3,6)--(6,0);
\draw[black,solid] (0,0)--(6,0);
\node[text=gray,above] at (3,6) {$(0,0,6)$};
\node[text=gray,below left] at (0,0) {$(6,0,0)$};
\node[text=gray,below right] at (6,0) {$(0,6,0)$};
\node[below] at (3,0) {$\{c=0\}\subset\partial S$};
\node[{text=blue,label={[text=blue]left:$D_5$}}] (D5) at (2.5,5) {$\times$};
\node[text=blue,label={[text=blue]right:$R_5$}
] (R5) at (3.5,5) {$\times$};
\node[text=blue,label={[text=blue]left:$D_4$}] (D4) at (2,4) {$\times$};        
\node (v114) at (3,4) {$(1,1,4)$};
\node[text=blue,label={[text=blue]right:$R_4$}] (R4) at (4,4) {$\times$};
\node[text=blue] (D3) at (1.5,3) {$D_3$};
\node (v213) at (2.5,3) {$(2,1,3)$};
\node (v123) at (3.5,3) {$(1,2,3)$};
\node[text=blue] (R3) at (4.5,3) {$R_3$};
\node[text=blue,label={[text=blue]left:$D_2$}] (D2) at (1,2) {$\times$};
\node (v312) at (2,2) {$(3,1,2)$};
\node (v222) at (3,2) {$(2,2,2)$};
\node (v132) at (4,2) {$(1,3,2)$};
\node[text=blue,label={[text=blue]right:$R_2$}] (R2) at (5,2) {$\times$};
\node[text=blue,label={[text=blue]left:$D_1$}] (D1) at (0.5,1) {$\times$};
\node[text=red,label={[text=red]below:$(4,1,1)$}] (v411) at (1.5,1) {$\times$};
\node (v321) at (2.5,1) {$(3,2,1)$};
\node (v231) at (3.5,1) {$(2,3,1)$};
\node[text=red,label={[text=red]below:$(1,4,1)$}] (v141) at (4.5,1) {$\times$};
\node[text=blue,label={[text=blue]right:$R_1$}] (R1) at (5.5,1) {$\times$};
\draw[->] (R3)--(v123);
\draw[->] (v114)--(v123);
\draw[->] (v123)--(v132);
\draw[->] (v123)--(v213);
\draw[->] (v132)--(R3);
\draw[->] (v132)--(v222);
\draw[->] (v213)--(v114);
\draw[->] (v213)--(v222);
\draw[->] (v213)--(D3);
\draw[->] (v222)--(v123);
\draw[->] (v222)--(v231);
\draw[->] (v222)--(v312);
\draw[->] (v231)--(v132);
\draw[->] (v231)--(v321);
\draw[->] (D3)--(v312);
\draw[->] (v312)--(v213);
\draw[->] (v312)--(v321);
\draw[->] (v321)--(v222);
\end{tikzpicture}
\caption{The hypotheses of Lemma~\ref{lem:triangle-shell} for $n=6$ and $p=2$.
We have \red{$C_2=\{(4,1,1),(1,4,1)\}$}
and the blue vertices $D_k,R_k$ have been deleted to depth $2$.}
\label{fig:triangle-notation}
\end{figure}
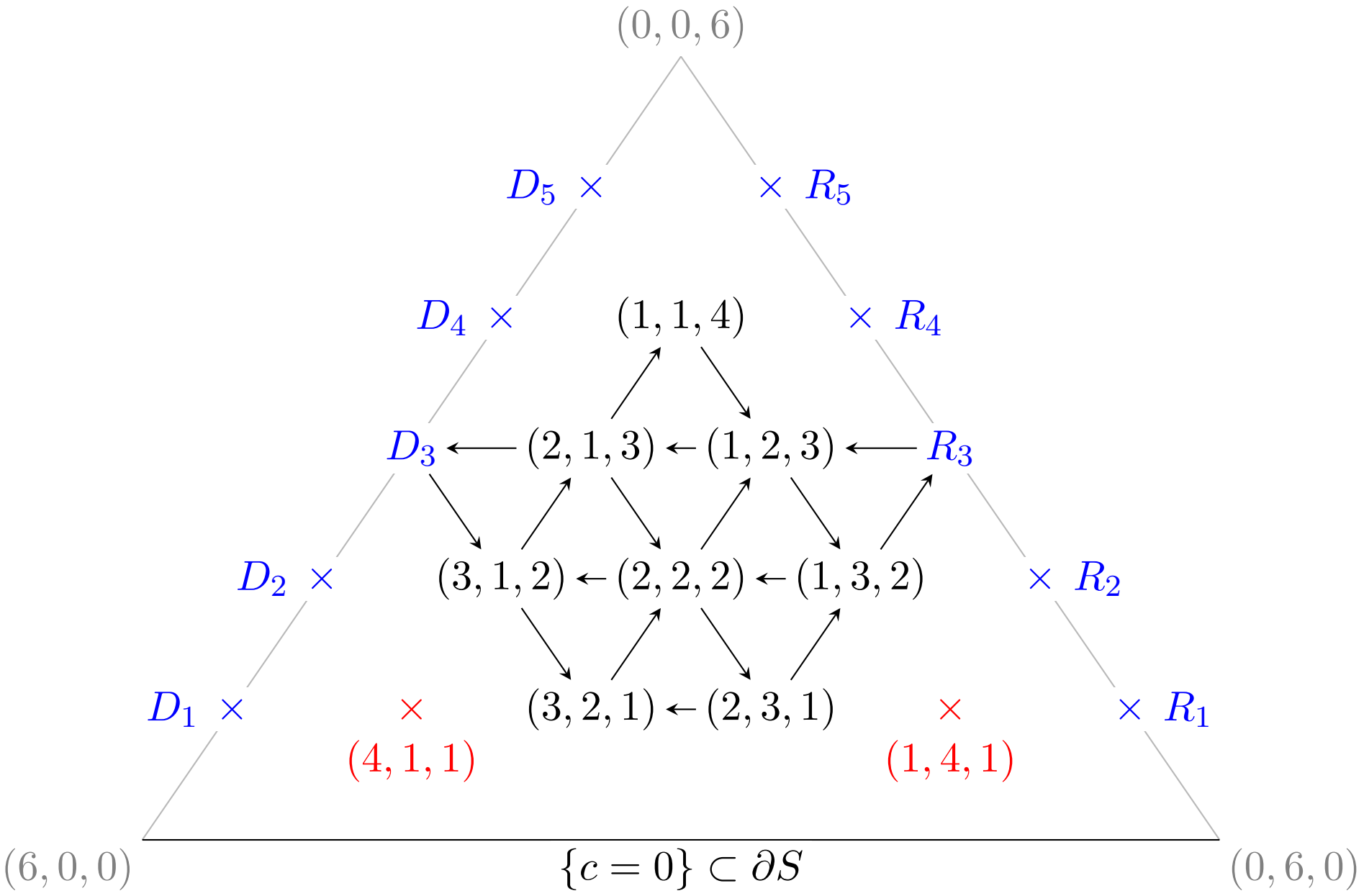

For an edge with vertices numbered $1,\ldots,n-1$, deleting to
depth $p$ means deleting precisely the vertices with
$k\le p$ or $k\ge n-p$. The definition does not depend on the
orientation used to number the edge.
In a triangle whose side $c=0$ lies on $\partial S$, put
\begin{equation}\label{eq:triangle-corners}
 C_p=\{(a,b,c)\in\Z_{>0}^3:a+b+c=n,
              \ a+c\le p\ \text{or}\ b+c\le p\}.
\end{equation}
In particular $C_0=C_1=\varnothing$.
The other two sides have vertices $D_k=(n-k,0,k)$ and
$R_k=(0,n-k,k)$ for $1\le k\le n-1$.
Figure~\ref{fig:triangle-notation} illustrates the hypotheses of Lemma~\ref{lem:triangle-shell} for $n=6$ and $p=2$.

For a quiver $Q$ with matrix $B$, write $\rho_Q(v)=B_{v,*}$,
with $B_{uv}>0$ for an arrow $u\to v$.
We abbreviate $\rho_Q$ to $\rho_{\mathrm{pre}}$ just before the indicated mutation,
and to $\rho_{\mathrm{final}}$ after the mutation sequence, before deletion.
Let $\mathbf e_v$ be the coordinate vector of $v$.
Terms indexed by absent vertices are zero.

\begin{lemma}\label{lem:triangle-shell}
Let $1 \le p \le \frac{n-1}{2}$.
In a triangle with boundary side $\{c=0\}$, suppose that the other
two edges have been deleted to depth $p$ and that $C_p$ has been deleted.
Then $C_{p+1}\setminus C_p$ can be removed by a sequence of deletion moves,
with all mutations and deleted vertices in the interior of the triangle.
\end{lemma}
\begin{proof}
We delete the vertices with $a+c=p+1$
in decreasing order of $a$.
Exchanging $a$ and $b$ gives the deletion moves for $b+c=p+1$.
For $n=3$ and $p=1$, the only interior vertex is isolated,
so no mutation is needed. Assume $n \ge 4$ and put $b_0=n-p-1$.

For $a=p,p-1,\ldots,1$, use the mutation sequence
\begin{equation}\label{eq:triangle-word}
 s_a=(a,b_0,p+1-a),\qquad
 v_1,\ldots,v_{b_0-1},\qquad v_b=(a,b,n-a-b).
\end{equation}
The vertices $v_b$ lie in the interior and satisfy
\[
 a+(n-a-b)=n-b\ge p+2,\qquad
 b+(n-a-b)=n-a\ge n-p\ge p+1.
\]
They lie outside $C_p$ and differ from the vertices $s_{a'}$ already deleted for $a'>a$.

Put
\[
 A_b=(a-1,b+1,n-a-b),\qquad
 B_b=(a+1,b,n-a-b-1).
\]
The incoming arrows at $v_1$ are $A_1\to v_1$ and $B_1\to v_1$,
with two exceptions. If $a=1$ and $p\ge2$, then
$A_1=(0,2,n-2)$ has been deleted from the edge.
If $n=2p+1$ and $a=p$, then $B_1=(p+1,1,p-1)\in C_p$
has been deleted. The third incoming neighbour in the original triangle,
$(a,0,n-a)$, has been deleted from the edge because $a\le p$.
For $2\le b\le b_0-1$,  before mutating $v_b$ we have
\begin{equation}\label{eq:triangle-pivot}
 \rho_{\mathrm{pre}}(v_b)
 =\mathbf e_{v_{b-1}}+\mathbf e_{v_{b+1}}
  -\mathbf e_{A_b}-\mathbf e_{B_b},
 \qquad v_{b_0}=s_a.
\end{equation}
To check \eqref{eq:triangle-pivot}, use the cycles in \eqref{eq:small-triangle}.
The mutations at $v_1,\ldots,v_{b-2}$ leave the arrows at $v_b$ unchanged.
Mutation at $v_{b-1}$ cancels the outgoing arrows from $v_b$
off the line $a=\mathrm{constant}$ and reverses $v_{b-1}\to v_b$.

Before the last mutation,
$
 \rho_{\mathrm{pre}}(s_a)=\mathbf e_{A_{b_0-1}}+\mathbf e_{B_{b_0-1}}
                  -\mathbf e_{v_{b_0-1}}.
$
Indeed, the neighbours
$
 (a-1,b_0+1,p+1-a)$ and $(a,b_0+1,p-a)$
have the sum of their first and third coordinates equal to $p$.
Each lies in $C_p$, on the deleted part of $a=0$, or on $c=0$.
The remaining possible incoming neighbour is $s_{a+1}$, already deleted,
or a vertex on $c=0$ when $a=p$.
Mutating $v_{b_0-1}$ cancels the arrows
$s_a\to A_{b_0-1}$ and $s_a\to B_{b_0-1}$, and reverses
$v_{b_0-1}\to s_a$. Thus
$
 \rho_{\mathrm{final}}(s_a)=\mathbf e_{v_{b_0-1}}.
$
Delete $s_a$ and undo the mutations in \eqref{eq:triangle-word}.

For $b+c=p+1$, the deletion moves use vertices with $b\le p$.
The deleted vertices $s_p,\ldots,s_1$ have $b=b_0\ge p$.
If $b_0=p$ (that is, $n=2p+1$), we skip $(p,p,1)$,
which has already been deleted.
\end{proof}

\begin{figure}
\centering
\begin{tikzpicture}[x=.66cm,y=.66cm,font=\small,
    every node/.style={inner sep=2pt}]

\begin{scope}[xshift=0.0cm]
\draw[gray!65] (0,0)--(7,0);
\draw[blue!60!black] (0,0)--(3.5,6.0622)--(7,0);
\node[below] at (3.5,0) {$\{c=0\}\subset\partial S$};
\node[left] at (1.35,3) {$b=0$};
\node[right] at (5.65,3) {$a=0$};
\draw[gray!65,line width=.7pt] (3.9100,5.1062)--(4.0900,5.2862) (3.9100,5.2862)--(4.0900,5.1062);
\draw[gray!65,line width=.7pt] (4.4100,4.2401)--(4.5900,4.4201) (4.4100,4.4201)--(4.5900,4.2401);
\fill[blue!60!black] (5.0000,3.4641) circle (2pt);
\fill[blue!60!black] (5.5000,2.5981) circle (2pt);
\draw[gray!65,line width=.7pt] (5.9100,1.6421)--(6.0900,1.8221) (5.9100,1.8221)--(6.0900,1.6421);
\draw[gray!65,line width=.7pt] (6.4100,0.7760)--(6.5900,0.9560) (6.4100,0.9560)--(6.5900,0.7760);
\draw[gray!65,line width=.7pt] (2.9100,5.1062)--(3.0900,5.2862) (2.9100,5.2862)--(3.0900,5.1062);
\fill[black] (3.5000,4.3301) circle (2pt);
\fill[black] (4.0000,3.4641) circle (2pt);
\fill[black] (4.5000,2.5981) circle (2pt);
\fill[orange!85!black] (5.0000,1.7321) circle (2pt);
\draw[gray!65,line width=.7pt] (5.4100,0.7760)--(5.5900,0.9560) (5.4100,0.9560)--(5.5900,0.7760);
\draw[gray!65,line width=.7pt] (2.4100,4.2401)--(2.5900,4.4201) (2.4100,4.4201)--(2.5900,4.2401);
\fill[black] (3.0000,3.4641) circle (2pt);
\fill[black] (3.5000,2.5981) circle (2pt);
\fill[black] (4.0000,1.7321) circle (2pt);
\fill[orange!85!black] (4.5000,0.8660) circle (2pt);
\fill[blue!60!black] (2.0000,3.4641) circle (2pt);
\fill[black] (2.5000,2.5981) circle (2pt);
\fill[black] (3.0000,1.7321) circle (2pt);
\fill[black] (3.5000,0.8660) circle (2pt);
\fill[blue!60!black] (1.5000,2.5981) circle (2pt);
\fill[orange!85!black] (2.0000,1.7321) circle (2pt);
\fill[orange!85!black] (2.5000,0.8660) circle (2pt);
\draw[gray!65,line width=.7pt] (0.9100,1.6421)--(1.0900,1.8221) (0.9100,1.8221)--(1.0900,1.6421);
\draw[gray!65,line width=.7pt] (1.4100,0.7760)--(1.5900,0.9560) (1.4100,0.9560)--(1.5900,0.7760);
\draw[gray!65,line width=.7pt] (0.4100,0.7760)--(0.5900,0.9560) (0.4100,0.9560)--(0.5900,0.7760);
\fill[blue!70!black] (3.0000,3.4641) circle (2.3pt);
\node[left,text=blue!70!black] at (3.0000,3.4641) {$v_1$};
\draw[->,>=stealth,dashed,blue!70!black,shorten <=4pt,shorten >=4pt] (3.0000,3.4641)--(3.5000,2.5981);
\fill[blue!70!black] (3.5000,2.5981) circle (2.3pt);
\node[left,text=blue!70!black] at (3.5000,2.5981) {$v_2$};
\draw[->,>=stealth,dashed,blue!70!black,shorten <=4pt,shorten >=4pt] (3.5000,2.5981)--(4.0000,1.7321);
\fill[blue!70!black] (4.0000,1.7321) circle (2.3pt);
\node[left,text=blue!70!black] at (4.0000,1.7321) {$v_3$};
\node[right,text=orange!85!black] at (4.5,.866) {$s_2$};
\node[right,text=orange!85!black] at (5,1.732) {$s_1$};
\node at (3.5,-1) {Edges at depth $2$; $C_2$ deleted};
\end{scope}
\begin{scope}[xshift=6.5cm]
\draw[gray!65] (0,0)--(7,0);
\draw[blue!60!black] (0,0)--(3.5,6.0622)--(7,0);
\node[below] at (3.5,0) {$\{c=0\}\subset\partial S$};
\node[left] at (1.35,3) {$b=0$};
\node[right] at (5.65,3) {$a=0$};
\draw[gray!65,line width=.7pt] (3.9100,5.1062)--(4.0900,5.2862) (3.9100,5.2862)--(4.0900,5.1062);
\draw[gray!65,line width=.7pt] (4.4100,4.2401)--(4.5900,4.4201) (4.4100,4.4201)--(4.5900,4.2401);
\fill[blue!60!black] (5.0000,3.4641) circle (2pt);
\fill[blue!60!black] (5.5000,2.5981) circle (2pt);
\draw[gray!65,line width=.7pt] (5.9100,1.6421)--(6.0900,1.8221) (5.9100,1.8221)--(6.0900,1.6421);
\draw[gray!65,line width=.7pt] (6.4100,0.7760)--(6.5900,0.9560) (6.4100,0.9560)--(6.5900,0.7760);
\draw[gray!65,line width=.7pt] (2.9100,5.1062)--(3.0900,5.2862) (2.9100,5.2862)--(3.0900,5.1062);
\fill[black] (3.5000,4.3301) circle (2pt);
\fill[black] (4.0000,3.4641) circle (2pt);
\fill[black] (4.5000,2.5981) circle (2pt);
\draw[orange!85!black,line width=.7pt] (4.9100,1.6421)--(5.0900,1.8221) (4.9100,1.8221)--(5.0900,1.6421);
\draw[gray!65,line width=.7pt] (5.4100,0.7760)--(5.5900,0.9560) (5.4100,0.9560)--(5.5900,0.7760);
\draw[gray!65,line width=.7pt] (2.4100,4.2401)--(2.5900,4.4201) (2.4100,4.4201)--(2.5900,4.2401);
\fill[black] (3.0000,3.4641) circle (2pt);
\fill[black] (3.5000,2.5981) circle (2pt);
\fill[black] (4.0000,1.7321) circle (2pt);
\draw[orange!85!black,line width=.7pt] (4.4100,0.7760)--(4.5900,0.9560) (4.4100,0.9560)--(4.5900,0.7760);
\fill[blue!60!black] (2.0000,3.4641) circle (2pt);
\fill[black] (2.5000,2.5981) circle (2pt);
\fill[black] (3.0000,1.7321) circle (2pt);
\fill[black] (3.5000,0.8660) circle (2pt);
\fill[blue!60!black] (1.5000,2.5981) circle (2pt);
\draw[orange!85!black,line width=.7pt] (1.9100,1.6421)--(2.0900,1.8221) (1.9100,1.8221)--(2.0900,1.6421);
\draw[orange!85!black,line width=.7pt] (2.4100,0.7760)--(2.5900,0.9560) (2.4100,0.9560)--(2.5900,0.7760);
\draw[gray!65,line width=.7pt] (0.9100,1.6421)--(1.0900,1.8221) (0.9100,1.8221)--(1.0900,1.6421);
\draw[gray!65,line width=.7pt] (1.4100,0.7760)--(1.5900,0.9560) (1.4100,0.9560)--(1.5900,0.7760);
\draw[gray!65,line width=.7pt] (0.4100,0.7760)--(0.5900,0.9560) (0.4100,0.9560)--(0.5900,0.7760);
\node at (3.5,-1) {After four deletion moves};
\end{scope}
\begin{scope}[shift={(1,-4.2)},x=1cm,y=1cm]
  \node (A) at (0,.9) {$A_3$};
  \node (B) at (0,-.9) {$B_3$};
  \node[text=blue!70!black] (v) at (1.5,0) {$v_3$};
  \node[text=orange!85!black] (s) at (3,0) {$s_2$};
  \draw[->,>=stealth,line width=.7pt] (A)--(v);
  \draw[->,>=stealth,line width=.7pt] (B)--(v);
  \draw[->,>=stealth,line width=.7pt] (v)--(s);
  \draw[->,>=stealth,line width=.7pt] (s) to[bend right=18] (A);
  \draw[->,>=stealth,line width=.7pt] (s) to[bend left=18] (B);
  \draw[->,>=stealth,line width=.8pt] (3.8,0)--node[above] {$\mu_{v_3}$} (5.3,0);
  \node[text=blue!70!black] (vv) at (6.2,0) {$v_3$};
  \node[text=orange!85!black] (ss) at (7.8,0) {$s_2$};
  \draw[->,>=stealth] (ss)--(vv);
  \node at (7,-.7) {$s_2$ is a source};
\end{scope}
\end{tikzpicture}
\caption{Lemma~\ref{lem:triangle-shell} for $n=7$, $p=2$.
Orange vertices form $C_3\setminus C_2$. Dashed arrows indicate the mutation order.}
\label{fig:shell}
\end{figure}
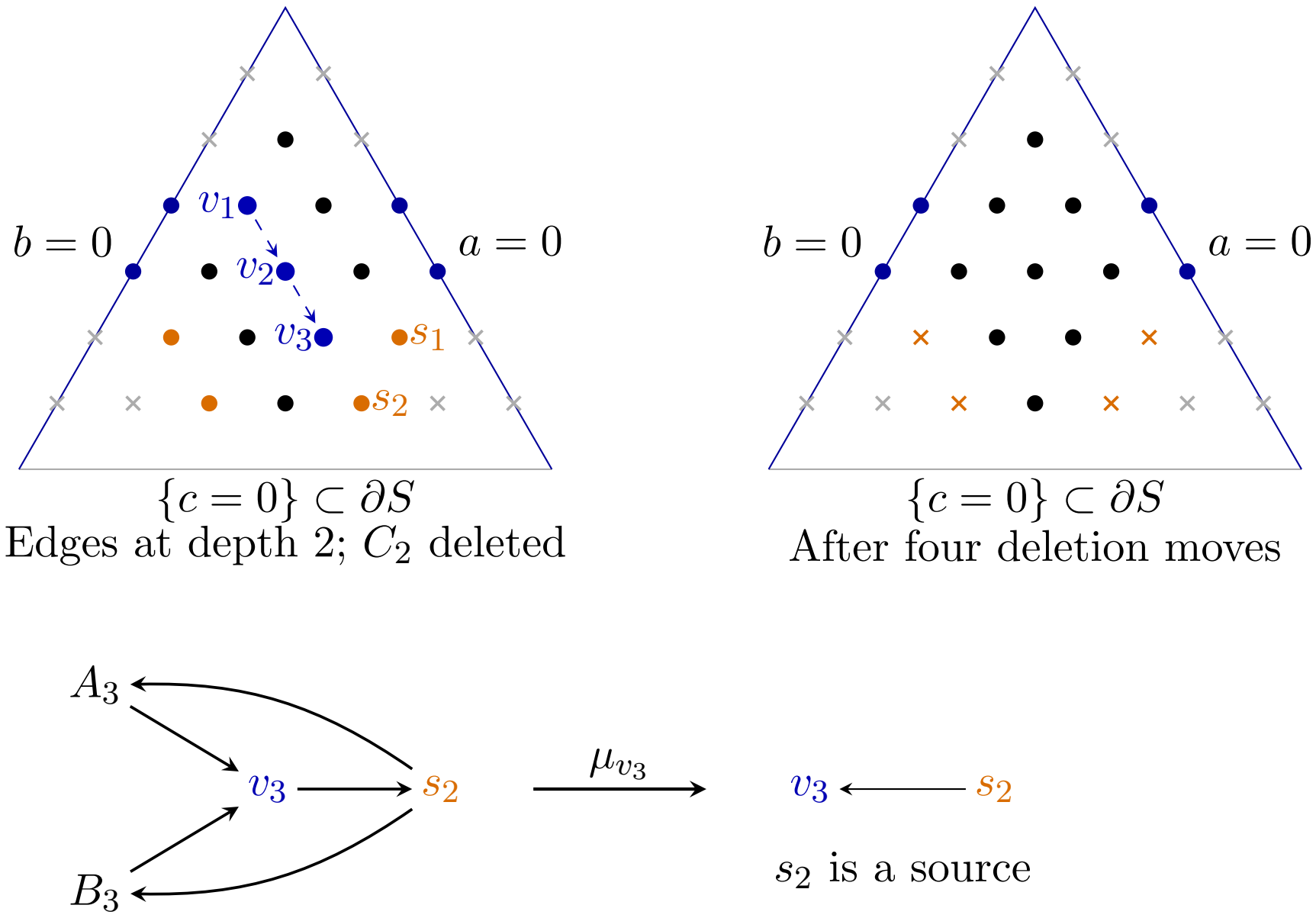

In the initial configuration shown in the upper left panel of Figure~\ref{fig:shell},
both sloping edges have been deleted to depth $2$, and $C_2$ has been deleted.

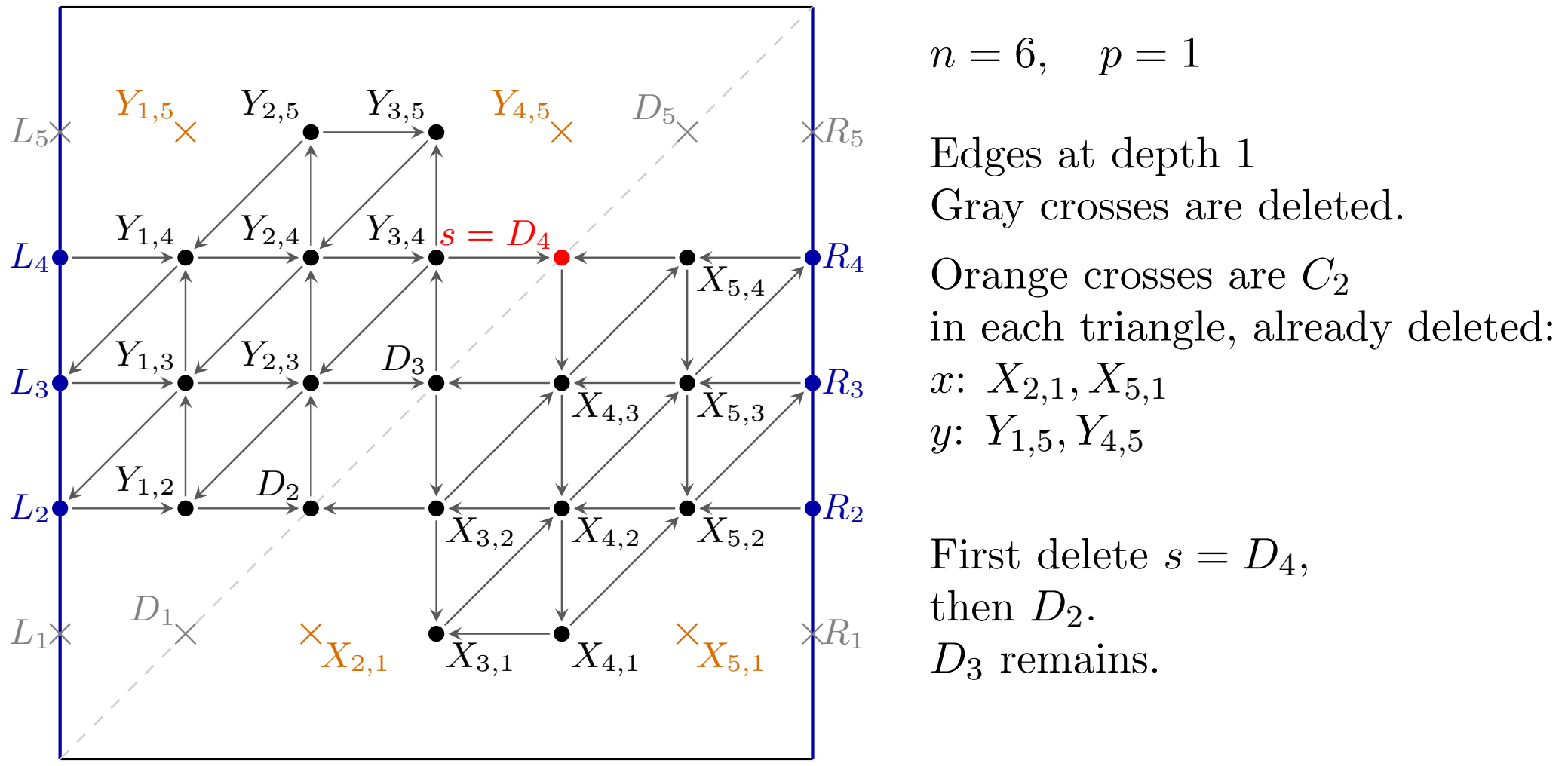
\begin{figure}
\centering
\begin{tikzpicture}[x=1.1cm,y=1.1cm,line width=.45pt,font=\small,>=stealth]
\draw[black] (0,0)--(6,0) (0,6)--(6,6);
\draw[blue!65!black,line width=.8pt] (0,0)--(0,6) (6,0)--(6,6);
\draw[dashed,gray!40] (0,0)--(6,6);
\draw[->,black!65,shorten >=3pt,shorten <=3pt] (0,2)--(1,2);
\draw[->,black!65,shorten >=3pt,shorten <=3pt] (0,3)--(1,3);
\draw[->,black!65,shorten >=3pt,shorten <=3pt] (0,4)--(1,4);
\draw[->,black!65,shorten >=3pt,shorten <=3pt] (1,2)--(1,3);
\draw[->,black!65,shorten >=3pt,shorten <=3pt] (1,2)--(2,2);
\draw[->,black!65,shorten >=3pt,shorten <=3pt] (1,3)--(0,2);
\draw[->,black!65,shorten >=3pt,shorten <=3pt] (1,3)--(1,4);
\draw[->,black!65,shorten >=3pt,shorten <=3pt] (1,3)--(2,3);
\draw[->,black!65,shorten >=3pt,shorten <=3pt] (1,4)--(0,3);
\draw[->,black!65,shorten >=3pt,shorten <=3pt] (1,4)--(2,4);
\draw[->,black!65,shorten >=3pt,shorten <=3pt] (2,2)--(2,3);
\draw[->,black!65,shorten >=3pt,shorten <=3pt] (2,3)--(1,2);
\draw[->,black!65,shorten >=3pt,shorten <=3pt] (2,3)--(2,4);
\draw[->,black!65,shorten >=3pt,shorten <=3pt] (2,3)--(3,3);
\draw[->,black!65,shorten >=3pt,shorten <=3pt] (2,4)--(1,3);
\draw[->,black!65,shorten >=3pt,shorten <=3pt] (2,4)--(2,5);
\draw[->,black!65,shorten >=3pt,shorten <=3pt] (2,4)--(3,4);
\draw[->,black!65,shorten >=3pt,shorten <=3pt] (2,5)--(1,4);
\draw[->,black!65,shorten >=3pt,shorten <=3pt] (2,5)--(3,5);
\draw[->,black!65,shorten >=3pt,shorten <=3pt] (3,1)--(4,2);
\draw[->,black!65,shorten >=3pt,shorten <=3pt] (3,2)--(2,2);
\draw[->,black!65,shorten >=3pt,shorten <=3pt] (3,2)--(3,1);
\draw[->,black!65,shorten >=3pt,shorten <=3pt] (3,2)--(4,3);
\draw[->,black!65,shorten >=3pt,shorten <=3pt] (3,3)--(3,2);
\draw[->,black!65,shorten >=3pt,shorten <=3pt] (3,3)--(3,4);
\draw[->,black!65,shorten >=3pt,shorten <=3pt] (3,4)--(2,3);
\draw[->,black!65,shorten >=3pt,shorten <=3pt] (3,4)--(3,5);
\draw[->,black!65,shorten >=3pt,shorten <=3pt] (3,4)--(4,4);
\draw[->,black!65,shorten >=3pt,shorten <=3pt] (3,5)--(2,4);
\draw[->,black!65,shorten >=3pt,shorten <=3pt] (4,1)--(3,1);
\draw[->,black!65,shorten >=3pt,shorten <=3pt] (4,1)--(5,2);
\draw[->,black!65,shorten >=3pt,shorten <=3pt] (4,2)--(3,2);
\draw[->,black!65,shorten >=3pt,shorten <=3pt] (4,2)--(4,1);
\draw[->,black!65,shorten >=3pt,shorten <=3pt] (4,2)--(5,3);
\draw[->,black!65,shorten >=3pt,shorten <=3pt] (4,3)--(3,3);
\draw[->,black!65,shorten >=3pt,shorten <=3pt] (4,3)--(4,2);
\draw[->,black!65,shorten >=3pt,shorten <=3pt] (4,3)--(5,4);
\draw[->,black!65,shorten >=3pt,shorten <=3pt] (4,4)--(4,3);
\draw[->,black!65,shorten >=3pt,shorten <=3pt] (5,2)--(4,2);
\draw[->,black!65,shorten >=3pt,shorten <=3pt] (5,2)--(6,3);
\draw[->,black!65,shorten >=3pt,shorten <=3pt] (5,3)--(4,3);
\draw[->,black!65,shorten >=3pt,shorten <=3pt] (5,3)--(5,2);
\draw[->,black!65,shorten >=3pt,shorten <=3pt] (5,3)--(6,4);
\draw[->,black!65,shorten >=3pt,shorten <=3pt] (5,4)--(4,4);
\draw[->,black!65,shorten >=3pt,shorten <=3pt] (5,4)--(5,3);
\draw[->,black!65,shorten >=3pt,shorten <=3pt] (6,2)--(5,2);
\draw[->,black!65,shorten >=3pt,shorten <=3pt] (6,3)--(5,3);
\draw[->,black!65,shorten >=3pt,shorten <=3pt] (6,4)--(5,4);
\fill[blue!65!black] (0,2) circle (2pt);
\node[left,text=blue!65!black,font=\scriptsize,inner sep=2pt] at (0,2) {$L_2$};
\fill[blue!65!black] (0,3) circle (2pt);
\node[left,text=blue!65!black,font=\scriptsize,inner sep=2pt] at (0,3) {$L_3$};
\fill[blue!65!black] (0,4) circle (2pt);
\node[left,text=blue!65!black,font=\scriptsize,inner sep=2pt] at (0,4) {$L_4$};
\fill[black] (1,2) circle (2pt);
\node[above left,text=black,font=\scriptsize,inner sep=2pt] at (1,2) {$Y_{1,2}$};
\fill[black] (1,3) circle (2pt);
\node[above left,text=black,font=\scriptsize,inner sep=2pt] at (1,3) {$Y_{1,3}$};
\fill[black] (1,4) circle (2pt);
\node[above left,text=black,font=\scriptsize,inner sep=2pt] at (1,4) {$Y_{1,4}$};
\fill[black] (2,2) circle (2pt);
\node[above left,text=black,font=\scriptsize,inner sep=2pt] at (2,2) {$D_2$};
\fill[black] (2,3) circle (2pt);
\node[above left,text=black,font=\scriptsize,inner sep=2pt] at (2,3) {$Y_{2,3}$};
\fill[black] (2,4) circle (2pt);
\node[above left,text=black,font=\scriptsize,inner sep=2pt] at (2,4) {$Y_{2,4}$};
\fill[black] (2,5) circle (2pt);
\node[above left,text=black,font=\scriptsize,inner sep=2pt] at (2,5) {$Y_{2,5}$};
\fill[black] (3,1) circle (2pt);
\node[below right,text=black,font=\scriptsize,inner sep=2pt] at (3,1) {$X_{3,1}$};
\fill[black] (3,2) circle (2pt);
\node[below right,text=black,font=\scriptsize,inner sep=2pt] at (3,2) {$X_{3,2}$};
\fill[black] (3,3) circle (2pt);
\node[above left,text=black,font=\scriptsize,inner sep=2pt] at (3,3) {$D_3$};
\fill[black] (3,4) circle (2pt);
\node[above left,text=black,font=\scriptsize,inner sep=2pt] at (3,4) {$Y_{3,4}$};
\fill[black] (3,5) circle (2pt);
\node[above left,text=black,font=\scriptsize,inner sep=2pt] at (3,5) {$Y_{3,5}$};
\fill[black] (4,1) circle (2pt);
\node[below right,text=black,font=\scriptsize,inner sep=2pt] at (4,1) {$X_{4,1}$};
\fill[black] (4,2) circle (2pt);
\node[below right,text=black,font=\scriptsize,inner sep=2pt] at (4,2) {$X_{4,2}$};
\fill[black] (4,3) circle (2pt);
\node[below right,text=black,font=\scriptsize,inner sep=2pt] at (4,3) {$X_{4,3}$};
\fill[red] (4,4) circle (2pt);
\node[above left,text=red,font=\scriptsize,inner sep=2pt] at (4,4) {$s=D_4$};
\fill[black] (5,2) circle (2pt);
\node[below right,text=black,font=\scriptsize,inner sep=2pt] at (5,2) {$X_{5,2}$};
\fill[black] (5,3) circle (2pt);
\node[below right,text=black,font=\scriptsize,inner sep=2pt] at (5,3) {$X_{5,3}$};
\fill[black] (5,4) circle (2pt);
\node[below right,text=black,font=\scriptsize,inner sep=2pt] at (5,4) {$X_{5,4}$};
\fill[blue!65!black] (6,2) circle (2pt);
\node[right,text=blue!65!black,font=\scriptsize,inner sep=2pt] at (6,2) {$R_2$};
\fill[blue!65!black] (6,3) circle (2pt);
\node[right,text=blue!65!black,font=\scriptsize,inner sep=2pt] at (6,3) {$R_3$};
\fill[blue!65!black] (6,4) circle (2pt);
\node[right,text=blue!65!black,font=\scriptsize,inner sep=2pt] at (6,4) {$R_4$};
\node[text=gray,inner sep=0pt] at (0,1) {$\times$};
\node[left,text=gray,font=\scriptsize,inner sep=2pt] at (0,1) {$L_1$};
\node[text=gray,inner sep=0pt] at (0,5) {$\times$};
\node[left,text=gray,font=\scriptsize,inner sep=2pt] at (0,5) {$L_5$};
\node[text=gray,inner sep=0pt] at (1,1) {$\times$};
\node[above left,text=gray,font=\scriptsize,inner sep=2pt] at (1,1) {$D_1$};
\node[text=gray,inner sep=0pt] at (5,5) {$\times$};
\node[above left,text=gray,font=\scriptsize,inner sep=2pt] at (5,5) {$D_5$};
\node[text=gray,inner sep=0pt] at (6,1) {$\times$};
\node[right,text=gray,font=\scriptsize,inner sep=2pt] at (6,1) {$R_1$};
\node[text=gray,inner sep=0pt] at (6,5) {$\times$};
\node[right,text=gray,font=\scriptsize,inner sep=2pt] at (6,5) {$R_5$};
\node[text=orange!85!black,inner sep=0pt] at (1,5) {$\times$};
\node[above left,text=orange!85!black,font=\scriptsize,inner sep=2pt] at (1,5) {$Y_{1,5}$};
\node[text=orange!85!black,inner sep=0pt] at (2,1) {$\times$};
\node[below right,text=orange!85!black,font=\scriptsize,inner sep=2pt] at (2,1) {$X_{2,1}$};
\node[text=orange!85!black,inner sep=0pt] at (4,5) {$\times$};
\node[above left,text=orange!85!black,font=\scriptsize,inner sep=2pt] at (4,5) {$Y_{4,5}$};
\node[text=orange!85!black,inner sep=0pt] at (5,1) {$\times$};
\node[below right,text=orange!85!black,font=\scriptsize,inner sep=2pt] at (5,1) {$X_{5,1}$};
\node[anchor=west] at (6.8,5.6) {$n=6,\quad p=1$};
\node[anchor=west,align=left] at (6.8,4.6) {Edges at depth $1$\\Gray crosses are deleted.};
\node[anchor=west,align=left] at (6.8,3.2) {Orange crosses are $C_2$\\in each triangle, already deleted:\\$x$: $X_{2,1},X_{5,1}$\\$y$: $Y_{1,5},Y_{4,5}$};
\node[anchor=west,align=left] at (6.8,1.2) {First delete $s=D_4$,\\then $D_2$.\\$D_3$ remains.};
\end{tikzpicture}
\caption{The starting quiver in Lemma~\ref{lem:edge-layer} for $n=6$, $p=1$.
Here $C_2=\{(4,1,1),(1,4,1)\}$ in the coordinates with boundary side $c=0$.}
\label{fig:triangle-band}
\label{fig:triangle-band-band}
\end{figure}

For two consecutive triangles, write $x(a,b,c)$ and $y(a,b,c)$
for the vertices labelled $(a,b,c)$ in the first and second triangles, respectively.
Place them in the square $[0,n]^2$ by
\[
 x(a,b,c)=(b+c,c),\qquad y(a,b,c)=(b,b+c).
\]
The arrows in each triangle are given by \eqref{eq:small-triangle}.
The boundary sides are $c=0$ in $x$ and $a=0$ in $y$.
We use the labels
\begin{equation}\label{eq:grid-labels}
 \begin{split}
 D_k&=x(n-k,0,k)=y(n-k,k,0)=(k,k),\\
 X_{i,j}&=x(n-i,i-j,j)=(i,j),\qquad 1\le j<i\le n-1,\\
 Y_{i,j}&=y(n-j,i,j-i)=(i,j),\qquad 1\le i<j\le n-1,\\
 L_j&=y(n-j,0,j)=(0,j),\qquad
 R_j=x(0,n-j,j)=(n,j).
 \end{split}
\end{equation}
Here $1\le k,j\le n-1$ for the edge vertices.
Set $X_{n,j}=R_j$.
For the $y$ triangle, use the cyclically reordered coordinates
$(b,c,a)$ in \eqref{eq:triangle-corners}.
Rotation of the square by $180^\circ$ about $(n/2,n/2)$ sends
$(i,j)$ to $(n-i,n-j)$. It preserves arrows and interchanges $D_k$ and $D_{n-k}$.
Every consecutive pair of strip triangles has this quiver
or its arrow reversal.
Figure~\ref{fig:triangle-band} illustrates the hypotheses of Lemma~\ref{lem:edge-layer}.

\begin{lemma}\label{lem:edge-layer}
Let $0\le p<\lfloor n/2\rfloor$.
In the two-triangle configuration \eqref{eq:grid-labels}, suppose that
all three edges have been deleted to depth $p$ and that $C_{p+1}$
has been deleted in each triangle.
Then the common edge can be deleted to depth $p+1$
by a sequence of deletion moves.
All mutations and deletions take place in the interiors of the two triangles
or on their common edge.
\end{lemma}
\begin{proof}
We delete $D_{n-p-1}$, then use the rotated square to delete $D_{p+1}$.
If $n=2p+2$, then $D_{n-p-1}=D_{p+1}$ and
\[
 \rho_Q(D_{p+1})=-\mathbf e_{x(p,1,p+1)}-\mathbf e_{y(p+1,p,1)}.
\]
Delete $D_{p+1}$ without mutations. Assume $n\ge2p+3$ and put
$r=p+2$, $t=n-p-1$, and $s=D_t$.

For the deletion move at $s$, use the sequence
\begin{equation}\label{eq:pruned-word}
 D_{p+1},\ldots,D_{t-1};\quad
 X_{r,r-1},\ldots,X_{r,1};\quad\cdots;\quad
 X_{t,t-1},\ldots,X_{t,1}.
\end{equation}
The deleted interior vertices in this triangle satisfy
$i\le p+1$ or $i-j\ge t$.
The interior vertices in \eqref{eq:pruned-word} satisfy $r\le i\le t$ and $1\le j<i$,
so none has been deleted. The sequence avoids $s$ and all vertices $L_k,R_k$.

For $r\le i<t$, just before mutating $X_{i,i-1}$, we have
\begin{equation}\label{eq:pruned-top}
 \rho_{\mathrm{pre}}(X_{i,i-1})
 =-\mathbf e_{D_{i-1}}+\mathbf e_{D_i}
  -\mathbf e_{X_{i+1,i-1}}+\mathbf e_{X_{i,i-2}}.
\end{equation}
For $i=t$, we have
\begin{equation}\label{eq:pruned-last}
 \begin{split}
 \rho_{\mathrm{pre}}(X_{t,t-1})={}&-\mathbf e_{D_{t-1}}-\mathbf e_s
                 -\mathbf e_{X_{t+1,t-1}}\\
              &+\mathbf e_{X_{t+1,t}}+\mathbf e_{Y_{t-1,t}}
                 +\mathbf e_{X_{t,t-2}}.
 \end{split}
\end{equation}
For $r\le i\le t$ and $1\le j<i-1$, let $K_{ij}$ denote the vertex
\[
 K_{ij}=\begin{cases}D_{p+1},&i=r,\\X_{i-1,j},&i>r.\end{cases}
\]
Just before mutating $X_{i,j}$, we have
\begin{equation}\label{eq:pruned-lower}
 \rho_{\mathrm{pre}}(X_{i,j})=-\mathbf e_{K_{ij}}
        -\mathbf e_{X_{i+1,j}}+\mathbf e_{X_{i,j+1}}
        +\mathbf e_{X_{i,j-1}}-\mathbf1_{i=t}\mathbf e_s.
\end{equation}
 
Among the diagonal vertices, only $D_{i-1}$ and $D_i$ are adjacent to $X_{i,i-1}$.
The cycles in \eqref{eq:small-triangle} give \eqref{eq:pruned-top}
after mutating both, and \eqref{eq:pruned-last} when $i=t$, since $D_t$ is not mutated.
The intervening mutations at $X_{h,j}$ with $h<i$ leave the arrows at $X_{i,i-1}$ unchanged.
After the mutations at $D_{p+1},\ldots,D_{t-1}$, for $1\le j<i-1$ we have
\[
 \begin{aligned}
 \rho_Q(X_{i,j})={}&\mathbf e_{X_{i+1,j+1}}-\mathbf e_{X_{i+1,j}}
 -\mathbf e_{X_{i,j+1}}\\
 &+\mathbf e_{X_{i,j-1}}+\mathbf e_{X_{i-1,j}}
 -\mathbf e_{X_{i-1,j-1}}.
 \end{aligned}
\]
For $i=r$, mutation at $X_{r,j+1}$ cancels the arrow
to $X_{r+1,j+1}$ and creates the arrow from $D_{p+1}$.
For $i>r$, mutation at $X_{i-1,j}$ cancels the arrow from $X_{i-1,j-1}$
and creates an arrow to $X_{i-1,j+1}$ (or $D_{i-1}$ if $j=i-2$).
Mutation at $X_{i,j+1}$ cancels this new arrow and the arrow to $X_{i+1,j+1}$.
These are the only earlier mutations that change the arrows at $X_{i,j}$,
and they give \eqref{eq:pruned-lower}.

Only mutations at $X_{t,j}$ affect $s$.
Before mutating $X_{t,t-1}$, we have
\[
 \rho_{\mathrm{pre}}(s)=\mathbf e_{X_{t,t-1}}-\mathbf e_{X_{t+1,t}}
                              -\mathbf e_{Y_{t-1,t}}.
\]
The possible neighbour $Y_{t,t+1}$ lies on the boundary
if $p=0$, and in $C_{p+1}$ if $p\ge1$, so it is absent.
Mutation at $X_{t,t-1}$ cancels the arrows $X_{t+1,t}\to s$ and $Y_{t-1,t}\to s$.
After mutating $X_{t,j}$, the current quiver $Q$ has
\[
 \rho_Q(s)=-\mathbf e_{X_{t,j}}+\mathbf e_{X_{t,j-1}}
 \qquad(j=t-1,t-2,\ldots,1).
\]
Since $X_{t,0}$ is absent,
\begin{equation}\label{eq:pruned-final}
 \rho_{\mathrm{final}}(s)=-\mathbf e_{X_{t,1}}.
\end{equation}
Delete $s$ and undo \eqref{eq:pruned-word}.

To construct a deletion move at $D_{p+1}$, rotate the square by $180^\circ$
and relabel its vertices by \eqref{eq:grid-labels}.
The original $D_{p+1}$ is now labelled $D_t$, and the deleted original $D_t$
is now labelled $D_{p+1}$. Omit the first mutation in \eqref{eq:pruned-word}.
If $r<t$, replace $D_{p+1}$ by $Y_{p+1,p+2}$ in
\eqref{eq:pruned-top} and \eqref{eq:pruned-lower} for $i=r$.
The other formulas, including \eqref{eq:pruned-final}, are unchanged.
Indeed, the mutation at $D_r$ leaves an arrow from $Y_{p+1,p+2}$
which was cancelled by the omitted mutation.
The calculation above then applies for $j=r-1,\ldots,1$.
If $r=t$, omit
$-\mathbf e_{D_{t-1}}+\mathbf e_{Y_{t-1,t}}$
from \eqref{eq:pruned-last} and omit $-\mathbf e_{K_{tj}}$
from \eqref{eq:pruned-lower}.
Returning to the original labels, these two cases give
\[
 \rho_{\mathrm{final}}(D_{p+1})=
 \begin{cases}
 -\mathbf e_{Y_{p+1,n-1}},&r<t,\\
 -\mathbf e_{Y_{p+1,n-1}}-\mathbf e_{X_{p+2,p+1}},&r=t.
 \end{cases}
\]
Delete $D_{p+1}$ and undo the mutations.
\end{proof}

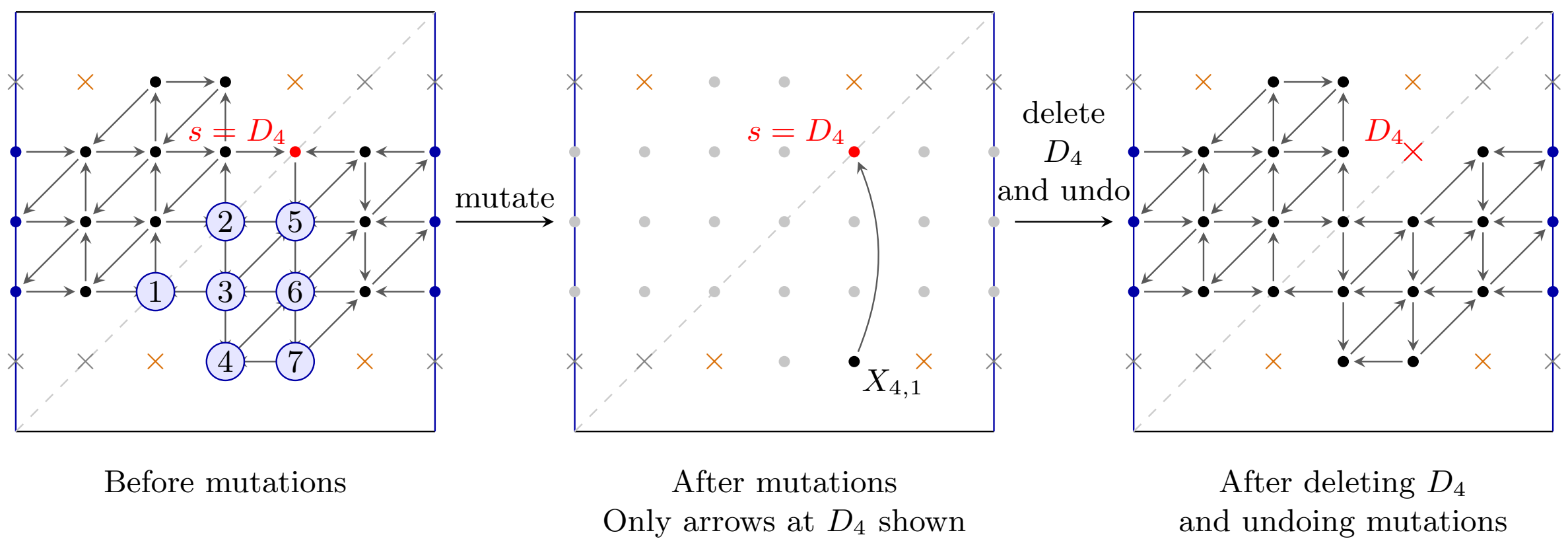
\begin{figure}
\centering
\begin{tikzpicture}[x=.7cm,y=.7cm,line width=.45pt,font=\small,>=stealth]
\begin{scope}[shift={(0,0)}]
\draw[black] (0,0)--(6,0) (0,6)--(6,6);
\draw[blue!65!black] (0,0)--(0,6) (6,0)--(6,6);
\draw[dashed,gray!40] (0,0)--(6,6);
\draw[->,black!65,shorten >=3pt,shorten <=3pt] (0,2)--(1,2);
\draw[->,black!65,shorten >=3pt,shorten <=3pt] (0,3)--(1,3);
\draw[->,black!65,shorten >=3pt,shorten <=3pt] (0,4)--(1,4);
\draw[->,black!65,shorten >=3pt,shorten <=3pt] (1,2)--(1,3);
\draw[->,black!65,shorten >=3pt,shorten <=3pt] (1,2)--(2,2);
\draw[->,black!65,shorten >=3pt,shorten <=3pt] (1,3)--(0,2);
\draw[->,black!65,shorten >=3pt,shorten <=3pt] (1,3)--(1,4);
\draw[->,black!65,shorten >=3pt,shorten <=3pt] (1,3)--(2,3);
\draw[->,black!65,shorten >=3pt,shorten <=3pt] (1,4)--(0,3);
\draw[->,black!65,shorten >=3pt,shorten <=3pt] (1,4)--(2,4);
\draw[->,black!65,shorten >=3pt,shorten <=3pt] (2,2)--(2,3);
\draw[->,black!65,shorten >=3pt,shorten <=3pt] (2,3)--(1,2);
\draw[->,black!65,shorten >=3pt,shorten <=3pt] (2,3)--(2,4);
\draw[->,black!65,shorten >=3pt,shorten <=3pt] (2,3)--(3,3);
\draw[->,black!65,shorten >=3pt,shorten <=3pt] (2,4)--(1,3);
\draw[->,black!65,shorten >=3pt,shorten <=3pt] (2,4)--(2,5);
\draw[->,black!65,shorten >=3pt,shorten <=3pt] (2,4)--(3,4);
\draw[->,black!65,shorten >=3pt,shorten <=3pt] (2,5)--(1,4);
\draw[->,black!65,shorten >=3pt,shorten <=3pt] (2,5)--(3,5);
\draw[->,black!65,shorten >=3pt,shorten <=3pt] (3,1)--(4,2);
\draw[->,black!65,shorten >=3pt,shorten <=3pt] (3,2)--(2,2);
\draw[->,black!65,shorten >=3pt,shorten <=3pt] (3,2)--(3,1);
\draw[->,black!65,shorten >=3pt,shorten <=3pt] (3,2)--(4,3);
\draw[->,black!65,shorten >=3pt,shorten <=3pt] (3,3)--(3,2);
\draw[->,black!65,shorten >=3pt,shorten <=3pt] (3,3)--(3,4);
\draw[->,black!65,shorten >=3pt,shorten <=3pt] (3,4)--(2,3);
\draw[->,black!65,shorten >=3pt,shorten <=3pt] (3,4)--(3,5);
\draw[->,black!65,shorten >=3pt,shorten <=3pt] (3,4)--(4,4);
\draw[->,black!65,shorten >=3pt,shorten <=3pt] (3,5)--(2,4);
\draw[->,black!65,shorten >=3pt,shorten <=3pt] (4,1)--(3,1);
\draw[->,black!65,shorten >=3pt,shorten <=3pt] (4,1)--(5,2);
\draw[->,black!65,shorten >=3pt,shorten <=3pt] (4,2)--(3,2);
\draw[->,black!65,shorten >=3pt,shorten <=3pt] (4,2)--(4,1);
\draw[->,black!65,shorten >=3pt,shorten <=3pt] (4,2)--(5,3);
\draw[->,black!65,shorten >=3pt,shorten <=3pt] (4,3)--(3,3);
\draw[->,black!65,shorten >=3pt,shorten <=3pt] (4,3)--(4,2);
\draw[->,black!65,shorten >=3pt,shorten <=3pt] (4,3)--(5,4);
\draw[->,black!65,shorten >=3pt,shorten <=3pt] (4,4)--(4,3);
\draw[->,black!65,shorten >=3pt,shorten <=3pt] (5,2)--(4,2);
\draw[->,black!65,shorten >=3pt,shorten <=3pt] (5,2)--(6,3);
\draw[->,black!65,shorten >=3pt,shorten <=3pt] (5,3)--(4,3);
\draw[->,black!65,shorten >=3pt,shorten <=3pt] (5,3)--(5,2);
\draw[->,black!65,shorten >=3pt,shorten <=3pt] (5,3)--(6,4);
\draw[->,black!65,shorten >=3pt,shorten <=3pt] (5,4)--(4,4);
\draw[->,black!65,shorten >=3pt,shorten <=3pt] (5,4)--(5,3);
\draw[->,black!65,shorten >=3pt,shorten <=3pt] (6,2)--(5,2);
\draw[->,black!65,shorten >=3pt,shorten <=3pt] (6,3)--(5,3);
\draw[->,black!65,shorten >=3pt,shorten <=3pt] (6,4)--(5,4);
\fill[blue!65!black] (0,2) circle (1.6pt);
\fill[blue!65!black] (0,3) circle (1.6pt);
\fill[blue!65!black] (0,4) circle (1.6pt);
\fill[black] (1,2) circle (1.6pt);
\fill[black] (1,3) circle (1.6pt);
\fill[black] (1,4) circle (1.6pt);
\fill[black] (2,2) circle (1.6pt);
\fill[black] (2,3) circle (1.6pt);
\fill[black] (2,4) circle (1.6pt);
\fill[black] (2,5) circle (1.6pt);
\fill[black] (3,1) circle (1.6pt);
\fill[black] (3,2) circle (1.6pt);
\fill[black] (3,3) circle (1.6pt);
\fill[black] (3,4) circle (1.6pt);
\fill[black] (3,5) circle (1.6pt);
\fill[black] (4,1) circle (1.6pt);
\fill[black] (4,2) circle (1.6pt);
\fill[black] (4,3) circle (1.6pt);
\fill[red] (4,4) circle (1.6pt);
\fill[black] (5,2) circle (1.6pt);
\fill[black] (5,3) circle (1.6pt);
\fill[black] (5,4) circle (1.6pt);
\fill[blue!65!black] (6,2) circle (1.6pt);
\fill[blue!65!black] (6,3) circle (1.6pt);
\fill[blue!65!black] (6,4) circle (1.6pt);
\node[text=gray,font=\scriptsize,inner sep=0pt] at (0,1) {$\times$};
\node[text=gray,font=\scriptsize,inner sep=0pt] at (0,5) {$\times$};
\node[text=gray,font=\scriptsize,inner sep=0pt] at (1,1) {$\times$};
\node[text=gray,font=\scriptsize,inner sep=0pt] at (5,5) {$\times$};
\node[text=gray,font=\scriptsize,inner sep=0pt] at (6,1) {$\times$};
\node[text=gray,font=\scriptsize,inner sep=0pt] at (6,5) {$\times$};
\node[text=orange!85!black,font=\scriptsize,inner sep=0pt] at (1,5) {$\times$};
\node[text=orange!85!black,font=\scriptsize,inner sep=0pt] at (2,1) {$\times$};
\node[text=orange!85!black,font=\scriptsize,inner sep=0pt] at (4,5) {$\times$};
\node[text=orange!85!black,font=\scriptsize,inner sep=0pt] at (5,1) {$\times$};
\node[circle,draw=blue!65!black,fill=blue!10,inner sep=1.2pt,font=\scriptsize] at (2,2) {1};
\node[circle,draw=blue!65!black,fill=blue!10,inner sep=1.2pt,font=\scriptsize] at (3,3) {2};
\node[circle,draw=blue!65!black,fill=blue!10,inner sep=1.2pt,font=\scriptsize] at (3,2) {3};
\node[circle,draw=blue!65!black,fill=blue!10,inner sep=1.2pt,font=\scriptsize] at (3,1) {4};
\node[circle,draw=blue!65!black,fill=blue!10,inner sep=1.2pt,font=\scriptsize] at (4,3) {5};
\node[circle,draw=blue!65!black,fill=blue!10,inner sep=1.2pt,font=\scriptsize] at (4,2) {6};
\node[circle,draw=blue!65!black,fill=blue!10,inner sep=1.2pt,font=\scriptsize] at (4,1) {7};
\node[above left,text=red,font=\scriptsize,inner sep=2pt] at (4,4) {$s=D_4$};
\node[below,align=center,font=\scriptsize] at (3,-.35) {Before mutations};
\end{scope}
\begin{scope}[shift={(8,0)}]
\draw[black] (0,0)--(6,0) (0,6)--(6,6);
\draw[blue!65!black] (0,0)--(0,6) (6,0)--(6,6);
\draw[dashed,gray!40] (0,0)--(6,6);
\draw[->,black!65,shorten >=3pt,shorten <=3pt] (4,1) to[bend right=22] (4,4);
\fill[gray!45] (0,2) circle (1.6pt);
\fill[gray!45] (0,3) circle (1.6pt);
\fill[gray!45] (0,4) circle (1.6pt);
\fill[gray!45] (1,2) circle (1.6pt);
\fill[gray!45] (1,3) circle (1.6pt);
\fill[gray!45] (1,4) circle (1.6pt);
\fill[gray!45] (2,2) circle (1.6pt);
\fill[gray!45] (2,3) circle (1.6pt);
\fill[gray!45] (2,4) circle (1.6pt);
\fill[gray!45] (2,5) circle (1.6pt);
\fill[gray!45] (3,1) circle (1.6pt);
\fill[gray!45] (3,2) circle (1.6pt);
\fill[gray!45] (3,3) circle (1.6pt);
\fill[gray!45] (3,4) circle (1.6pt);
\fill[gray!45] (3,5) circle (1.6pt);
\fill[black] (4,1) circle (1.6pt);
\fill[gray!45] (4,2) circle (1.6pt);
\fill[gray!45] (4,3) circle (1.6pt);
\fill[red] (4,4) circle (1.6pt);
\fill[gray!45] (5,2) circle (1.6pt);
\fill[gray!45] (5,3) circle (1.6pt);
\fill[gray!45] (5,4) circle (1.6pt);
\fill[gray!45] (6,2) circle (1.6pt);
\fill[gray!45] (6,3) circle (1.6pt);
\fill[gray!45] (6,4) circle (1.6pt);
\node[text=gray,font=\scriptsize,inner sep=0pt] at (0,1) {$\times$};
\node[text=gray,font=\scriptsize,inner sep=0pt] at (0,5) {$\times$};
\node[text=gray,font=\scriptsize,inner sep=0pt] at (1,1) {$\times$};
\node[text=gray,font=\scriptsize,inner sep=0pt] at (5,5) {$\times$};
\node[text=gray,font=\scriptsize,inner sep=0pt] at (6,1) {$\times$};
\node[text=gray,font=\scriptsize,inner sep=0pt] at (6,5) {$\times$};
\node[text=orange!85!black,font=\scriptsize,inner sep=0pt] at (1,5) {$\times$};
\node[text=orange!85!black,font=\scriptsize,inner sep=0pt] at (2,1) {$\times$};
\node[text=orange!85!black,font=\scriptsize,inner sep=0pt] at (4,5) {$\times$};
\node[text=orange!85!black,font=\scriptsize,inner sep=0pt] at (5,1) {$\times$};
\node[above left,text=red,font=\scriptsize,inner sep=2pt] at (4,4) {$s=D_4$};
\node[below right,font=\scriptsize,inner sep=2pt] at (4,1) {$X_{4,1}$};
\node[below,align=center,font=\scriptsize] at (3,-.35) {After mutations\\Only arrows at $D_4$ shown};
\end{scope}
\begin{scope}[shift={(16,0)}]
\draw[black] (0,0)--(6,0) (0,6)--(6,6);
\draw[blue!65!black] (0,0)--(0,6) (6,0)--(6,6);
\draw[dashed,gray!40] (0,0)--(6,6);
\draw[->,black!65,shorten >=3pt,shorten <=3pt] (0,2)--(1,2);
\draw[->,black!65,shorten >=3pt,shorten <=3pt] (0,3)--(1,3);
\draw[->,black!65,shorten >=3pt,shorten <=3pt] (0,4)--(1,4);
\draw[->,black!65,shorten >=3pt,shorten <=3pt] (1,2)--(1,3);
\draw[->,black!65,shorten >=3pt,shorten <=3pt] (1,2)--(2,2);
\draw[->,black!65,shorten >=3pt,shorten <=3pt] (1,3)--(0,2);
\draw[->,black!65,shorten >=3pt,shorten <=3pt] (1,3)--(1,4);
\draw[->,black!65,shorten >=3pt,shorten <=3pt] (1,3)--(2,3);
\draw[->,black!65,shorten >=3pt,shorten <=3pt] (1,4)--(0,3);
\draw[->,black!65,shorten >=3pt,shorten <=3pt] (1,4)--(2,4);
\draw[->,black!65,shorten >=3pt,shorten <=3pt] (2,2)--(2,3);
\draw[->,black!65,shorten >=3pt,shorten <=3pt] (2,3)--(1,2);
\draw[->,black!65,shorten >=3pt,shorten <=3pt] (2,3)--(2,4);
\draw[->,black!65,shorten >=3pt,shorten <=3pt] (2,3)--(3,3);
\draw[->,black!65,shorten >=3pt,shorten <=3pt] (2,4)--(1,3);
\draw[->,black!65,shorten >=3pt,shorten <=3pt] (2,4)--(2,5);
\draw[->,black!65,shorten >=3pt,shorten <=3pt] (2,4)--(3,4);
\draw[->,black!65,shorten >=3pt,shorten <=3pt] (2,5)--(1,4);
\draw[->,black!65,shorten >=3pt,shorten <=3pt] (2,5)--(3,5);
\draw[->,black!65,shorten >=3pt,shorten <=3pt] (3,1)--(4,2);
\draw[->,black!65,shorten >=3pt,shorten <=3pt] (3,2)--(2,2);
\draw[->,black!65,shorten >=3pt,shorten <=3pt] (3,2)--(3,1);
\draw[->,black!65,shorten >=3pt,shorten <=3pt] (3,2)--(4,3);
\draw[->,black!65,shorten >=3pt,shorten <=3pt] (3,3)--(3,2);
\draw[->,black!65,shorten >=3pt,shorten <=3pt] (3,3)--(3,4);
\draw[->,black!65,shorten >=3pt,shorten <=3pt] (3,4)--(2,3);
\draw[->,black!65,shorten >=3pt,shorten <=3pt] (3,4)--(3,5);
\draw[->,black!65,shorten >=3pt,shorten <=3pt] (3,5)--(2,4);
\draw[->,black!65,shorten >=3pt,shorten <=3pt] (4,1)--(3,1);
\draw[->,black!65,shorten >=3pt,shorten <=3pt] (4,1)--(5,2);
\draw[->,black!65,shorten >=3pt,shorten <=3pt] (4,2)--(3,2);
\draw[->,black!65,shorten >=3pt,shorten <=3pt] (4,2)--(4,1);
\draw[->,black!65,shorten >=3pt,shorten <=3pt] (4,2)--(5,3);
\draw[->,black!65,shorten >=3pt,shorten <=3pt] (4,3)--(3,3);
\draw[->,black!65,shorten >=3pt,shorten <=3pt] (4,3)--(4,2);
\draw[->,black!65,shorten >=3pt,shorten <=3pt] (4,3)--(5,4);
\draw[->,black!65,shorten >=3pt,shorten <=3pt] (5,2)--(4,2);
\draw[->,black!65,shorten >=3pt,shorten <=3pt] (5,2)--(6,3);
\draw[->,black!65,shorten >=3pt,shorten <=3pt] (5,3)--(4,3);
\draw[->,black!65,shorten >=3pt,shorten <=3pt] (5,3)--(5,2);
\draw[->,black!65,shorten >=3pt,shorten <=3pt] (5,3)--(6,4);
\draw[->,black!65,shorten >=3pt,shorten <=3pt] (5,4)--(5,3);
\draw[->,black!65,shorten >=3pt,shorten <=3pt] (6,2)--(5,2);
\draw[->,black!65,shorten >=3pt,shorten <=3pt] (6,3)--(5,3);
\draw[->,black!65,shorten >=3pt,shorten <=3pt] (6,4)--(5,4);
\fill[blue!65!black] (0,2) circle (1.6pt);
\fill[blue!65!black] (0,3) circle (1.6pt);
\fill[blue!65!black] (0,4) circle (1.6pt);
\fill[black] (1,2) circle (1.6pt);
\fill[black] (1,3) circle (1.6pt);
\fill[black] (1,4) circle (1.6pt);
\fill[black] (2,2) circle (1.6pt);
\fill[black] (2,3) circle (1.6pt);
\fill[black] (2,4) circle (1.6pt);
\fill[black] (2,5) circle (1.6pt);
\fill[black] (3,1) circle (1.6pt);
\fill[black] (3,2) circle (1.6pt);
\fill[black] (3,3) circle (1.6pt);
\fill[black] (3,4) circle (1.6pt);
\fill[black] (3,5) circle (1.6pt);
\fill[black] (4,1) circle (1.6pt);
\fill[black] (4,2) circle (1.6pt);
\fill[black] (4,3) circle (1.6pt);
\fill[black] (5,2) circle (1.6pt);
\fill[black] (5,3) circle (1.6pt);
\fill[black] (5,4) circle (1.6pt);
\fill[blue!65!black] (6,2) circle (1.6pt);
\fill[blue!65!black] (6,3) circle (1.6pt);
\fill[blue!65!black] (6,4) circle (1.6pt);
\node[text=gray,font=\scriptsize,inner sep=0pt] at (0,1) {$\times$};
\node[text=gray,font=\scriptsize,inner sep=0pt] at (0,5) {$\times$};
\node[text=gray,font=\scriptsize,inner sep=0pt] at (1,1) {$\times$};
\node[text=gray,font=\scriptsize,inner sep=0pt] at (5,5) {$\times$};
\node[text=gray,font=\scriptsize,inner sep=0pt] at (6,1) {$\times$};
\node[text=gray,font=\scriptsize,inner sep=0pt] at (6,5) {$\times$};
\node[text=orange!85!black,font=\scriptsize,inner sep=0pt] at (1,5) {$\times$};
\node[text=orange!85!black,font=\scriptsize,inner sep=0pt] at (2,1) {$\times$};
\node[text=orange!85!black,font=\scriptsize,inner sep=0pt] at (4,5) {$\times$};
\node[text=orange!85!black,font=\scriptsize,inner sep=0pt] at (5,1) {$\times$};
\node[text=red,inner sep=0pt] at (4,4) {$\times$};
\node[above left,text=red,font=\scriptsize,inner sep=2pt] at (4,4) {$D_4$};
\node[below,align=center,font=\scriptsize] at (3,-.35) {After deleting $D_4$\\and undoing mutations};
\end{scope}
\draw[->] (6.3,3)--node[above,font=\scriptsize] {mutate} (7.7,3);
\draw[->] (14.3,3)--(15.7,3);
\node[above,align=center,font=\scriptsize] at (15,3.1) {delete\\$D_4$\\and undo};
\end{tikzpicture}
\caption{A deletion move at $D_4$ for $n=6$, $p=1$. Numbers indicate the mutation order.}
\label{fig:sweep}
\end{figure}

For $n=6$ and $p=1$, start with the quiver in Figure~\ref{fig:triangle-band}.
The deletion move at $D_4$ uses the mutation sequence
\[
 D_2,D_3,X_{3,2},X_{3,1},X_{4,3},X_{4,2},X_{4,1}.
\]
After the first four mutations, the remaining three change $\rho_Q(D_4)$ as follows.
\[
 \begin{aligned}
 &\mathbf e_{X_{4,3}}-\mathbf e_{X_{5,4}}-\mathbf e_{Y_{3,4}}\\
 &\quad\xrightarrow{\ \mu_{X_{4,3}}\ }-\mathbf e_{X_{4,3}}+\mathbf e_{X_{4,2}}
 \xrightarrow{\ \mu_{X_{4,2}}\ }-\mathbf e_{X_{4,2}}+\mathbf e_{X_{4,1}}
 \xrightarrow{\ \mu_{X_{4,1}}\ }-\mathbf e_{X_{4,1}}.
 \end{aligned}
\]
Figure~\ref{fig:sweep} shows this deletion move.

\begin{proof}[Proof of Proposition~\ref{prop:strip}]
Use the deletion moves of Lemma~\ref{lem:edge-layer} with $p=0$
on $e_1,\ldots,e_{2N-1}$, in this order.
Since $C_1=\varnothing$, this deletes each edge to depth $1$.
Each adjacent pair of triangles has the quiver
\eqref{eq:grid-labels} or its arrow reversal.
Previous deletion moves on a neighbouring edge remove only vertices
labelled $L_k$ or $R_k$ in the next pair.
These vertices are neither mutated nor deleted by the moves for this pair,
so restricting to the remaining vertices
preserves the calculation and the source or sink condition.

After round $p$, we have deleted
\begin{equation}\label{eq:strip-induction}
 \begin{aligned}
 &\{k:k\le p\text{ or }k\ge n-p\}\text{ on }e_j,
       &&p\le j\le2N-p,\\
 &C_p\text{ in }\tau_i,
       &&p\le i\le2N-p+1.
 \end{aligned}
\end{equation}
The first round gives \eqref{eq:strip-induction} for $p=1$.
Suppose $1\le p<\lfloor n/2\rfloor$.
For $p+1\le i\le2N-p$, the edges $e_{i-1},e_i$ have both
been deleted to depth $p$.
Use the deletion moves of Lemma~\ref{lem:triangle-shell} in these triangles
to remove $C_{p+1}\setminus C_p$.
For $p+1\le j\le2N-p-1$, the pair $\tau_j,\tau_{j+1}$ now
satisfies the hypotheses of Lemma~\ref{lem:edge-layer}.
Perform its deletion moves on $e_j$, in increasing order of $j$.
As in the first round, extra deletions on $e_{j-1}$ and $e_{j+1}$ do
not affect the calculation. This proves \eqref{eq:strip-induction}
with $p$ replaced by $p+1$.

Put $q=\lfloor n/2\rfloor$. Since $N\ge q$, the middle edge
$\ell=e_N$ lies among $e_q,\ldots,e_{2N-q}$.
After round $q$, its remaining vertices would have indices
$q<k<n-q$, and there are no such integers.
The end edges $e_0,e_{2N}$ never serve as the common edge of these
two-triangle configurations, so their vertices are neither mutated nor deleted.
If $Z$ is the set of vertices removed by these deletion moves,
the resulting quiver is $Q_S\setminus Z$ by \eqref{eq:delete-undo}.
\end{proof}

\begin{figure}[htbp]
\centering
\begin{tikzpicture}[scale=1.15,transform shape,x=1cm,y=1cm,font=\small,>=stealth]
\begin{scope}[yshift=-0cm]
\draw[black] (0,0)--(8,0) (0,2)--(8,2);
\draw[gray!60] (0,0)--(2,2) (2,0)--(4,2) (4,0)--(6,2) (6,0)--(8,2)
    (2,0)--(2,2) (6,0)--(6,2);
\draw[Green,line width=1pt] (0,0)--(0,2) (8,0)--(8,2);
\draw[red,line width=1pt] (4,0)--(4,2);
\foreach \u/\v in {
    1/1, 1/2, 1/3, 1/4, 1/5, 1/6, 2/1, 2/2, 2/3, 2/4, 2/5, 2/6,
    3/1, 3/2, 3/3, 3/4, 3/5, 3/6, 4/1, 4/2, 4/3, 4/4, 4/5, 4/6,
    5/1, 5/2, 5/3, 5/4, 5/5, 5/6, 6/1, 6/2, 6/3, 6/4, 6/5, 6/6,
    7/1, 7/2, 7/3, 7/4, 7/5, 7/6, 8/1, 8/2, 8/3, 8/4, 8/5, 8/6,
    9/1, 9/2, 9/3, 9/4, 9/5, 9/6, 10/1, 10/2, 10/3, 10/4, 10/5, 10/6,
    11/1, 11/2, 11/3, 11/4, 11/5, 11/6, 12/1, 12/2, 12/3, 12/4, 12/5, 12/6,
    13/1, 13/2, 13/3, 13/4, 13/5, 13/6, 14/1, 14/2, 14/3, 14/4, 14/5, 14/6,
    15/1, 15/2, 15/3, 15/4, 15/5, 15/6, 16/1, 16/2, 16/3, 16/4, 16/5, 16/6,
    17/1, 17/2, 17/3, 17/4, 17/5, 17/6, 18/1, 18/2, 18/3, 18/4, 18/5, 18/6,
    19/1, 19/2, 19/3, 19/4, 19/5, 19/6, 20/1, 20/2, 20/3, 20/4, 20/5, 20/6,
    21/1, 21/2, 21/3, 21/4, 21/5, 21/6, 22/1, 22/2, 22/3, 22/4, 22/5, 22/6,
    23/1, 23/2, 23/3, 23/4, 23/5, 23/6, 24/1, 24/2, 24/3, 24/4, 24/5, 24/6,
    25/1, 25/2, 25/3, 25/4, 25/5, 25/6, 26/1, 26/2, 26/3, 26/4, 26/5, 26/6,
    27/1, 27/2, 27/3, 27/4, 27/5, 27/6
} {
    \fill ({2*\u/7},{2*\v/7}) circle (1.4pt);
}
\foreach \v in {1,...,6} {
    \fill[Green] (0,{2*\v/7}) circle (1.4pt);
    \fill[Green] (8,{2*\v/7}) circle (1.4pt);
}
\node[left,text=Green] at (0,1) {$e_0$};
\node[right,text=Green] at (8,1) {$e_8$};
\draw[gray!60] (2,2)--(2.2,2.2);
\node[above right,font=\scriptsize,inner sep=1pt] at (2.2,2.2) {$e_1$};
\node[above] at (1,2) {$\tau_1$};
\node[below] at (1,0) {$\tau_2$};
\draw[gray!60] (4,2)--(4.2,2.2);
\node[above right,font=\scriptsize,inner sep=1pt] at (4.2,2.2) {$e_3$};
\node[above] at (3,2) {$\tau_3$};
\node[below] at (3,0) {$\tau_4$};
\draw[gray!60] (6,2)--(6.2,2.2);
\node[above right,font=\scriptsize,inner sep=1pt] at (6.2,2.2) {$e_5$};
\node[above] at (5,2) {$\tau_5$};
\node[below] at (5,0) {$\tau_6$};
\draw[gray!60] (8,2)--(8.2,2.2);
\node[above right,font=\scriptsize,inner sep=1pt] at (8.2,2.2) {$e_7$};
\node[above] at (7,2) {$\tau_7$};
\node[below] at (7,0) {$\tau_8$};
\node[below,font=\scriptsize] at (2,0) {$e_2$};
\node[below,text=red,font=\scriptsize] at (4,0) {$e_4=\ell$};
\node[below,font=\scriptsize] at (6,0) {$e_6$};
\node[anchor=west,align=left] at (9.1,1) {Initial state\\$n=7$, $N=4$\\$C_1=\varnothing$};
\end{scope}
\draw[->] (9.3,-0)--node[right] {Lemma~\ref{lem:edge-layer}} (9.3,-0.8);
\begin{scope}[yshift=-3.1cm]
\draw[black] (0,0)--(8,0) (0,2)--(8,2);
\draw[gray!60] (0,0)--(2,2) (2,0)--(4,2) (4,0)--(6,2) (6,0)--(8,2)
    (2,0)--(2,2) (6,0)--(6,2);
\draw[Green,line width=1pt] (0,0)--(0,2) (8,0)--(8,2);
\draw[red,line width=1pt] (4,0)--(4,2);
\foreach \u/\v in {
    1/2, 1/3, 1/4, 1/5, 1/6, 2/1, 2/2, 2/3, 2/4, 2/5, 2/6, 3/1,
    3/2, 3/3, 3/4, 3/5, 3/6, 4/1, 4/2, 4/3, 4/4, 4/5, 4/6, 5/1,
    5/2, 5/3, 5/4, 5/5, 5/6, 6/1, 6/2, 6/3, 6/4, 6/5, 7/2, 7/3,
    7/4, 7/5, 8/2, 8/3, 8/4, 8/5, 8/6, 9/1, 9/2, 9/3, 9/4, 9/5,
    9/6, 10/1, 10/2, 10/3, 10/4, 10/5, 10/6, 11/1, 11/2, 11/3, 11/4, 11/5,
    11/6, 12/1, 12/2, 12/3, 12/4, 12/5, 12/6, 13/1, 13/2, 13/3, 13/4, 13/5,
    14/2, 14/3, 14/4, 14/5, 15/2, 15/3, 15/4, 15/5, 15/6, 16/1, 16/2, 16/3,
    16/4, 16/5, 16/6, 17/1, 17/2, 17/3, 17/4, 17/5, 17/6, 18/1, 18/2, 18/3,
    18/4, 18/5, 18/6, 19/1, 19/2, 19/3, 19/4, 19/5, 19/6, 20/1, 20/2, 20/3,
    20/4, 20/5, 21/2, 21/3, 21/4, 21/5, 22/2, 22/3, 22/4, 22/5, 22/6, 23/1,
    23/2, 23/3, 23/4, 23/5, 23/6, 24/1, 24/2, 24/3, 24/4, 24/5, 24/6, 25/1,
    25/2, 25/3, 25/4, 25/5, 25/6, 26/1, 26/2, 26/3, 26/4, 26/5, 26/6, 27/1,
    27/2, 27/3, 27/4, 27/5
} {
    \fill ({2*\u/7},{2*\v/7}) circle (1.4pt);
}
\foreach \v in {1,...,6} {
    \fill[Green] (0,{2*\v/7}) circle (1.4pt);
    \fill[Green] (8,{2*\v/7}) circle (1.4pt);
}
\foreach \u/\v in {
    1/1, 6/6, 7/1, 7/6, 8/1, 13/6, 14/1, 14/6, 15/1, 20/6, 21/1, 21/6,
    22/1, 27/6
} {
    \node[text=gray,font=\scriptsize,inner sep=0pt] at ({2*\u/7},{2*\v/7}) {$\times$};
}
\node[left,text=Green] at (0,1) {$e_0$};
\node[right,text=Green] at (8,1) {$e_8$};
\draw[gray!60] (2,2)--(2.2,2.2);
\node[above right,font=\scriptsize,inner sep=1pt] at (2.2,2.2) {$e_1$};
\node[above] at (1,2) {$\tau_1$};
\node[below] at (1,0) {$\tau_2$};
\draw[gray!60] (4,2)--(4.2,2.2);
\node[above right,font=\scriptsize,inner sep=1pt] at (4.2,2.2) {$e_3$};
\node[above] at (3,2) {$\tau_3$};
\node[below] at (3,0) {$\tau_4$};
\draw[gray!60] (6,2)--(6.2,2.2);
\node[above right,font=\scriptsize,inner sep=1pt] at (6.2,2.2) {$e_5$};
\node[above] at (5,2) {$\tau_5$};
\node[below] at (5,0) {$\tau_6$};
\draw[gray!60] (8,2)--(8.2,2.2);
\node[above right,font=\scriptsize,inner sep=1pt] at (8.2,2.2) {$e_7$};
\node[above] at (7,2) {$\tau_7$};
\node[below] at (7,0) {$\tau_8$};
\node[below,font=\scriptsize] at (2,0) {$e_2$};
\node[below,text=red,font=\scriptsize] at (4,0) {$e_4=\ell$};
\node[below,font=\scriptsize] at (6,0) {$e_6$};
\node[anchor=west,align=left] at (9.1,1) {After round $1$\\$e_1,\ldots,e_7$: depth $1$};
\end{scope}
\draw[->] (9.3,-3.1)--node[right] {Lemma~\ref{lem:triangle-shell}} (9.3,-3.9);
\begin{scope}[yshift=-6.2cm]
\draw[black] (0,0)--(8,0) (0,2)--(8,2);
\draw[gray!60] (0,0)--(2,2) (2,0)--(4,2) (4,0)--(6,2) (6,0)--(8,2)
    (2,0)--(2,2) (6,0)--(6,2);
\draw[Green,line width=1pt] (0,0)--(0,2) (8,0)--(8,2);
\draw[red,line width=1pt] (4,0)--(4,2);
\foreach \u/\v in {
    1/2, 1/3, 1/4, 1/5, 1/6, 2/2, 2/3, 2/4, 2/5, 2/6, 3/1, 3/2,
    3/3, 3/4, 3/5, 3/6, 4/1, 4/2, 4/3, 4/4, 4/5, 4/6, 5/1, 5/2,
    5/3, 5/4, 5/5, 5/6, 6/2, 6/3, 6/4, 6/5, 7/2, 7/3, 7/4, 7/5,
    8/2, 8/3, 8/4, 8/5, 9/2, 9/3, 9/4, 9/5, 9/6, 10/1, 10/2, 10/3,
    10/4, 10/5, 10/6, 11/1, 11/2, 11/3, 11/4, 11/5, 11/6, 12/1, 12/2, 12/3,
    12/4, 12/5, 13/2, 13/3, 13/4, 13/5, 14/2, 14/3, 14/4, 14/5, 15/2, 15/3,
    15/4, 15/5, 16/2, 16/3, 16/4, 16/5, 16/6, 17/1, 17/2, 17/3, 17/4, 17/5,
    17/6, 18/1, 18/2, 18/3, 18/4, 18/5, 18/6, 19/1, 19/2, 19/3, 19/4, 19/5,
    20/2, 20/3, 20/4, 20/5, 21/2, 21/3, 21/4, 21/5, 22/2, 22/3, 22/4, 22/5,
    23/1, 23/2, 23/3, 23/4, 23/5, 23/6, 24/1, 24/2, 24/3, 24/4, 24/5, 24/6,
    25/1, 25/2, 25/3, 25/4, 25/5, 25/6, 26/1, 26/2, 26/3, 26/4, 26/5, 27/1,
    27/2, 27/3, 27/4, 27/5
} {
    \fill ({2*\u/7},{2*\v/7}) circle (1.4pt);
}
\foreach \v in {1,...,6} {
    \fill[Green] (0,{2*\v/7}) circle (1.4pt);
    \fill[Green] (8,{2*\v/7}) circle (1.4pt);
}
\foreach \u/\v in {
    1/1, 2/1, 6/1, 6/6, 7/1, 7/6, 8/1, 8/6, 9/1, 12/6, 13/1, 13/6,
    14/1, 14/6, 15/1, 15/6, 16/1, 19/6, 20/1, 20/6, 21/1, 21/6, 22/1, 22/6,
    26/6, 27/6
} {
    \node[text=gray,font=\scriptsize,inner sep=0pt] at ({2*\u/7},{2*\v/7}) {$\times$};
}
\node[left,text=Green] at (0,1) {$e_0$};
\node[right,text=Green] at (8,1) {$e_8$};
\draw[gray!60] (2,2)--(2.2,2.2);
\node[above right,font=\scriptsize,inner sep=1pt] at (2.2,2.2) {$e_1$};
\node[above] at (1,2) {$\tau_1$};
\node[below] at (1,0) {$\tau_2$};
\draw[gray!60] (4,2)--(4.2,2.2);
\node[above right,font=\scriptsize,inner sep=1pt] at (4.2,2.2) {$e_3$};
\node[above] at (3,2) {$\tau_3$};
\node[below] at (3,0) {$\tau_4$};
\draw[gray!60] (6,2)--(6.2,2.2);
\node[above right,font=\scriptsize,inner sep=1pt] at (6.2,2.2) {$e_5$};
\node[above] at (5,2) {$\tau_5$};
\node[below] at (5,0) {$\tau_6$};
\draw[gray!60] (8,2)--(8.2,2.2);
\node[above right,font=\scriptsize,inner sep=1pt] at (8.2,2.2) {$e_7$};
\node[above] at (7,2) {$\tau_7$};
\node[below] at (7,0) {$\tau_8$};
\node[below,font=\scriptsize] at (2,0) {$e_2$};
\node[below,text=red,font=\scriptsize] at (4,0) {$e_4=\ell$};
\node[below,font=\scriptsize] at (6,0) {$e_6$};
\node[anchor=west,align=left] at (9.1,1) {Round $2$\\Triangle deletion moves\\$C_2$ removed in $\tau_2,\ldots,\tau_7$};
\end{scope}
\draw[->] (9.3,-6.2)--node[right] {Lemma~\ref{lem:edge-layer}} (9.3,-7);
\begin{scope}[yshift=-9.3cm]
\draw[black] (0,0)--(8,0) (0,2)--(8,2);
\draw[gray!60] (0,0)--(2,2) (2,0)--(4,2) (4,0)--(6,2) (6,0)--(8,2)
    (2,0)--(2,2) (6,0)--(6,2);
\draw[Green,line width=1pt] (0,0)--(0,2) (8,0)--(8,2);
\draw[red,line width=1pt] (4,0)--(4,2);
\foreach \u/\v in {
    1/2, 1/3, 1/4, 1/5, 1/6, 2/2, 2/3, 2/4, 2/5, 2/6, 3/1, 3/2,
    3/3, 3/4, 3/5, 3/6, 4/1, 4/2, 4/3, 4/4, 4/5, 4/6, 5/1, 5/2,
    5/3, 5/4, 5/5, 5/6, 6/2, 6/3, 6/4, 6/5, 7/3, 7/4, 8/2, 8/3,
    8/4, 8/5, 9/3, 9/4, 9/5, 9/6, 10/1, 10/2, 10/3, 10/4, 10/5, 10/6,
    11/1, 11/2, 11/3, 11/4, 11/5, 11/6, 12/1, 12/2, 12/3, 12/4, 13/2, 13/3,
    13/4, 13/5, 14/3, 14/4, 15/2, 15/3, 15/4, 15/5, 16/3, 16/4, 16/5, 16/6,
    17/1, 17/2, 17/3, 17/4, 17/5, 17/6, 18/1, 18/2, 18/3, 18/4, 18/5, 18/6,
    19/1, 19/2, 19/3, 19/4, 20/2, 20/3, 20/4, 20/5, 21/3, 21/4, 22/2, 22/3,
    22/4, 22/5, 23/1, 23/2, 23/3, 23/4, 23/5, 23/6, 24/1, 24/2, 24/3, 24/4,
    24/5, 24/6, 25/1, 25/2, 25/3, 25/4, 25/5, 25/6, 26/1, 26/2, 26/3, 26/4,
    26/5, 27/1, 27/2, 27/3, 27/4, 27/5
} {
    \fill ({2*\u/7},{2*\v/7}) circle (1.4pt);
}
\foreach \v in {1,...,6} {
    \fill[Green] (0,{2*\v/7}) circle (1.4pt);
    \fill[Green] (8,{2*\v/7}) circle (1.4pt);
}
\foreach \u/\v in {
    1/1, 2/1, 6/1, 6/6, 7/1, 7/2, 7/5, 7/6, 8/1, 8/6, 9/1, 9/2,
    12/5, 12/6, 13/1, 13/6, 14/1, 14/2, 14/5, 14/6, 15/1, 15/6, 16/1, 16/2,
    19/5, 19/6, 20/1, 20/6, 21/1, 21/2, 21/5, 21/6, 22/1, 22/6, 26/6, 27/6
} {
    \node[text=gray,font=\scriptsize,inner sep=0pt] at ({2*\u/7},{2*\v/7}) {$\times$};
}
\node[left,text=Green] at (0,1) {$e_0$};
\node[right,text=Green] at (8,1) {$e_8$};
\draw[gray!60] (2,2)--(2.2,2.2);
\node[above right,font=\scriptsize,inner sep=1pt] at (2.2,2.2) {$e_1$};
\node[above] at (1,2) {$\tau_1$};
\node[below] at (1,0) {$\tau_2$};
\draw[gray!60] (4,2)--(4.2,2.2);
\node[above right,font=\scriptsize,inner sep=1pt] at (4.2,2.2) {$e_3$};
\node[above] at (3,2) {$\tau_3$};
\node[below] at (3,0) {$\tau_4$};
\draw[gray!60] (6,2)--(6.2,2.2);
\node[above right,font=\scriptsize,inner sep=1pt] at (6.2,2.2) {$e_5$};
\node[above] at (5,2) {$\tau_5$};
\node[below] at (5,0) {$\tau_6$};
\draw[gray!60] (8,2)--(8.2,2.2);
\node[above right,font=\scriptsize,inner sep=1pt] at (8.2,2.2) {$e_7$};
\node[above] at (7,2) {$\tau_7$};
\node[below] at (7,0) {$\tau_8$};
\node[below,font=\scriptsize] at (2,0) {$e_2$};
\node[below,text=red,font=\scriptsize] at (4,0) {$e_4=\ell$};
\node[below,font=\scriptsize] at (6,0) {$e_6$};
\node[anchor=west,align=left] at (9.1,1) {After round $2$\\$e_2,\ldots,e_6$: depth $2$};
\end{scope}
\draw[->] (9.3,-9.3)--node[right] {Lemma~\ref{lem:triangle-shell}} (9.3,-10.1);
\begin{scope}[yshift=-12.4cm]
\draw[black] (0,0)--(8,0) (0,2)--(8,2);
\draw[gray!60] (0,0)--(2,2) (2,0)--(4,2) (4,0)--(6,2) (6,0)--(8,2)
    (2,0)--(2,2) (6,0)--(6,2);
\draw[Green,line width=1pt] (0,0)--(0,2) (8,0)--(8,2);
\draw[red,line width=1pt] (4,0)--(4,2);
\foreach \u/\v in {
    1/2, 1/3, 1/4, 1/5, 1/6, 2/2, 2/3, 2/4, 2/5, 2/6, 3/1, 3/2,
    3/3, 3/4, 3/5, 3/6, 4/1, 4/2, 4/3, 4/4, 4/5, 4/6, 5/1, 5/2,
    5/3, 5/4, 5/5, 5/6, 6/2, 6/3, 6/4, 6/5, 7/3, 7/4, 8/2, 8/3,
    8/4, 9/3, 9/4, 9/5, 10/3, 10/4, 10/5, 10/6, 11/1, 11/2, 11/3, 11/4,
    12/2, 12/3, 12/4, 13/3, 13/4, 13/5, 14/3, 14/4, 15/2, 15/3, 15/4, 16/3,
    16/4, 16/5, 17/3, 17/4, 17/5, 17/6, 18/1, 18/2, 18/3, 18/4, 19/2, 19/3,
    19/4, 20/3, 20/4, 20/5, 21/3, 21/4, 22/2, 22/3, 22/4, 22/5, 23/1, 23/2,
    23/3, 23/4, 23/5, 23/6, 24/1, 24/2, 24/3, 24/4, 24/5, 24/6, 25/1, 25/2,
    25/3, 25/4, 25/5, 25/6, 26/1, 26/2, 26/3, 26/4, 26/5, 27/1, 27/2, 27/3,
    27/4, 27/5
} {
    \fill ({2*\u/7},{2*\v/7}) circle (1.4pt);
}
\foreach \v in {1,...,6} {
    \fill[Green] (0,{2*\v/7}) circle (1.4pt);
    \fill[Green] (8,{2*\v/7}) circle (1.4pt);
}
\foreach \u/\v in {
    1/1, 2/1, 6/1, 6/6, 7/1, 7/2, 7/5, 7/6, 8/1, 8/5, 8/6, 9/1,
    9/2, 9/6, 10/1, 10/2, 11/5, 11/6, 12/1, 12/5, 12/6, 13/1, 13/2, 13/6,
    14/1, 14/2, 14/5, 14/6, 15/1, 15/5, 15/6, 16/1, 16/2, 16/6, 17/1, 17/2,
    18/5, 18/6, 19/1, 19/5, 19/6, 20/1, 20/2, 20/6, 21/1, 21/2, 21/5, 21/6,
    22/1, 22/6, 26/6, 27/6
} {
    \node[text=gray,font=\scriptsize,inner sep=0pt] at ({2*\u/7},{2*\v/7}) {$\times$};
}
\node[left,text=Green] at (0,1) {$e_0$};
\node[right,text=Green] at (8,1) {$e_8$};
\draw[gray!60] (2,2)--(2.2,2.2);
\node[above right,font=\scriptsize,inner sep=1pt] at (2.2,2.2) {$e_1$};
\node[above] at (1,2) {$\tau_1$};
\node[below] at (1,0) {$\tau_2$};
\draw[gray!60] (4,2)--(4.2,2.2);
\node[above right,font=\scriptsize,inner sep=1pt] at (4.2,2.2) {$e_3$};
\node[above] at (3,2) {$\tau_3$};
\node[below] at (3,0) {$\tau_4$};
\draw[gray!60] (6,2)--(6.2,2.2);
\node[above right,font=\scriptsize,inner sep=1pt] at (6.2,2.2) {$e_5$};
\node[above] at (5,2) {$\tau_5$};
\node[below] at (5,0) {$\tau_6$};
\draw[gray!60] (8,2)--(8.2,2.2);
\node[above right,font=\scriptsize,inner sep=1pt] at (8.2,2.2) {$e_7$};
\node[above] at (7,2) {$\tau_7$};
\node[below] at (7,0) {$\tau_8$};
\node[below,font=\scriptsize] at (2,0) {$e_2$};
\node[below,text=red,font=\scriptsize] at (4,0) {$e_4=\ell$};
\node[below,font=\scriptsize] at (6,0) {$e_6$};
\node[anchor=west,align=left] at (9.1,1) {Round $3$\\Triangle deletion moves\\$C_3$ removed in $\tau_3,\ldots,\tau_6$};
\end{scope}
\draw[->] (9.3,-12.4)--node[right] {Lemma~\ref{lem:edge-layer}} (9.3,-13.2);
\begin{scope}[yshift=-15.5cm]
\draw[black] (0,0)--(8,0) (0,2)--(8,2);
\draw[gray!60] (0,0)--(2,2) (2,0)--(4,2) (4,0)--(6,2) (6,0)--(8,2)
    (2,0)--(2,2) (6,0)--(6,2);
\draw[Green,line width=1pt] (0,0)--(0,2) (8,0)--(8,2);
\draw[red,line width=1pt] (4,0)--(4,2);
\foreach \u/\v in {
    1/2, 1/3, 1/4, 1/5, 1/6, 2/2, 2/3, 2/4, 2/5, 2/6, 3/1, 3/2,
    3/3, 3/4, 3/5, 3/6, 4/1, 4/2, 4/3, 4/4, 4/5, 4/6, 5/1, 5/2,
    5/3, 5/4, 5/5, 5/6, 6/2, 6/3, 6/4, 6/5, 7/3, 7/4, 8/2, 8/3,
    8/4, 9/3, 9/4, 9/5, 10/4, 10/5, 10/6, 11/1, 11/2, 11/3, 12/2, 12/3,
    12/4, 13/3, 13/4, 13/5, 15/2, 15/3, 15/4, 16/3, 16/4, 16/5, 17/4, 17/5,
    17/6, 18/1, 18/2, 18/3, 19/2, 19/3, 19/4, 20/3, 20/4, 20/5, 21/3, 21/4,
    22/2, 22/3, 22/4, 22/5, 23/1, 23/2, 23/3, 23/4, 23/5, 23/6, 24/1, 24/2,
    24/3, 24/4, 24/5, 24/6, 25/1, 25/2, 25/3, 25/4, 25/5, 25/6, 26/1, 26/2,
    26/3, 26/4, 26/5, 27/1, 27/2, 27/3, 27/4, 27/5
} {
    \fill ({2*\u/7},{2*\v/7}) circle (1.4pt);
}
\foreach \v in {1,...,6} {
    \fill[Green] (0,{2*\v/7}) circle (1.4pt);
    \fill[Green] (8,{2*\v/7}) circle (1.4pt);
}
\foreach \u/\v in {
    1/1, 2/1, 6/1, 6/6, 7/1, 7/2, 7/5, 7/6, 8/1, 8/5, 8/6, 9/1,
    9/2, 9/6, 10/1, 10/2, 10/3, 11/4, 11/5, 11/6, 12/1, 12/5, 12/6, 13/1,
    13/2, 13/6, 14/1, 14/2, 14/3, 14/4, 14/5, 14/6, 15/1, 15/5, 15/6, 16/1,
    16/2, 16/6, 17/1, 17/2, 17/3, 18/4, 18/5, 18/6, 19/1, 19/5, 19/6, 20/1,
    20/2, 20/6, 21/1, 21/2, 21/5, 21/6, 22/1, 22/6, 26/6, 27/6
} {
    \node[text=gray,font=\scriptsize,inner sep=0pt] at ({2*\u/7},{2*\v/7}) {$\times$};
}
\node[left,text=Green] at (0,1) {$e_0$};
\node[right,text=Green] at (8,1) {$e_8$};
\draw[gray!60] (2,2)--(2.2,2.2);
\node[above right,font=\scriptsize,inner sep=1pt] at (2.2,2.2) {$e_1$};
\node[above] at (1,2) {$\tau_1$};
\node[below] at (1,0) {$\tau_2$};
\draw[gray!60] (4,2)--(4.2,2.2);
\node[above right,font=\scriptsize,inner sep=1pt] at (4.2,2.2) {$e_3$};
\node[above] at (3,2) {$\tau_3$};
\node[below] at (3,0) {$\tau_4$};
\draw[gray!60] (6,2)--(6.2,2.2);
\node[above right,font=\scriptsize,inner sep=1pt] at (6.2,2.2) {$e_5$};
\node[above] at (5,2) {$\tau_5$};
\node[below] at (5,0) {$\tau_6$};
\draw[gray!60] (8,2)--(8.2,2.2);
\node[above right,font=\scriptsize,inner sep=1pt] at (8.2,2.2) {$e_7$};
\node[above] at (7,2) {$\tau_7$};
\node[below] at (7,0) {$\tau_8$};
\node[below,font=\scriptsize] at (2,0) {$e_2$};
\node[below,text=red,font=\scriptsize] at (4,0) {$e_4=\ell$};
\node[below,font=\scriptsize] at (6,0) {$e_6$};
\node[anchor=west,align=left] at (9.1,1) {After round $3$\\$e_3,e_4,e_5$: depth $3$\\No vertices remain on $\ell$};
\end{scope}
\end{tikzpicture}
\caption{Proposition~\ref{prop:strip} for $n=7$ and $N=4$.
Crosses are the deleted vertices. No vertex on the green end edges is mutated or deleted.}
\label{fig:buffer}
\end{figure}
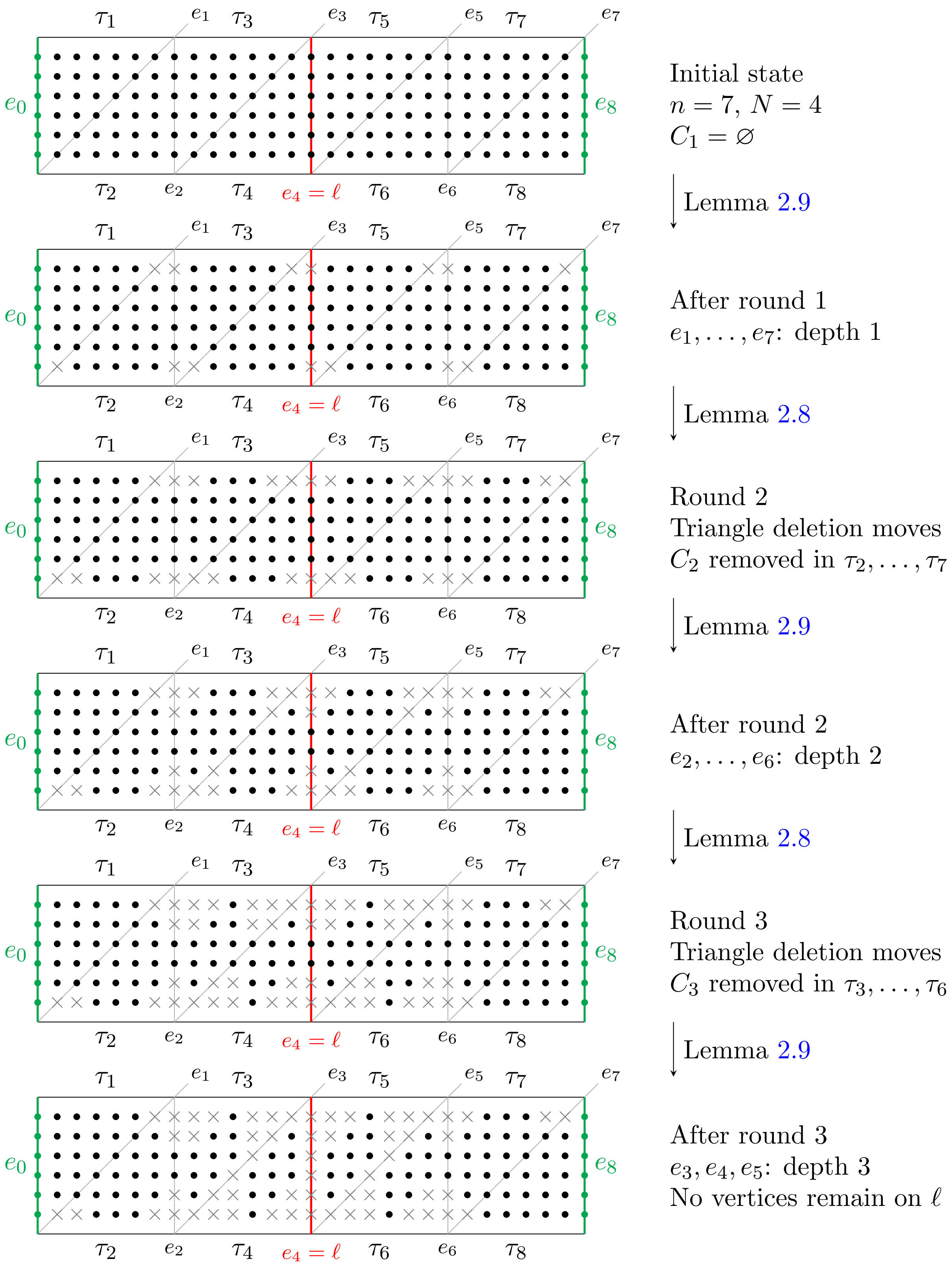

Figure~\ref{fig:buffer} illustrates the proof for $n=7$ and $N=4$.

\clearpage

\section{When localizations are localizations}
\label{sec:localization}

Let $A$ be a (not derived) commutative ring, and let $\Mod_A$ be its unbounded
derived category.  All functors below are derived.

A full stable subcategory is called {\em localizing} if it is closed under coproducts; such subcategories are characterized as the kernels of {\em localizations}, i.e. maps $q:\Mod_A\to\mathcal C$, with fully faithful right
adjoint $r$. 
Given such a localization, we will say an object of $\Mod_A$ is {\em $q$-local} if it belongs to the image of $r$.
A localization is called {\em  smashing} if $r$ preserves coproducts, or equivalently if $q$ preserves compact objects.

For Noetherian $A$, Neeman shows that localizing subcategories of $\Mod_A$ are in bijection with
subsets of $\Spec A$ \cite{Neeman}; to a subset $Z \subset \Spec A$ is associated the subcategory generated by the
residue fields $k(\mathfrak p)$ for $\mathfrak p\in Z$.
Neeman further shows that said category is smashing exactly when $Z$ is specialization-closed, i.e. a union of closed subsets, in which case its objects are those modules supported on $Z$.  

For possibly non-Noetherian $A$, Thomason shows that localizing subcategories generated by compact objects (the localizations for which are necessarily smashing) in $\Mod_A$ are similarly in bijection with subsets $Z\subset\Spec A$ which are unions of closed subsets with quasi-compact complement \cite[Theorem~3.15]{Thomason}.  For Noetherian $A$, this recovers Neeman's result.  In the Noetherian case, all smashing localizations are given as quotients by compactly generated subcategories \cite[Corollary~3.4]{Neeman}; in general this is false \cite[Section~2, Proposition]{Keller-smashing}.
(We do not need anything non-Noetherian in this article, but we wish to be clear about when this hypothesis becomes necessary.)

Fix such a $Z$ and let $q: \Mod_A \to \mathcal{C}$ be the corresponding smashing localization by a compactly generated subcategory.  We write $R := \Hom_{\mathcal{C}}(q(A), q(A))$.  Since $q(A)$ is evidently a compact generator of $\mathcal{C}$, we have $\mathcal{C} \simeq \Mod_R$.
Under this identification, 
$q$ is extension of scalars and $r$ is
restriction of scalars. The counit $qr(R)\to R$ is an equivalence,
since $r$ is fully faithful. It is the multiplication
$R\otimes_A R\to R$.

\begin{lemma}\label{lem:smashing-unit}
The algebra $R$ is a commutative derived $A$-algebra.
\end{lemma}
\begin{proof}
Every localizing subcategory of $\Mod_A$ is a tensor ideal.
Indeed, for an object $N$ in that subcategory, the objects $M$ with
$M\otimes_A N$ in the subcategory form a localizing subcategory containing $A$.
The quotient therefore inherits a symmetric monoidal structure with
unit $q(A)$, whose endomorphism algebra $R$ is therefore commutative.
\end{proof}

\begin{lemma}\label{lem:smashing-open-colimit}
Let $\mathcal F_Z$ be the poset of closed subsets $F\subset Z$
such that $\Spec A\setminus F$ is quasi-compact, ordered by inclusion. Then for $M \in \Mod_A$,
\[
 R\otimes_A M\simeq\operatorname*{colim}_{F\in\mathcal F_Z}
       R\Gamma(\Spec A\setminus F,\widetilde M).
\]
When $M = A$, the equivalence is one of derived commutative algebras. 
\end{lemma}
\begin{proof}
Write $\mathcal K_Z=\ker q$ and
$C_F=R\Gamma(\Spec A\setminus F,\OO)$.  
For $F=V(f_1,\ldots,f_n)$, use the affine open cover
$\Spec A\setminus F=\bigcup_{i=1}^n D(f_i)$.
Its intersections are $D(f_{i_0}\cdots f_{i_p})$, so
\[
 C_F\simeq\left[
 \bigoplus_i A[f_i^{-1}]\longrightarrow
 \bigoplus_{i<j}A[(f_if_j)^{-1}]\longrightarrow\cdots
 \longrightarrow A[(f_1\cdots f_n)^{-1}]
 \right],
\]
in degrees $0,\ldots,n-1$, with the alternating restriction maps.
The terms are flat over $A$, hence the non-derived tensor product
$C_F\otimes_A M$ computes  $R\Gamma(\Spec A\setminus F,\widetilde M)$.

The fiber of $A\to C_F$ is the complex
\[
 \left[A\longrightarrow\bigoplus_i A[f_i^{-1}]
 \longrightarrow\bigoplus_{i<j}A[(f_if_j)^{-1}]
 \longrightarrow\cdots\right]
\]
in degrees $0,\ldots,n$.
If $K_A(f_1^r,\ldots,f_n^r)$ denotes the Koszul complex in degrees
$-n,\ldots,0$, then
\[
 \fib(M\to C_F\otimes_A M)
 \simeq\operatorname*{colim}_{r\ge1}
       K_A(f_1^r,\ldots,f_n^r)^\vee\otimes_A M
\]
by \cite[Section~36.9]{Stacks-section23}.
Each dual Koszul complex is perfect and supported on $F$.
Since $\mathcal K_Z$ is a localizing tensor ideal, this fiber belongs to $\mathcal K_Z$.
The same holds after taking the filtered colimit over $F$.

If $P$ is perfect and $\Supp P\subset Z$, the sets
$F\supset\Supp P$ are cofinal and satisfy
\[
 \Hom_A(P,C_F\otimes_A M)
 \simeq R\Gamma(\Spec A\setminus F,
                    \widetilde{P^\vee\otimes_A M})=0.
\]
Compactness of $P$ passes this vanishing to the colimit.
Thomason's classification \cite[Theorem~3.15]{Thomason} says that
these $P$ generate $\mathcal K_Z$.
Hence the colimit is $q$-local and its map from $M$ has fiber in
$\mathcal K_Z$, which identifies it with $rq(M)$.
For $M=A$, the restriction maps are maps of commutative algebras,
and this identification preserves their units and multiplication.
\end{proof}

\begin{corollary} \label{cor: B truncation}
    $H^j(R)=0$ for $j<0$, so truncation gives a natural map $H^0(R) \to R$.
\end{corollary}
\begin{corollary}  \label{cor: frac field embedding}
    If $A$ is  integral and $\mathcal{C} \ne 0$,  
    one has canonically 
    $A \subset H^0(R) \subset \mathrm{Frac}(A)$.
\end{corollary}
\begin{proof}
Let $Z$ be the specialization closed subset giving the support of $\ker q$. It cannot contain the generic point, since then
$Z=\Spec A$ and $q=0$.
Lemma~\ref{lem:smashing-open-colimit} therefore gives
$H^0(R)=\operatorname*{colim}_{F\in\mathcal F_Z}
       H^0(\Spec A\setminus F,\OO)\subset \operatorname{Frac}(A)$.
The map $A \to H^0(R)$ then gives the other inclusion.
\end{proof}

\begin{remark}
    If one only assumes $q$ is smashing and not that $\ker q$ is compactly generated, then Corollary \ref{cor: frac field embedding} is false, and Keller's example \cite[Section~2, Proposition]{Keller-smashing} gives a counterexample.
\end{remark}

\begin{corollary} \label{lem:discrete-smashing-flat}
If $H^0(R) \simeq R$, then the map of (not derived) commutative rings $A\to R$ is a flat ring epimorphism, and $\operatorname{image}(\Spec R\to\Spec A)=\Spec A\setminus Z$.
\end{corollary}
\begin{proof}
    For every $A$-module $M$, Lemma~\ref{lem:smashing-open-colimit} gives
$R\otimes_A M\in\Mod_A^{\ge0}$.
Thus $\operatorname{Tor}^A_i(R,M)=H^{-i}(R\otimes_A M)=0$ for $i>0$,
so $R$ is flat.  The fact that multiplication $R \otimes_A R \to R$ is an equivalence is a classical criterion for a ring morphism to be an epimorphism.
For $\mathfrak p\in\Spec A$, the same \v{C}ech complexes give
\[
 R\otimes_A k(\mathfrak p)\cong
 \begin{cases}
  0,&\mathfrak p\in Z,\\
  k(\mathfrak p),&\mathfrak p\notin Z.
 \end{cases}
\]
This proves the assertion about the image on spectra.
\end{proof}
\begin{remark}
For a flat epimorphism of (not derived) commutative rings $A\to T$ and a multiplicative subset $S\subset A$,
\[
 T\cong S^{-1}A\ \text{as $A$-algebras}
 \quad\Longleftrightarrow\quad
 \operatorname{image}(\Spec T\to\Spec A)
       =\{\mathfrak p:\mathfrak p\cap S=\varnothing\}.
\]
The forward implication is immediate.
For the converse, each $s\in S$ avoids every prime of $T$, so is a unit.
The induced map $S^{-1}A\to T$ is a flat epimorphism and is surjective
on spectra. It is therefore faithfully flat, hence an isomorphism
\cite[Tag~04VU]{Stacks-section23}.
\end{remark}

\begin{lemma}\label{lem: local implies affine}
The map $H^0(R)\to R$ is an equivalence if and only if $H^0(R)$ is $q$-local.
\end{lemma}
\begin{proof}
If said map is an equivalence, then $H^0(R)$ is $q$-local since $R$ is.  For the converse,
because $H^{< 0}R = 0$, the map $A \to R$ factors as $A \to H^{0}(R) \to R$.  If $H^{0}(R)$ is local, then we have $\Hom_A(A, H^0(R)) = \Hom_R(R, H^0(R))$; taking the image of the above map gives a map $R \to H^0(R)$.  It is an exercise with adjunction to see that it is inverse to $H^0(R) \to R$.
\end{proof}

\begin{lemma}\label{lem: perfect implies affine}
    If $R$ is perfect over $H^0(R)$, then the map $H^0(R) \to R$ is an equivalence.
\end{lemma}
\begin{proof}
Evidently the support of $R$ as an $H^0(R)$ module is all of $\Spec \, H^0(R)$.
Thomason's support theorem \cite[Lemma~3.14]{Thomason} shows that $H^0(R)$ is in the thick subcategory of $\Perf_{H^0(R)}$ generated by $R$.
Restricting scalars, $H^0(R)$ is in the thick subcategory of $\Mod_A$ generated by $R$, so $H^0(R)$ is local.  The result now follows by Lemma \ref{lem: local implies affine}.
\end{proof}

\begin{example}\label{ex:punctured-affine-localization}
Take
$A=\C[u,v]$.
The localization $\Mod_A\to\QCoh(\bA^2\setminus\{0\})$ preserves compact objects,
and its unit has endomorphism algebra $R:=R\Gamma(\bA^2\setminus\{0\},\OO)$.
The \v{C}ech complex for the cover by $\bA^1 \times (\bA^1 \setminus 0)$ and $(\bA^1 \setminus 0) \times \bA^1$ gives
\[
 H^0(R)=A,\qquad
 H^1(R)=A[(uv)^{-1}]/(A[u^{-1}]+A[v^{-1}])\ne0.
\]
The classes $u^{-a}v^{-b}$, $a,b\ge1$, show that $H^1(R)$ is not
a finite $A$-module. 
\end{example}

We will use the following criterion: 

\begin{theorem}\label{thm:recognition}
Let $k$ be a commutative ring and let $A$ be a Noetherian commutative $k$-algebra
of finite Krull dimension.
Let $q:\Mod_A\to\mathcal C$ be a localization preserving compact
objects, and set $R=\End_{\mathcal C}(q(A))$.
Suppose $\HH^0(\mathcal C^c/k)$ is Noetherian of finite global dimension
and $\HH^j(\mathcal C^c/k)$ is finite over it for every $j$.
Then the natural map $H^0(R)\to R$ is an equivalence.
\end{theorem}
\begin{proof}
Since $A$ is Noetherian, the kernel of $q$ is compactly generated
\cite[Corollary~3.4]{Neeman}.
The preceding results therefore apply to $R$.
By \cite[Corollaries~8.1--8.2]{Toen}, Hochschild cohomology is the
cohomology of the derived center
\[
 \mathfrak Z(\mathcal C^c)
 :=\RHom_{\operatorname{Fun}^{\mathrm{ex}}_k(\mathcal C^c,\mathcal C^c)}
           (\id,\id).
\]
The tensor unit of $\mathcal C^c\simeq\Perf(R)$ is $q(A)$.
Tensoring an endomorphism of $q(A)$ with each object defines
a natural endomorphism of the identity functor.
Evaluation at $q(A)$ gives a retraction
\[
 R\xrightarrow{s}\mathfrak Z(\mathcal C^c)
   \xrightarrow{e}R,\qquad
 s(a)_M=a\otimes\id_M,\quad e(\eta)=\eta_{q(A)},\quad es=\id_R.
\]
Both maps are $R$-linear for the action defined by $s$.

To identify the degree-zero rings, use $R\otimes_A R\simeq R$.
Extension and restriction of scalars give
\[
 \begin{aligned}
 (R\otimes_k R)\otimes_{A\otimes_k A}A&\simeq R,\\
 \mathfrak Z(\mathcal C^c)
 &\simeq\RHom_{R\otimes_k R}(R,R)\\
 &\simeq\RHom_{A\otimes_k A}(A,R)\\
 &\simeq\RHom_A(A\otimes_{A\otimes_k A}A,R).
 \end{aligned}
\]
The derived tensor product $A\otimes_{A\otimes_k A}A$ lies in
$\Mod_A^{\le0}$ and has degree-zero cohomology $A$.
Since $R\in\Mod_A^{\ge0}$, evaluation induces
\[
 \HH^0(\mathcal C^c/k)
 \cong\Hom_A(A,H^0(R))\cong H^0(R).
\]
Consequently every $H^j(R)$ is a direct summand of
$\HH^j(\mathcal C^c/k)$ as an $H^0(R)$-module, and is finite.

Lemma~\ref{lem:smashing-open-colimit} and Grothendieck vanishing give $H^j(R)=0$ for $j>\dim A$.
Since $H^0(R)$ has finite global dimension, the bounded complex $R$
with finite cohomology is perfect over $H^0(R)$.
Lemma~\ref{lem: perfect implies affine} gives the result.
\end{proof}

\section{Localization  of affine mirror equivalences}
\label{sec:affine-freezing}

The following is a translation of the general facts from Section \ref{sec:localization} about when categorical localizations are geometric localizations to the setting of wrapped Fukaya categories, for which we use $\C$-coefficients and integer gradings:

\begin{theorem}\label{thm:affine-freezing}
Let $D_0\subset D$ be a Liouville subdomain inclusion, with
 Weinstein completions $X_0$ and $X$.
Fix compatible grading/orientation data.  Consider the Viterbo functor 
$\mathrm{Vit}:\Fuk(X)\to\Fuk(X_0)$.

Suppose given a smooth finitely generated integral commutative $\C$-algebra $A$,
and an equivalence $\Phi:\Perf(A)\xrightarrow{\sim}\Fuk(X)$. 
We denote $E:= \Phi(A)$ and $E_0 = \mathrm{Vit}(E)$.

Suppose there is a multiplicative subset $S\subset A\setminus\{0\}$, such that: 
\begin{enumerate}[label=\textup{(H\arabic*)}]
\item\label{hyp:restriction-rings}
There is a nonempty family of subrings  $R_\alpha$ with $S^{-1} A\subset R_\alpha\subset \mathrm{Frac}(A)$ such that $\bigcap_\alpha R_\alpha=S^{-1}A$,
and  functors $\rho_\alpha:\Fuk(X_0)\to\Perf(R_\alpha)$, such that
\begin{equation}\label{eq:restriction-input}
 \rho_\alpha\circ\mathrm{Vit} \circ \Phi
 \simeq R_\alpha\otimes_A-
\end{equation}
\item\label{hyp:denominator-units}
For every $s\in S$, the endomorphism
\[
 \mathrm{Vit}\bigl(\Phi(s:A\to A)\bigr):E_0\to E_0
\]
is invertible in $H^0\Fuk(X_0)$.
\item\label{hyp:restriction-finite}
For every $j\in\Z$, $\SH^j(X_0)$ is finite over $\SH^0(X_0)$,
in the specified grading.
\end{enumerate}
Then there is an equivalence $\Phi_0:\Perf(S^{-1}A)\xrightarrow{\sim}\Fuk(X_0)$
with $\Phi_0(S^{-1}A )=E_0$, making the square
\begin{equation}\label{eq:affine-freezing-square}
 \begin{tikzcd}[column sep=large]
 \Perf(A)\arrow[r,"\Phi","\sim"']
       \arrow[d,"S^{-1}A\otimes_A-"']
    &\Fuk(X)\arrow[d,"\mathrm{Vit}"]\\
 \Perf(S^{-1}A)\arrow[r,"\Phi_0","\sim"']&\Fuk(X_0)
 \end{tikzcd}
\end{equation}
commute. It also satisfies
\begin{equation}\label{eq:restriction-output}
 \rho_\alpha \circ \Phi_0 \simeq R_\alpha\otimes_{S^{-1} A}-
\end{equation}
\end{theorem}

\begin{proof} 
By either \cite{GPS} or \cite{Sylvan}, the Viterbo functor is a localization (without any ind-completion), so defines a functor preserving compact objects: 
\[
 q=\Ind(\mathrm{Vit})\circ\Ind\Phi:
 \Mod_A\longrightarrow\Ind\Fuk(X_0),\qquad q(A)=E_0.
\]

We take $R = \End(E_0)$, and are in the setting of Section \ref{sec:localization}.

 Corollary \ref{cor: frac field embedding} gives $A \subset H^0(R) \subset \mathrm{Frac}(A)$.  By \ref{hyp:restriction-rings}, we get $A \subset H^0(R) \subset \bigcap R_\alpha = S^{-1} A$.  Then \ref{hyp:denominator-units} shows $H^0(R) = S^{-1} A$.
Since $A$ is smooth of finite type, $S^{-1} A$ is Noetherian of finite
global dimension.

The closed--open ring homomorphism
$\SH^*(X_0)\cong\HH^*(\Fuk(X_0))$ \cite{Ganatra} is an isomorphism for Weinstein manifolds (the nondegeneracy hypothesis of
\cite[Theorem~1.1]{Ganatra} having been verified by \cite{CDRGG} or \cite{GPS}). 
Thus hypothesis~\ref{hyp:restriction-finite} and
Theorem \ref{thm:recognition} give $R \simeq H^0(R) = S^{-1} A$.  This constructs $\Phi_0$. The remaining assertions are immediate.
\end{proof}

\begin{remark}
    If $R$ is smooth and finitely generated,
$\Perf(R)=\Perf(\Spec R)$, but in the above argument we do not assume finite generation of $S^{-1}A$. 
\end{remark}

In our work here, we will give only application-specific ways of verifying \ref{hyp:restriction-rings} and \ref{hyp:denominator-units}.  By contrast, there is a general criterion for verifying \ref{hyp:restriction-finite} -- finiteness of symplectic cohomology -- coming from work of Ganatra and Pomerleano \cite{GP, Pom}. 

By a {\em  positive normal-crossings compactification}, they mean a smooth projective
compactification whose simple normal-crossings boundary supports
an ample effective divisor.  Given such a compactification, Ganatra and Pomerleano give a multiplicative spectral sequence from (what they call) the logarithmic cohomology of $(\overline X,D)$ to symplectic cohomology of $X$
\cite[Theorem 1.1]{GP}.  One finiteness result they establish is that if the boundary divisor itself is ample and all iterated intersections of boundary strata are connected, then the symplectic cohomology is finitely generated over its degree zero part \cite[Theorem 1.3]{GP}  (which implies finiteness of $\SH^i$ over $\SH^0$). 

We will be interested in the setting where $X$ carries a nowhere-vanishing algebraic volume form, and a positive normal-crossings compactification with
\[
 \operatorname{div}(\eta)=-\sum_i a_iD_i,\qquad a_i\le1.
\]
Pomerleano refers to this as the {\em log nef} setting, and the situation when all $a_i=1$ as {\em log Calabi-Yau}.  He shows 
in the log nef case that $\SH^0(X)$ is a finitely generated $\C$-algebra \cite[Corollary~4.22]{Pom}, and 
that when all $a_i = 1$, in addition
$\SH^*(X)$ is a finite $\SH^0(X)$ module
\cite[Theorem~1.1(a),(b), Lemma~4.24]{Pom}.

\begin{warning}
    In many references, e.g. \cite{GHK}, a variety $X$ is said to be {\em log Calabi-Yau} if, on a smooth
normal-crossings compactification $(\overline X,D)$,
\[
 H^0\bigl(\overline X,\omega_{\overline X}(D)^{\otimes m}\bigr)
       =\C\eta^{\otimes m}\qquad(m\ge1)
\]
for a nowhere-vanishing volume form $\eta$ on $X$
\cite[Definition~1.1]{GHK}.
This permits $\operatorname{div}(\eta)+D$ to be nonzero and effective.
Pomerleano uses {\em log Calabi-Yau} to mean $D$ itself is
anticanonical
\cite[Section~1 and Definition~3.1]{Pom}.
\end{warning}

Under the weaker log nef assumption, we have a corresponding weakening of \cite[Lemma 4.24]{Pom}: 

\begin{lemma} \label{lem:degreewise} 
Let $X$ be smooth affine with a nowhere-vanishing algebraic volume
form $\eta$. Suppose a positive normal-crossings compactification
$(\overline X,D)$ satisfies
\[
 \operatorname{div}(\eta)=-\sum_i a_iD_i,\qquad a_i\le1.
\]
For the affine Weinstein structure and the grading determined by $\eta$, every $\SH^j(X)$ is finite over
$\SH^0(X)$.
\end{lemma}
\begin{proof}
The log cohomology has summands indexed by the intersections of the components of the boundary divisor; writing $\mathring S_I$ for the real normal torus bundle of the open locus of such a stratum, the summands are
\[
 t^{\mathbf1_I}H^*(\mathring S_I;\C)[t_i:i\in I],\qquad
 \deg(\alpha t^{\mathbf v})
 =\deg\alpha+2\sum_i(1-a_i)v_i.
\]
The summand for $I=\varnothing$ is $H^*(X;\C)$.

The degree-zero part consists of permanent cycles
\cite[Theorem~4.21 and Corollary~4.22]{Pom}.
It receives an algebra map
\[
 P=\C[T_i:a_i=1]\longrightarrow H^0_{\log}(\overline X,D),
 \qquad T_i\longmapsto 1_{\mathring S_{\{i\}}}t_i.
\]
This map need not be an isomorphism.
On a summand supported on $I$, multiplication by $T_i$ for $i\in I$
increases $v_i$ by one.
In fixed degree $j$, the degree formula bounds every $v_i$ with
$a_i<1$.
Each $H^*(\mathring S_I;\C)$ is finite-dimensional, and there are
finitely many $I$.
Thus a finite set of classes with $v_i=1$ for $i\in I$ and $a_i=1$
generates $H^j_{\log}(\overline X,D)$ over $P$, hence over
$H^0_{\log}(\overline X,D)$.
Taking $j=0$ shows that the degree-zero ring is finite over $P$,
and hence Noetherian.
Every degree-$j$ term on the limiting page is therefore finite over it.
Lifting generators through the bounded-below weight filtration, as in
\cite[Lemma~4.24]{Pom}, proves that $\SH^j(X)$ is finite over $\SH^0(X)$.
\end{proof}

For us, the necessary volume forms come from the following:

\begin{lemma}\label{lem:log-volume} 
Suppose $X$ is smooth affine and contains a dense algebraic torus
with coordinates $z_1,\ldots,z_d$. If
$\eta=c\,\dd\log z_1\wedge\cdots\wedge\dd\log z_d$, $c\ne0$,
extends to a nowhere-vanishing volume form on $X$, then $(X, \eta)$ is log nef in the sense of Pomerleano. 
\end{lemma}
\begin{proof}
This is \cite[Lemma~1.4(1)]{GHK}, applied to
$(\C^*)^d\subset X$.  Or explicitly:
resolve a projective compactification to obtain a positive
normal-crossings compactification.
At a boundary divisor with uniformizer $u$, write $z_i=u^{m_i}v_i$
with $v_i$ a unit at its generic point. Then
\[
 \dd\log z_i=m_i\frac{\dd u}{u}+\frac{\dd v_i}{v_i}.
\]
In their wedge at most one factor $\dd u/u$ occurs. Thus $\eta$
has pole order at most one along every boundary divisor,
including the exceptional divisors.
\end{proof}

\section{Freezing principal cluster mirror symmetry}

\subsection{Freezing of principal cluster algebras}
Recall our conventions from Section~\ref{sec:intro}. Let \(B\in \Mat_m(\bZ)\) be a skew-symmetrizable matrix, and \(D\) be its right symmetrizer:
\begin{equation}\label{eq:right-symmetrizer}
 D=\diag(d_1,\ldots,d_m),\qquad BD=-DB^{\mathsf T},\qquad
 B^\vee=-B^{\mathsf T}=D^{-1}BD.
\end{equation}
The dual matrix \(B^\vee\) has right symmetrizer $D^{-1}$.

For a seed $\fs$, let $L_\fs$ be its Laurent ring, including the inverted
frozen variables. The upper cluster algebra is
$\U(\widetilde B)=\bigcap_\fs L_\fs$ inside the common rational function
field. If \(\widetilde B\) is full rank, then by \cite[Corollary 1.9]{BFZ}, 
\begin{equation}\label{eq:finite-star}
 \U(\widetilde B)=L_0\cap\bigcap_{i=1}^m L_i,
\end{equation}
where $L_i$ is obtained by one mutation at $i$. 
We write 
\[
 \mathcal{U}(\widetilde B)=\Spec\U(\widetilde B),\qquad
 \mathcal{U}_{\pr}(B)=\Spec\U\binom B{I_m}.
\]

If $J$ is a set of mutable indices, write $\fr_J\widetilde B$ for the
seed obtained by freezing them. Its matrix is obtained by deleting
the columns indexed by $J$, but not the corresponding rows.

\begin{lemma}\label{lem:upper-freezing}
 If $\widetilde B$ has full rank and $w=\prod_{j\in J}x_j$, then
\begin{equation}\label{eq:upper-freezing}
 \U(\widetilde B)[w^{-1}]=\U(\fr_J\widetilde B).
\end{equation}
\end{lemma}
\begin{proof}
By~\eqref{eq:finite-star} we have
$
 \U(\widetilde B)[w^{-1}]
 =L_0[w^{-1}]\cap\bigcap_{i=1}^m L_i[w^{-1}]$,
since localization commutes with finite intersections.
The variables $x_j$, $j\in J$, are already invertible in $L_0$
and in $L_i$ for $i\notin J$. Thus
$L_0[w^{-1}]=L_0$ and  $L_i[w^{-1}]=L_i\quad(i\notin J)$.
For $j\in J$, the exchange relation $x_jx_j'=P_j$ gives
$L_j[w^{-1}]=L_j[x_j^{-1}]=L_0[P_j^{-1}]\supseteq L_0.$
 Consequently,
$
 \U(\widetilde B)[w^{-1}]
 =L_0\cap\bigcap_{i\notin J}L_i$.
Observe that \(\fr_J \widetilde B\) is still of full rank, so
\eqref{eq:finite-star} still holds after freezing. This proves the claim. 
\end{proof}

For \(J \subseteq \{1,\ldots,m\}\) we write \(x_J = \prod_{j\in J} x_j\), \(\sU(\widetilde B)_J = \sU(\widetilde B)_{x_J}\), and \( \cU (\widetilde B)_J \coloneqq \Spec \sU(\widetilde B)_J.\)
Recall the following lemma.
\begin{lemma} \cite[Lemma 3.2]{HS} \ \label{lem:coefficients}
Let $\widetilde B_0=(B;C)$ have full rank. Let $P$ be an integral
$f\times m$ matrix and $Q$ an integral $f\times f$ matrix with
$\det Q\ne0$. Then
\begin{equation}\label{eq:coefficient-change}
 \U(B;C)\otimes_{\C[y^{\pm1}]}\C[v^{\pm1}]
       \cong\U(B;PB+QC),
\end{equation}
where $y^a\mapsto v^{Qa}$ and $x_i\mapsto u_i v^{P_i}$.
\end{lemma}

Consider the cluster Weil--Petersson form \cite{GSV05, FG09},
given by
\begin{equation}\label{eq:principal-omega}
 \Omega_{B,D}=\frac12\sum_{i,j}B_{ij}d_j\,
      \dd\log x_i\wedge\dd\log x_j
       +\sum_i d_i\,\dd\log y_i\wedge\dd\log x_i.
\end{equation}
We abbreviate it to $\Omega_B$ when the symmetrizer is fixed.

\begin{lemma} \label{lem: form}
If $\cU_{\pr}(B)$ is smooth and of finite type, then $\Omega_B$
extends uniquely to a regular algebraic symplectic form on
$\cU_{\pr}(B)$.
\end{lemma}
\begin{proof}
Put $X=\cU_{\pr}(B)$ and let
$X^\circ=\bigcup_{\fs}T_{\fs}$ be the union of its seed tori.
The Laurent phenomenon identifies each $T_{\fs}=\Spec L_{\fs}$
with the principal open subset of $X$ where its mutable variables
are invertible. By the definition of the upper cluster algebra,
\[
 \Gamma(X^\circ,\cO_{X^\circ})
 =\bigcap_{\fs}L_{\fs}
 =\Gamma(X,\cO_X).
\]
Then $\operatorname{codim}_X(X\setminus X^\circ)\ge2$
by \cite[Lemma~4.4.1]{Matherne-Muller}.

The form~\eqref{eq:principal-omega} is preserved by mutation
\cite[Corollary~6.5]{FG09}.
Thus $\Omega_B$ is regular, closed, and nondegenerate on $X^\circ$.
Since $X$ is smooth, it extends uniquely to a regular algebraic
symplectic form on $X$.
\end{proof}

In the setting of principal coefficients, we use the following completion:
\begin{equation}\label{eq:completion}
 E_B=\begin{pmatrix}B&-I\\ I&0\end{pmatrix},\qquad
 \widetilde D=\diag(D,D),\qquad
 K_B=E_B\widetilde D=\begin{pmatrix}BD&-D\\D&0\end{pmatrix}.
\end{equation}
The matrix $E_B$ is integral unimodular and $K_B$ is skew-symmetric.

\begin{lemma} \label{lem: commute mutation and principal}
    The above change of frozen variables gives an isomorphism between mutation of the upper principal cluster algebra of $B$ and the upper principal cluster algebra of a mutation of $B$, respecting the cluster Weil--Petersson form and the positive locus.
\end{lemma}
\begin{proof}
    Mutating \eqref{eq:completion} gives
\[
 E'=\begin{pmatrix}
 B'&-D(C')^{\mathsf T}D^{-1}\\ C'&0
 \end{pmatrix}.
\]
The frozen block remains zero because the columns of $C'$ are
sign-coherent \cite[Theorem~2.4]{CL}.
Mutation preserves the determinant by \cite[Lemma~4.11]{HS}, so
$(\det C')^2=\det E'=1$.
Lemma~\ref{lem:coefficients} therefore gives the isomorphism
\[
 (x';y)\longmapsto(x';Y),\qquad
 Y_i=\prod_a y_a^{C'_{ai}}.
\]
Substituting $\dd\log Y=(C')^{\mathsf T}\dd\log y$ identifies
the logarithmic forms. The monomial map preserves the positive locus.
\end{proof}

\begin{lemma}\label{lem:principal-product}
Let $I\sqcup J=\{1,\ldots,m\}$.
There is an isomorphism
\begin{equation}\label{eq:principal-product}
 F:\mathcal{U}_{\pr}(B)_{J}
 \xrightarrow{\sim}\mathcal{U}_{\pr}(B_{I,I})\times(\C^*)^{2|J|}.
\end{equation}
On the initial seed torus, $F$ is given by
\[
 (x;y)\longmapsto
 \bigl((x_i;Y_i)_{i\in I},(x_j;s_j)_{j\in J}\bigr),
\]
where
\begin{equation}\label{eq:principal-freezing-coordinates}
 Y_i=y_i\prod_{j\in J}x_j^{B_{ji}},\qquad
 s_j=y_j\prod_{\ell\in J}x_\ell^{-H_{j\ell}},\qquad
 H_{ab}=\begin{cases}-B_{ba},&a<b,\\0,&a\ge b.\end{cases}
\end{equation}
Here $a,b\in J$, with the order inherited from $\{1,\ldots,m\}$.
As rational differential forms,
\begin{equation}\label{eq:product-omega}
 \Omega_{B,D}
 =F^*\left(\Omega_{B_{I,I},D_I}
       +\sum_{j\in J}d_j\,\dd\log s_j\wedge\dd\log x_j\right).
\end{equation}
\end{lemma}
\begin{proof}
The inverse coordinate substitution is
\[
 y_i=Y_i\prod_{j\in J}x_j^{-B_{ji}}\quad(i\in I),\qquad
 y_j=s_j\prod_{\ell\in J}x_\ell^{H_{j\ell}}\quad(j\in J),
\]
with all $x_k$ unchanged. After freezing $J$, order the old frozen variables
as $(y_i)_{i\in I},(x_j)_{j\in J},(y_j)_{j\in J}$,
and the new ones as $(Y_i)_{i\in I},(x_j)_{j\in J},(s_j)_{j\in J}$.
In the notation of Lemma~\ref{lem:coefficients}, this substitution has $P=0$ and
\[
 C=\begin{pmatrix}I_{|I|}\\B_{J,I}\\0\end{pmatrix},\qquad
 Q=\begin{pmatrix}
 I_{|I|}&0&0\\-B_{J,I}&I_{|J|}&H^{\mathsf T}\\0&0&I_{|J|}
 \end{pmatrix},\qquad
 QC=\begin{pmatrix}I_{|I|}\\0\\0\end{pmatrix}.
\]
The matrices $\binom{B_{I,I}}{C}$ and $\binom{B_{I,I}}{QC}$ have full rank, and $\det Q=1$.
Lemmas~\ref{lem:upper-freezing} and~\ref{lem:coefficients} therefore give
\[
 \U\binom{B}{I_m}[x_J^{-1}]
 \xleftarrow[\ F^*\ ]{\sim}\U\binom{B_{I,I}}{I_{|I|}}
       [x_j^{\pm1},s_j^{\pm1}\mid j\in J].
\]
The Laurent variables on the right have zero exchange rows;
\eqref{eq:finite-star} gives the displayed Laurent extension.
Taking spectra proves~\eqref{eq:principal-product}.

For the differential forms, $B_{ab}d_b=-d_aB_{ba}$ gives
\begin{equation}\label{eq:weighted-shear}
 D_JH-H^{\mathsf T}D_J=B_{J,J}D_J.
\end{equation}
On the initial seed torus, \eqref{eq:principal-freezing-coordinates} and \eqref{eq:weighted-shear} give
\[
 \begin{aligned}
 \sum_{i\in I}d_i\,\dd\log Y_i\wedge\dd\log x_i
 &=\sum_{i\in I}d_i\,\dd\log y_i\wedge\dd\log x_i +\sum_{i\in I,\,j\in J}B_{ij}d_j\, \dd\log x_i\wedge\dd\log x_j,\\
 \sum_{j\in J}d_j\,\dd\log s_j\wedge\dd\log x_j
 &=\sum_{j\in J}d_j\,\dd\log y_j\wedge\dd\log x_j + \frac12\sum_{j,\ell\in J}B_{j\ell}d_\ell\, \dd\log x_j\wedge\dd\log x_\ell.
 \end{aligned}
\]
Summing these two identities and adding
$\frac12\sum_{i,k\in I}B_{ik}d_k\,
\dd\log x_i\wedge\dd\log x_k$ to both sides gives~\eqref{eq:product-omega}.
\end{proof}

\subsection{\(\cX\) versus \(\cA\) coordinates}\label{sec:mirror-coordinates}

A rank-one local system on the unitary seed torus $L_\fs$ is a point of
$\operatorname{Hom}(H_1(L_\fs;\Z),\C^*)
 =\Spec\C[H_1(L_\fs;\Z)].$
Its monodromies around the seed coordinate circles are the mirror \(\cX\)
coordinates $\chi_{\fs,i}$ in Definition~\ref{def: cms}.
For principal coefficients the mirror is also
$\cU_{\pr}(-B^{\mathsf T})$, with cluster coordinates $(a_i,b_i)$.
We use these $\cA$ coordinates for the freezing product of
Lemma~\ref{lem:principal-product}. 
Here we compute the integral change of circle basis that relates
$(\chi_{\fs,i})$ to $(a_i,b_i)$. 

Write $z_i=e^{u_i+\mathrm i\theta_i}$ on a seed torus, for $1\le i\le m+f$.
The unitary seed torus is parameterized by
\[
 (S^1)^{m+f}\longrightarrow L_{\fs},\qquad
 \theta\longmapsto(e^{\mathrm i\theta_1},\ldots,e^{\mathrm i\theta_{m+f}}).
\]
These parameterizations are chosen separately in each seed.
On the mirror quotient chart, order the coordinates as $(X,Z,U)$.
Its weight matrix and invariant characters are
\begin{equation}
\label{eq:change_A_X}
 W=[I_m\mid B^{\mathsf T}\mid C^{\mathsf T}],\qquad
 \chi_i=Z_i\prod_jX_j^{-B_{ij}},\qquad
 \chi_{m+a}=U_a\prod_jX_j^{-C_{aj}}.
\end{equation}
Here $1\leq i\leq m$ and $1\leq a\leq f$.
Consider the identification 
\[
 \Z^{m+f}\xrightarrow{\sim}\ker W,\qquad
 \binom vw\longmapsto
 \begin{pmatrix}-B^{\mathsf T}v-C^{\mathsf T}w\\v\\w\end{pmatrix}.
\]
Note that \(\ker W\) is the character lattice of the seed torus. Hence the monodromy of the $i$th coordinate loop corresponds to $\chi_i$. These \(\chi_i\) form the cluster \(\cX\)-coordinates~\cite{FG09}, and the mutation is given by
\[
 \chi'_k=\chi_k^{-1},\qquad
 \chi'_i=\chi_i\chi_k^{[-\widetilde B_{ik}]_+}
                  (1+\chi_k)^{\widetilde B_{ik}}\quad(i\ne k).
\]
Under~\eqref{eq:change_A_X}, this is compatible with the mutation of cluster \(\cA\) variables~\cite[Proposition 2.2]{FG09}.

For principal coefficients $C=I_m$, the frozen variables $U_i$ have
the standard basis weights. Setting $U_i=1$ in the quotient gives
\[
 \cX\left(\binom{B}{I_m}^{\vee}\right)
 \cong\cU_{\pr}(-B^{\mathsf T}),\qquad
 a_i=X_i/U_i,\qquad b_i=Z_i\prod_jU_j^{-B_{ij}},
\]
where $a_i$ and \(b_i\) are the mutable and frozen cluster \(\cA\) variables.
On the initial chart the two coordinate systems satisfy
\begin{equation}\label{eq:principal-XA}
 \chi_i=b_i\prod_j a_j^{-B_{ij}},\qquad
 \chi_{m+i}=a_i^{-1}.
\end{equation}
In particular this is an isomorphism of tori, with inverse
$a_i=\chi_{m+i}^{-1}$ and
$b_i=\chi_i\prod_j\chi_{m+j}^{-B_{ij}}$.

Use the completion $E_B$ as in \eqref{eq:completion}. On the initial seed torus, put
\[
 \lambda_{\mathrm{tor}}=u^{\mathsf T}K_B\,\dd\theta,
 \qquad \dd\lambda_{\mathrm{tor}}=\Im\Omega_B.
\]
Since $K_B=E_B\widetilde D=-\widetilde D E_B^{\mathsf T}$, the coordinates
\begin{equation}\label{eq:marking}
 \varphi=-E_B^{\mathsf T}\theta,\qquad p=\widetilde D u
\end{equation}
give $\lambda_{\mathrm{tor}}=p^{\mathsf T}\dd\varphi$.
The positive locus is the cotangent fiber $\varphi=0$.

Denote the base isomorphism by
\begin{align*}
 g:(S^1)^{2m}_{\varphi} &\longrightarrow(S^1)^{2m}_{\theta}\\
 \varphi & \longmapsto -(E_B^\T)^{-1} \varphi
\end{align*}
Its cotangent lift sends $(\varphi,p)$ to
$(-(E_B^{\mathsf T})^{-1}\varphi,-E_Bp)$.
For a local system with $\theta$-circle monodromies $\chi_j$, the
monodromies around the $\varphi_{x_i}$ and $\varphi_{y_i}$ circles are
\[
 \chi_{m+i}^{-1}=a_i,\qquad
 \chi_i\prod_j\chi_{m+j}^{-B_{ij}}=b_i,
\]
which is compatible with~\eqref{eq:principal-XA}.

After any mutation sequence, use the full mutated square matrix $E_{\fs}$.
Then
\[
 \varphi_{\fs}=-E_{\fs}^{\mathsf T}\theta_{\fs},\qquad
 \chi_{\fs,i}=\prod_j a_{\fs,j}^{-(E_{\fs})_{ij}}.
\]
Here $a_{\fs,1},\ldots,a_{\fs,2m}$ denote all the dual $\cA$ coordinates,
including the frozen ones.
Write $g_{\fs}$ for the base map
$\theta_{\fs}=-E_{\fs}^{-\mathsf T}\varphi_{\fs}$.

We now translate Definition~\ref{def: cms}.
Let $T_{\fs}^{\cA}$ be the dual seed torus with coordinates $a_{\fs,j}$,
and let $F_0=T^*_0(S^1)^{2m}_{\varphi_{\fs}}$ be the cotangent fiber over zero.
For $1\le j\le2m$, let $\gamma_j$ be the based loop along which
$\varphi_{\fs,j}$ increases from $0$ to $2\pi$, with all other angles zero.
Denote the equivalence~\cite{Ab} by
\[
 \beta_{\fs}^{\cA}:\Fuk(T^*(S^1)^{2m}_{\varphi_{\fs}})
 \xrightarrow{\sim}\Perf(T_{\fs}^{\cA}),\qquad
 K\longmapsto\RHom(F_0,K),
\]
where we identify
\[
 \operatorname{REnd}(F_0)\simeq\C[\pi_1((S^1)^{2m}_{\varphi_{\fs}},0)]
 \xrightarrow{\sim}\C[a_{\fs,1}^{\pm1},\ldots,a_{\fs,2m}^{\pm1}],
 \qquad [\gamma_j]\longmapsto a_{\fs,j}.
\]
Thus $\beta_{\fs}^{\cA}(F_0)=\OO_{T_{\fs}^{\cA}}$, and the mirror coordinate
$a_{\fs,j}$ records monodromy around the positively oriented
$\varphi_{\fs,j}$ circle. The monomial isomorphism is
\[
 h_{\fs}:T_{\fs}^{\cA}\xrightarrow{\sim}\Spec A_{\fs},\qquad
 h_{\fs}^*t_i=\prod_j a_{\fs,j}^{-(E_{\fs})_{ij}}.
\]
If $(g_{\fs})_*$ denotes the equivalence induced by the cotangent lift,
the group-algebra descriptions give
\[
 \beta_{\fs}^{\cA}(g_{\fs})_*^{-1}
 \simeq h_{\fs}^*\chi_{\fs}.
\]
Put $\mathrm{Vit}_{\fs}^{\cA}=(g_{\fs})_*^{-1}\mathrm{Vit}_{\fs}$
and $j_{\fs}^{\cA}=j_{\fs}h_{\fs}$, using the global principal
mirror identification above. The required comparison becomes
\[
 \chi_{\fs}\mathrm{Vit}_{\fs}\Phi\simeq j_{\fs}^*
 \quad\Longleftrightarrow\quad
 \beta_{\fs}^{\cA}\mathrm{Vit}_{\fs}^{\cA}\Phi\simeq(j_{\fs}^{\cA})^*.
\]
The specified isomorphism $\Phi(\OO)\simeq L_+$ is unchanged.
The base maps fix zero, so their cotangent lifts carry the fiber over zero
to the fiber over zero. Use the fiber gradings and Lie group spin structures.
The same translation transports an equivariant structure and its coherences.

\begin{example}
For $\widetilde B=\binom01$, write
\[
 \Nn=\{ab=1+z\}\subset\C^2\times\C^*,\qquad
 \Nn^\vee=\{a^\vee b^\vee=1+z^\vee\}\subset\C^2\times\C^*.
\]
The unitary seed tori are $L_a=\{|a|=|z|=1\}$ and
$L_b=\{|b|=|z|=1\}$ in $\Nn$.
The parameterizations in \eqref{eq:marking} are
\[
 \varphi_a=(-\arg z,\arg a),\qquad
 \varphi_b=(\arg z,-\arg b).
\]
Monodromies around these circles give the $\cA$ coordinates
$(a^\vee,z^\vee)$ and $(b^\vee,z^\vee)$ used in
\cite[Section~5.3]{HS}.

Let $(\chi_1,\chi_2)$ be the monodromies around $(\arg a,\arg z)$,
and $(\chi'_1,\chi'_2)$ those around $(\arg b,\arg z)$.
In coordinates the chart maps are given by
\[
 \begin{aligned}
 \operatorname{Hom}(H_1(L_a;\Z),\C^*)&\xrightarrow{\sim}D(a^\vee)\subset\Nn^\vee,\\
 (\chi_1,\chi_2)&\longmapsto
   \bigl(\chi_2^{-1},(1+\chi_1)\chi_2,\chi_1\bigr),\\[4pt]
 \operatorname{Hom}(H_1(L_b;\Z),\C^*)&\xrightarrow{\sim}D(b^\vee)\subset\Nn^\vee,\\
 (\chi'_1,\chi'_2)&\longmapsto
   \bigl((1+(\chi'_1)^{-1})/\chi'_2,\chi'_2,(\chi'_1)^{-1}\bigr).
 \end{aligned}
\]
On the overlap, the $\cX$ mutation
\[
 \chi'_1=\chi_1^{-1},\qquad \chi'_2=\chi_2(1+\chi_1)
\]
becomes the exchange relation $b^\vee=(1+z^\vee)/a^\vee$.

With $j_a:D(a^\vee)\hookrightarrow\Nn^\vee$ and
$j_b:D(b^\vee)\hookrightarrow\Nn^\vee$, cluster mirror symmetry in
$\cA$ coordinates takes the form
\[
 \Phi_{\Nn}:\Perf(\Nn^\vee)\xrightarrow{\sim}\Fuk(\Nn),\qquad
 \Phi_{\Nn}(\OO)\simeq L_+,\qquad
 \beta_a^{\cA}\mathrm{Vit}_a^{\cA}\Phi_{\Nn}\simeq j_a^*,\quad
 \beta_b^{\cA}\mathrm{Vit}_b^{\cA}\Phi_{\Nn}\simeq j_b^*.
\]

The equivalence is defined on all of $\Nn^\vee$, including $(0,0,-1)$,
which neither chart contains.
An element $\kappa\in\C^*$ acts on the mirror by
\[
(a^\vee,b^\vee,z^\vee)\longmapsto
(\kappa^{-1}a^\vee,\kappa b^\vee,z^\vee), \qquad \kappa \in \bC^*. 
\]
\end{example}

We now compute the monodromies under the product isomorphism $F$ of
Lemma~\ref{lem:principal-product}, using the $\cA$ coordinates above.

\begin{lemma}\label{lem:dual-product}
Let $F$ be the isomorphism of Lemma~\ref{lem:principal-product}.
Use~\eqref{eq:marking} to identify its initial seed tori with cotangent bundles,
applying the same rule to each factor of the product.
Write $\varphi$ and $\varphi_{\mathrm{prod}}$ for the original and product base angles.
Then $F^{-1}$ covers an isomorphism
\[
 g:(S^1)^{2m}_{\varphi_{\mathrm{prod}}}
       \longrightarrow (S^1)^{2m}_{\varphi}.
\]
All coordinate circles are based at zero and positively oriented.
Let $x_k^\vee,y_k^\vee$ be the monodromies of a rank-one complex local system
$\mathcal L$ around the original coordinate circles
$\varphi_{x_k},\varphi_{y_k}$.
The monodromies of $g^*\mathcal L$ around the product coordinate circles
$\varphi_{\mathrm{prod},x_k}$, $\varphi_{\mathrm{prod},Y_i}$,
and $\varphi_{\mathrm{prod},s_j}$ are, respectively,
$x_k^\vee$, $Y_i^\vee$, and $s_j^\vee$, where
\begin{equation}\label{eq:dual-product-coordinates}
 \begin{aligned}
 Y_i^\vee&=y_i^\vee\prod_{j\in J}(x_j^\vee)^{-B_{ij}}
                   &&(i\in I),\\
 s_j^\vee&=y_j^\vee\prod_{\ell\in J}(x_\ell^\vee)^{-H^\vee_{j\ell}}
                   &&(j\in J),\qquad H^\vee=D_JHD_J^{-1}.
 \end{aligned}
\end{equation}
These are the coordinate formulas of Lemma~\ref{lem:principal-product}
for $B^\vee=-B^{\mathsf T}$.
\end{lemma}
\begin{proof}
Order the original coordinates as $(x_1,\ldots,x_m;y_1,\ldots,y_m)$.
In the product coordinates, keep $x_k$ in position $k$ and replace $y_k$
by $Y_k$ for $k\in I$ and by $s_k$ for $k\in J$.
Writing $z=e^{u+\mathrm i\theta}$ in each coordinate, the formulas
of Lemma~\ref{lem:principal-product} give
\[
 u_{\mathrm{prod}}=Tu,\qquad\theta_{\mathrm{prod}}=T\theta,\qquad
 T=\begin{pmatrix}I_m&0\\ A&I_m\end{pmatrix},\qquad
 A_{ab}=\begin{cases}
 B_{ba},&a\in I,\ b\in J,\\
 -H_{ab},&a,b\in J,\\
 0,&\text{otherwise}.
 \end{cases}
\]
The product has completion
\[
 E_{\mathrm{prod}}=\begin{pmatrix}B_{\mathrm{prod}}&-I_m\\I_m&0\end{pmatrix},
 \qquad (B_{\mathrm{prod}})_{ab}
 =\begin{cases}B_{ab},&a,b\in I,\\0,&\text{otherwise}.
 \end{cases}
\]
Thus the cotangent coordinates on the two sides are
\[
 \begin{aligned}
 \varphi&=-E_B^{\mathsf T}\theta,&p&=\widetilde D u,\\
 \varphi_{\mathrm{prod}}&=-E_{\mathrm{prod}}^{\mathsf T}\theta_{\mathrm{prod}},
 &p_{\mathrm{prod}}&=\widetilde D u_{\mathrm{prod}}.
 \end{aligned}
\]
Equation~\eqref{eq:product-omega} gives
$E_B\widetilde D=T^{\mathsf T}E_{\mathrm{prod}}\widetilde D T$.
Since $E_B\widetilde D$ and $E_{\mathrm{prod}}\widetilde D$ are skew-symmetric,
\[
 E_B^{\mathsf T}
 =\widetilde D^{-1}T^{\mathsf T}\widetilde D
       E_{\mathrm{prod}}^{\mathsf T}T.
\]
Substituting $\theta_{\mathrm{prod}}=T\theta$ therefore gives
\begin{equation}\label{eq:angle-map}
 \varphi=M\varphi_{\mathrm{prod}},\qquad
 M=\widetilde D^{-1}T^{\mathsf T}\widetilde D
   =\begin{pmatrix}I_m&D^{-1}A^{\mathsf T}D\\0&I_m\end{pmatrix}.
\end{equation}
Also $p_{\mathrm{prod}}=M^{\mathsf T}p$.
Thus, in cotangent coordinates,
\[
 F^{-1}(\varphi_{\mathrm{prod}},p_{\mathrm{prod}})
   =(M\varphi_{\mathrm{prod}},M^{-\mathsf T}p_{\mathrm{prod}}).
\]

Let $e_{x_k},e_{y_k}$ be the original coordinate circle classes in $\pi_1$,
and denote the product classes by
$e^{\mathrm{prod}}_{x_k},e^{\mathrm{prod}}_{Y_i},e^{\mathrm{prod}}_{s_j}$.
Using $d_iB_{ji}/d_j=-B_{ij}$ in~\eqref{eq:angle-map}, we obtain
\[
 \begin{aligned}
 g_*(e^{\mathrm{prod}}_{x_k})&=e_{x_k}&&(1\le k\le m),\\
 g_*(e^{\mathrm{prod}}_{Y_i})&=e_{y_i}-\sum_{j\in J}B_{ij}e_{x_j}
                                      &&(i\in I),\\
 g_*(e^{\mathrm{prod}}_{s_j})&=e_{y_j}-\sum_{\ell\in J}H^\vee_{j\ell}e_{x_\ell}
                                      &&(j\in J).
 \end{aligned}
\]
For $a<b$, $H^\vee_{ab}=-(d_a/d_b)B_{ba}=B_{ab}$;
all other entries of $H^\vee$ vanish.
Thus $M$ and $M^{-1}$ are integral, so $g$ is an isomorphism of tori.
The rank-one local system $\mathcal L$ is on
$(S^1)^{2m}_{\varphi}$. For a based loop class $w$ in this torus,
write $\operatorname{mon}_{\mathcal L}(w)\in\C^*$ for the scalar
by which parallel transport around $w$ acts on the fiber
$\mathcal L_0$.
For $v\in\pi_1((S^1)^{2m}_{\varphi_{\mathrm{prod}}},0)$,
we have $g_*v=Mv$, so
\[
 \operatorname{mon}_{g^*\mathcal L}(v)
 =\operatorname{mon}_{\mathcal L}(Mv)
 =\prod_{k=1}^m(x_k^\vee)^{(Mv)_{x_k}}
                     (y_k^\vee)^{(Mv)_{y_k}}.
\]
Applying this formula to the three types of coordinate circle classes gives
\eqref{eq:dual-product-coordinates} and the unchanged monodromies $x_k^\vee$.
Finally, $B^\vee_{ji}=-B_{ij}$ and
$H^\vee_{ab}=-B^\vee_{ba}$ for $a<b$, as required by
Lemma~\ref{lem:principal-product} for $B^\vee$.
\end{proof}

\subsection{Liouville geometry}

Here we construct cluster-admissible Liouville forms on principal cluster varieties, and study their properties.  A sufficient condition in the following is that $B$ has a reddening sequence, since this gives a large cluster complex
per \cite[Lemma~8.29]{GHKK}.

\begin{lemma}\label{lem:quadratic}
Suppose $B$ has a large cluster complex in the sense of \cite{GHKK} and $\cU_{\pr}(B)$ is smooth affine
of finite type, with regular nondegenerate form $\Omega_B$.
Then there is a closed embedding
$\iota:\cU_{\pr}(B)\hookrightarrow\C^N$, defined over $\R$, such that
$
 \{\iota_\mu,\iota_\nu\}_{\Omega_B}\in
 \operatorname{span}_{\C} \langle 1,\iota_\rho,\iota_\rho\iota_\sigma\rangle.
$
Here $\{-,-\}_{\Omega_B}$ denotes the Poisson bracket induced by
$\Omega_B$. The set on the right consists of the constant function
$1$, the coordinate functions $\iota_\rho$, and their products
$\iota_\rho\iota_\sigma$, for $1\le\rho,\sigma\le N$.
For any algebraic torus acting diagonally on the seed and preserving
the exchange relations, the coordinates $\iota_\mu$ may be chosen
to be weight vectors.
\end{lemma}
\begin{proof}
 The  construction in
\cite[proof of Lemma~B.1]{HS} gives an exhaustive filtration
$F_\bullet$ of $\U_{\pr}(B)$ by finite-dimensional real,
torus-invariant subspaces, with finitely generated Rees algebra.
We use this construction before coefficient specialization;
the assumption $\A=\U$ in that lemma is needed for the
specialization, not for this principal-coefficient filtration.

The filtration is defined by monomial valuations on seed charts.
The bracket is log-canonical on each chart, so
$\{F_d,F_e\}_{\Omega_B}\subset F_{d+e}$ by
\cite[proof of Lemma~B.2]{HS}.
For sufficiently divisible $d$, the Veronese Rees algebra is
generated in degree one, and $F_{nd}=(F_d)^n$.
A real weight basis $1,\iota_1,\ldots,\iota_N$ of $F_d$ therefore
generates $\U_{\pr}(B)$ and satisfies
\[
 \{\iota_\mu,\iota_\nu\}_{\Omega_B}\in F_{2d}
 =(F_d)^2
 =\operatorname{span}_{\C}\{1,\iota_\rho,\iota_\rho\iota_\sigma\}.
 \qedhere
\]
\end{proof}

\begin{lemma}\label{lem:principal-geometry}
Assume that $X=\cU_{\pr}(B)$ is smooth and of finite type. Let  $\Omega_B$ be its regular nondegenerate form. Suppose $X$ has a real closed embedding
with quadratic coordinate Poisson brackets, equivariant for the action
\begin{equation}\label{eq:principal-compact-action}
 \T=(S^1)^m,\qquad
 e^{i\gamma}\cdot(x_i,y_i)
 =\left(e^{-i\gamma_i}x_i,
        e^{i\sum_a B_{ai}\gamma_a}y_i\right).
\end{equation}
Suppose $L_+$ is proper and there is a one-form $\alpha_0$ with
$\dd\alpha_0=\Im\Omega_B$. Assume that for each seed torus $T$
there is a smooth function $f_T:T\to\R$ such that
\[
 \alpha_0|_T=\lambda_{\mathrm{tor}}+\dd f_T.
\]
Then $X$ has a finite-type Morse--Bott Weinstein Liouville form
$\lambda$ with the following properties:
\begin{enumerate}
\item $\dd\lambda=\Im\Omega_B$ and, for each seed torus $T$,
      there is a smooth function $g_T:T\to\R$ such that
      $\lambda|_T=\lambda_{\mathrm{tor}}+\dd g_T$;
\item $\lambda|_{L_+}=0$;
\item $\lambda$ and its Lyapunov function are $\T$-invariant, and
      $Z_\lambda\log|y_a|=\log|y_a|$;
\item the critical components are isotropic $\T$-orbits.
\end{enumerate}
There is a graded symplectomorphism from $(X,\dd\lambda)$ to an
affine K\"ahler form on $X$. The grading on a seed torus is the
cotangent grading of \eqref{eq:marking}, normalized at $L_+$.
These conclusions also hold after a product with cotangent tori.
\end{lemma}
\begin{proof}
We average $\alpha_0$, construct an invariant affine Liouville form,
and compare the two forms by \cite{HSinv}.
We then compute the frozen moment maps and compare the gradings.

Let $\tau$ be complex conjugation. The algebraic torus containing
$\T$ acts diagonally on every seed. Indeed, the ratio of the two
exchange monomials at $i$ is
\[
 \widehat y_i=y_i\prod_a x_a^{B_{ai}},\qquad
 \widehat y_i(e^{i\gamma}\cdot(x,y))=\widehat y_i(x,y).
\]
Thus the exchange monomials have the same weight, and mutation
produces another weight vector. The action extends to the upper
cluster algebra and preserves $\Omega_B$.
It need not fix the coefficient coordinates $y_i$.

With normalized Haar measure on $\T$, put
\[
 \alpha=\int_{\T} \kappa^*\frac{\alpha_0-\tau^*\alpha_0}{2}\,\dd\kappa.
\]
For $\kappa\in\T$, we have
\[
 \dd\alpha=\Im\Omega_B,\qquad
 \kappa^*\alpha=\alpha,\qquad \tau^*\alpha=-\alpha.
\]
On every seed torus, $\kappa^*\lambda_{\mathrm{tor}}=
\lambda_{\mathrm{tor}}$ and
$\tau^*\lambda_{\mathrm{tor}}=-\lambda_{\mathrm{tor}}$.
Averaging therefore preserves exactness of
$\alpha-\lambda_{\mathrm{tor}}$ there. Equivalently,
\[
 \int_\gamma(\alpha-\lambda_{\mathrm{tor}})=0
 \quad\text{for every closed loop }\gamma\subset T.
\]
Since $\tau$ fixes the real locus, $\alpha|_{X(\R)}=0$.

Choose the assumed real weight-vector embedding and a
$\T$-invariant Hermitian norm. The function
$R_0=1+\frac14|\iota|^2$ is proper and strictly plurisubharmonic.
Its critical values form a semialgebraic set of measure zero,
hence a finite set. Properness implies that $\Crit R_0$ is compact.
By \cite[Density Lemma~4.8 and the relative remark after
Lemma~4.9]{Wasserman}, there is a
$\T\rtimes\langle\tau\rangle$-invariant $C^2$-small perturbation $R$
whose critical components are nondegenerate orbits.
Choose $R-R_0$ supported in a compact invariant neighborhood of
$\Crit R_0$. Then $R$ remains strictly plurisubharmonic, and each
connected critical component is a $\T$-orbit.

Put
\[
 \beta=-\dd R\circ I,\qquad
 g=\dd\beta(\,\cdot\,,I\,\cdot\,).
\]
The Liouville field is $Z_\beta=\nabla_gR$.
Outside the perturbation support,
$Z_\beta R=|\dd R|_g^2\le R-1$ by
\cite[Lemma~1.3]{HSinv}. On the remaining compact set,
$|Z_\beta R|$ is bounded.
An integral curve therefore stays in a compact sublevel of $R$
on every bounded time interval, so $Z_\beta$ is complete.
Above the largest critical value, $Z_\beta R>0$.
Thus $\beta$ is a finite-type Morse--Bott Weinstein form.
At a critical point, $\beta=0$; invariance gives
$\dd\beta(\xi_a,\xi_b)=0$ for the orbit tangent fields.
The critical orbits are therefore isotropic.
Strict plurisubharmonicity also gives
$\operatorname{ind}(C)+\dim C\le\dim_{\C}X$ for each critical component $C$.

The quadratic Poisson brackets give
$
 \| (\Im\Omega_B)^{-1}\|_g\le C R$ and $\Im\Omega_B(v,Iv)=0.$
The compact perturbation of $R_0$ changes only the constant in
the first bound. Apply \cite[Theorem~1.4 and proof of
Corollary~1.6]{HSinv} to obtain
$\Phi^*\beta=\alpha+\dd f$,
where $\Phi$ is isotopic to the identity.

The comparison can be made $\T$-equivariant and real.
The four paths in \cite[Section~4]{HSinv} are
\[
 \begin{aligned}
 \mu^{(1)}_t&=\alpha+t k(R)\beta, &
 \mu^{(2)}_t&=\alpha+((1-t)k(R)+tA(R))\beta,\\
 \mu^{(3)}_t&=(1-t)\alpha+A(R)\beta, &
 \mu^{(4)}_t&=((1-t)A(R)+t)\beta,
 \end{aligned}
 \qquad 0\le t\le1.
\]
Here $k,A$ are the functions constructed in that proof.
Each $\mu_t$ is $\T$-invariant and $\tau$-anti-invariant.
Uniqueness in
$\iota_{V_t}\dd\mu_t=-\dot\mu_t$
implies $\kappa_*V_t=V_t$ and $\tau_*V_t=V_t$.
The flows supplied by \cite{HSinv} therefore commute with $\T$ and $\tau$.
Their correction functions
$\int_0^t\phi_s^*(\mu_s(V_s))\,\dd s$ are $\tau$-odd.
Hence $f|_{L_+}=0$.
Properness makes $L_+$ a closed component of $X(\R)$,
so the isotopy preserves it.
Set $\lambda=\Phi^*\beta$, with Lyapunov function $R\circ\Phi$.
This proves (i), (ii), (iv), and the invariance in (iii).

Let $\xi_a$ generate the $a$th circle in
\eqref{eq:principal-compact-action}. Direct contraction, or
\cite[equation~(46)]{HS}, gives
\[
 \iota_{\xi_a}\Im\Omega_B=d_a\,\dd\log|y_a|.
\]
Cartan's formula implies
$\dd(\lambda(\xi_a)+d_a\log|y_a|)=0$.
For a $2\pi$-periodic orbit $\gamma_a$ in a seed torus,
\[
 0=\int_{\gamma_a}(\lambda-\lambda_{\mathrm{tor}})
   =2\pi\bigl(\lambda(\xi_a)+d_a\log|y_a|\bigr).
\]
Consequently
\[
 \lambda(\xi_a)=-d_a\log|y_a|,
 \qquad
 -d_a Z_\lambda\log|y_a|=\lambda(\xi_a),
\]
which proves (iii).

We use the symplectic grading of \cite[Section~5.2]{HS}.
A compact reduction of the complex symplectic frame bundle gives
a metric $g_*$ and a complex structure $J$ with $IJ=-JI$ and
\[
 \omega_J:=\Im\Omega_B=g_*(J\,\cdot\,,\,\cdot\,),\qquad
 \omega_I:=g_*(I\,\cdot\,,\,\cdot\,).
\]
The comparison path is $I$-taming except at its initial endpoint.
For the functions used above, this follows from
\cite[Lemmas~4.3 and~4.6]{HSinv} and the identity
\[
 \dd(f(R)\beta)(v,Iv)
 =f(R)g(v,v)+f'(R)
       \bigl(\dd R(v)^2+\dd R(Iv)^2\bigr).
\]
Indeed, outside a compact set $|\dd R|_g^2\le R$ and
$f>0$, $f+Rf'>0$; the functions are constant on the remaining
compact region.

For a real two-form $\eta$, write
$\eta^\pm=(\eta\pm I^*\eta)/2$.
Then $\eta^-(v,Iv)=0$, so pointwise taming is precisely positivity
of the $(1,1)$ part $\eta^+$.
This space is convex. Interpolating the positive parts and the
anti-invariant parts compares the path pointwise with
$\omega_J+t\omega_I$, whose compatible complex structure is
\[
 J_t=\frac{J+tI}{\sqrt{1+t^2}},\qquad
 (\omega_J+t\omega_I)(v,J_tw)=\sqrt{1+t^2}\,g_*(v,w).
\]
The compact reduction identifies the determinant trivializations
along this rotation from $J$ to $I$.
At the $I$ endpoint the trivialization is $\Omega_B^m/m!$.
Thus the symplectic grading agrees under $\Phi$ with the affine
volume grading.

On each seed torus, $\Omega_B$ has constant coefficients in
logarithmic coordinates. The compact reduction there is homotopic
to a constant one. The rotation and the change
$(u,\theta)\mapsto(p,\varphi)$ of \eqref{eq:marking} are then
constant linear transformations, so the determinant trivialization
has no winding along any seed loop.
Its integer shift is fixed by the positive fiber $\varphi=0$.
Products add the standard cotangent grading and the same construction
applies to them.
\end{proof}

The following Lemma will be later used to transport Liouville structures across mutations. 

\begin{lemma}\label{lem:reseeding}
The isomorphism of Lemma \ref{lem: commute mutation and principal} respects gradings, and the geometry of Lemma~\ref{lem:principal-geometry}
can be transported to $\cU_{\pr}(B')$.
\end{lemma}
\begin{proof}
For one mutation at $k$, \cite[Lemma~4.13]{HS} gives
\[
 \lambda'_{\mathrm{tor}}-\lambda_{\mathrm{tor}}
 =-d_k\,\dd\mathcal D(-\widehat y_k),
\]
where $\widehat y_k=\prod_a z_a^{E_{ak}}$ for the current seed
coordinates $z$ and completion $E$.
Here $\mathcal D$ is the Bloch--Wigner dilogarithm,
\[
 \mathcal D(z)=\Im\operatorname{Li}_2(z)
       +\log|z|\arg(1-z),\qquad
 \operatorname{Li}_2(z)=-\int_0^z\frac{\log(1-t)}{t}\,\dd t.
\]
The same lemma extends exactness from the intersection to either
seed torus. Applying it along the mutation sequence preserves the
seed-torus condition in Lemma~\ref{lem:principal-geometry}.

For the dual coefficient change, put
$C'^\vee=D^{-1}C'D$. This is the integral lower-left block of
$-E'^{\mathsf T}$, and $\det C'^\vee=\det C'=\pm1$.
The cotangent coordinates \eqref{eq:marking} satisfy
\[
 \varphi_{\mathrm{mut}}
 =\begin{pmatrix}I_m&0\\0&C'^\vee\end{pmatrix}
       \varphi_{\mathrm{prin}}.
\]
Thus pullback of local systems fixes $x_i^\vee$ and sends the
principal frozen monodromies to
\[
 Y_i^\vee=\prod_a(y_a^\vee)^{C'^\vee_{ai}}.
\]
This is the dual coefficient change. The cotangent lift preserves
the cotangent grading and the positive fiber with its Lie group
spin structure.
Let $F:\cU_{\pr}(B)\to\cU_{\pr}(B')$ denote the isomorphism above.
The transported Liouville form and positive locus are
\[
 \lambda'=(F^{-1})^*\lambda,\qquad L_+'=F(L_+),\qquad
 \dd\lambda'=\Im\Omega_{B'},\qquad \lambda'|_{L_+'}=0.
\]
The form $\lambda'$ is finite-type Morse--Bott Weinstein and
satisfies properties (1)--(4) of Lemma~\ref{lem:principal-geometry},
with the compact torus action transported by $F$.
The seed-torus grading is the cotangent grading normalized at $L_+'$,
and $L_+'$ carries the transported Lie group spin structure.
Transporting the graded symplectomorphism of
Lemma~\ref{lem:principal-geometry} gives a graded symplectomorphism
from $(\cU_{\pr}(B'),\dd\lambda')$ to an affine K\"ahler form.
\end{proof}

\begin{lemma}\label{lem:actual-transfer}
Let $U\subset Y$ carry finite-type Liouville forms
$\lambda_U,\lambda_Y$, with independently Weinstein completions.
Let $L_+\subset Y$ be a proper Lagrangian, and suppose
\[
 \lambda_Y|_{L_+}=0,\qquad
 \lambda_U-\lambda_Y|_U=\dd f,\qquad
 f|_{L_+\cap U}=0.
\]
For every compact $Q\subset U$, a compact exact deformation
$\lambda_{Y,t}$ gives Liouville domains $D_U\subset D_Y$ with
\[
 Q\subset\operatorname{int}D_U,\qquad
 \lambda_{Y,1}=\lambda_U\text{ near }D_U,\qquad
 \lambda_{Y,t}|_{L_+}=0\quad(0\le t\le1).
\]
The domain $D_U$ completes to $(U,\lambda_U)$.
The outer domain retains the original Weinstein completion up to
compact exact deformation.
\end{lemma}
\begin{proof}
Choose a sufficiently large Liouville domain $D_U$ for $U$ containing
$Q$ in its interior. Take $\chi\in C_c^\infty(U)$ equal to one
near $D_U$, extend $\chi f$ by zero to $Y$, and put
$
 \lambda_{Y,t}=\lambda_Y+t\dd(\chi f).
$
The symplectic form is unchanged, and
$\lambda_{Y,1}=\lambda_U$ near $D_U$.
Choose a Liouville domain $D_Y$ for $Y$ containing
$D_U\cup\supp\chi$ in its interior.
The form near $\partial D_Y$ is unchanged throughout the deformation.
Moreover,
$
 \lambda_{Y,t}|_{L_+}
 =\dd_{L_+}\bigl(t(\chi f)|_{L_+}\bigr)=0.
$
Thus strict exactness holds on all of $L_+$.
\end{proof}

\begin{lemma}\label{lem:change-primitive}
Let $\lambda_0,\lambda_1$ be finite-type weakly Weinstein Liouville
forms on $X$, with a fixed grading.  Let $L_+$ be a properly embedded Lagrangian with a fixed brane structure, and with $\lambda_0|_{L_+}=0$. 
Suppose 
$\lambda_1-\lambda_0=\dd f$ for some $f$ with $f|_{L_+}=0$.

Let $\{L_{\fs}\}$ be compact exact Lagrangians with cotangent neighborhoods
as in the definition of cluster-admissibility, for both forms.
Use the same parameterizations and brane structures on these neighborhoods.
There is an equivalence $E_{01}:\Fuk(X,\lambda_0)\xrightarrow{\sim}\Fuk(X,\lambda_1)$ and a choice of isomorphism $\eta_{01}:
 E_{01}(L_+)\simeq L_+$, along with
natural isomorphisms
$
 \Vit_{\fs,1}E_{01}\simeq\Vit_{\fs,0}
$
for every $\fs$, respecting the identifications of both sides with  positive
cotangent fibers. These comparisons commute with local-system twists
pulled back from $X$.

In particular, on a smooth cluster variety any two cluster-admissible
forms for the same $\omega$ have these properties. Cluster HMS for one
therefore gives cluster HMS for the other. 
\end{lemma}

\begin{proof}
The equivalence is \cite[Theorem~1.1]{HSinv}. We check that its construction respects the stated additional
compatibilities.
Choose domains
$
 A_0\Subset B_0\Subset A_1\Subset B_1,
$
where $A_i$ completes to $(X,\lambda_0)$ and $B_i$ to
$(X,\lambda_1)$. As in \cite[Section~3]{HSinv}, choose $\chi$ equal
to zero near $\partial A_i$, and to one near $\partial B_i$, and set
$\lambda_* =\lambda_0+\dd(\chi f)$.
The Viterbo functors
\[
 \Fuk(B_1,\lambda_*)\xrightarrow{F_1}
 \Fuk(A_1,\lambda_*)\xrightarrow{F_2}
 \Fuk(B_0,\lambda_*)\xrightarrow{F_3}
 \Fuk(A_0,\lambda_*)
\]
satisfy $F_2F_1\simeq\mathrm{id}$ and
$F_3F_2\simeq\mathrm{id}$ after the compact exact deformations
identifying the two $A_i$ categories and the two $B_i$ categories.
Thus $F_2$ gives $E_{01}$.

All deformation potentials are multiples of $f$, so they vanish on
$L_+$. The deformations keep this brane strictly exact, and Viterbo
restriction preserves it by \cite[Lemma~5.8]{HS}.
For a fixed finite collection of $L_{\fs}$, choose $A_0$ to contain their
cotangent neighborhoods and take $\chi=0$ there.
The same cotangent embedding can be used for both primitives.
Indeed, if $i_{\fs}^*(\lambda_0+\dd h_{\fs,0})=\lambda_{\mathrm{can}}$,
choose a compactly supported cutoff $\rho_{\fs}=1$ near that neighborhood and put
$ h_{\fs,1}=h_{\fs,0}-\rho_{\fs} f$.
Then $h_{\fs,1}|_{L_+}=0$ and
$i_{\fs}^*(\lambda_1+\dd h_{\fs,1})=\lambda_{\mathrm{can}}$.
The same formula along each deformation keeps the cotangent target fixed.
Applying the functoriality diagram of \cite[Section~2]{HSinv} with these
cotangent domains as the innermost domains gives
$\Vit_{\fs,1}E_{01}\simeq\Vit_{\fs,0}$, including the fiber identifications.
Enlarging the outer domains gives common refinements of these diagrams.
For fixed domains the cutoffs with the prescribed boundary values form
a convex space, so their interpolations give compatible comparison
isomorphisms under these refinements.
The cofinal system of such domains therefore gives a single comparison
compatible with every $L_{\fs}$, rather than a separately chosen equivalence
for each finite collection. The fixed grading and brane structures are
transported throughout. Tensoring with a local system from $X$ commutes
with every restriction and deformation in the diagram.

It remains to check exactness of the difference for cluster-admissible
forms. Put $\delta=\lambda_1-\lambda_0$. Then $\dd\delta=0$ and
$\delta|_{L_{\fs}}$ is exact. The complex seed torus retracts onto $L_{\fs}$.
Its inclusion in the smooth cluster variety is surjective on $H_1$,
since a loop can be moved off the algebraic complement, whose real
codimension is at least two. Hence every period of $\delta$ vanishes,
so $\delta=\dd f$. Both forms vanish on the connected positive locus,
and we normalize $f|_{L_+}=0$.
\end{proof}

\subsection{Ice and mirrors}

\begin{theorem}\label{thm:freezing}
Assume that $B$ has a reddening sequence.
Assume
$\cU_{\pr}(B)$ and $\cU_{\pr}(B^\vee)$ are smooth affine of finite type.
Assume that $\cU_{\pr}(B)$ satisfies the hypotheses of
Lemma~\ref{lem:principal-geometry}, and, for the resulting symplectic form, satisfies cluster homological mirror symmetry.

Then for $I=\{1,\ldots,m\}\setminus J$, the variety
$\cU_{\pr}(B_{I,I})$ also satisfies the hypotheses of Lemma~\ref{lem:principal-geometry} and satisfies
cluster homological mirror symmetry with mirror 
$\cU_{\pr}((B^\vee)_{I,I})$.
\end{theorem}

\begin{proof}
Write
\[
 \begin{aligned}
 X&=\cU_{\pr}(B),&
 U&=X_{x_J},& Z&=\cU_{\pr}(B_{I,I}),\\
 X^\vee&=\cU_{\pr}(B^\vee),&
 U^\vee&=X^\vee_{x_J^\vee},&
 Z^\vee&=\cU_{\pr}((B^\vee)_{I,I}).
 \end{aligned}
\]
Lemmas~\ref{lem:upper-freezing} and~\ref{lem:principal-product} give
\[
 \begin{aligned}
 U&=\cU\left(\fr_J\binom B{I_m}\right)
       \cong Z\times(\C^*)^{2|J|},\\
 U^\vee&=\cU\left(\fr_J\binom{B^\vee}{I_m}\right)
       \cong Z^\vee\times(\C^*)^{2|J|}.
 \end{aligned}
\]
Smoothness and finite type descend to $Z$ and $Z^\vee$ along
the faithfully flat projections from these products.
The product formula \eqref{eq:product-omega} gives the regular
holomorphic symplectic form $\Omega_{B_{I,I},D_I}$ on $Z$.

Let $\jmath:Z\hookrightarrow U$ be the slice $x_j=s_j=1$.
Restricting the ambient Liouville form gives
\[
 \alpha_{0,Z}=\jmath^*\lambda_X,\qquad
 \dd\alpha_{0,Z}=\Im\Omega_{B_{I,I},D_I},\qquad
 (\alpha_{0,Z}-\lambda_{\mathrm{tor},Z})|_{T_Z}
       \in\dd C^\infty(T_Z).
\]
The last assertion holds on every seed of $Z$, since its mutations
are mutations in $I$ of the frozen product seed.
The positive locus $L_Z$ is proper, since its image under $\jmath$
is the closed slice $x_j=s_j=1$ in $L_+$.
The matrix $B_{I,I}$ has a reddening sequence by
\cite[Lemma~4.9]{CL}.
Lemmas~\ref{lem:quadratic} and~\ref{lem:principal-geometry}
therefore give a Liouville form $\lambda_Z$ with the asserted geometry.

For $x_j=e^{u_j+i\theta_j}$ and $s_j=e^{v_j+i\psi_j}$, give $U$
the product form
\[
 \lambda_U=\lambda_Z+
       \sum_{j\in J}d_j(v_j\,\dd\theta_j-u_j\,\dd\psi_j).
\]
The difference $\lambda_U-\lambda_X|_U$ is closed.
It is exact on the initial seed torus because, for the matrix
$T$ in Lemma~\ref{lem:dual-product}, we have $
 K_B=T^{\mathsf T}K_{\mathrm{prod}}T$ and $
 u^{\mathsf T}K_B\,\dd\theta
 =u_{\mathrm{prod}}^{\mathsf T}K_{\mathrm{prod}}
       \,\dd\theta_{\mathrm{prod}}.
$

The algebraic complement of the seed torus has real codimension
at least two. Transversality to a stratification moves every loop
in $U$ into that torus. Thus its inclusion surjects on
$H_1(-;\R)$, and every period of $\lambda_U-\lambda_X|_U$ vanishes.
Consequently we may write 
$
 \lambda_U-\lambda_X|_U=\dd f$ where $ f|_{L_+\cap U}=0. 
$
Lemma~\ref{lem:actual-transfer} gives inclusions of the Liouville domains as needed by the setup of  Theorem~\ref{thm:affine-freezing}.

Under the product identification, the positive brane restricts to
$
 E_0=\mathrm{Vit}(L_+)=L_Z\boxtimes F,
$
where $F$ is the fiber at zero of $T^*(S^1)^{2|J|}$.

For the mirror rings, put
\[
 A=\OO(X^\vee),\qquad K=\operatorname{Frac}(A),\qquad
 S=\langle x_j^\vee:j\in J\rangle,\qquad
 R=S^{-1}A.
\]
The algebra $A$ is a smooth finitely generated domain, and
\[
 R=\OO(Z^\vee)[(x_j^\vee)^{\pm1},(s_j^\vee)^{\pm1}:j\in J].
\]
Let $R_0$ be the Laurent ring of the initial frozen seed,
and let $R_i$, $i\in I$, be the Laurent rings after one mutation
at $i$. Embed them in $K$ by the dual cluster transformations.
Choose $D_U$ to contain neighborhoods of the corresponding finitely
many unitary seed tori. Their Viterbo restrictions define
$
 \rho_a:\Fuk(U)\longrightarrow\Perf(R_a)$ for $
 (a\in\{0\}\cup I).
$
Write $\Phi:\Perf(A)\xrightarrow{\sim}\Fuk(X)$ for the given
equivalence, with the fixed isomorphism $\Phi(A)=L_+$ as required by cluster HMS.
Its seed compatibility and transitivity of Viterbo restriction
\cite[Section~2]{HSinv} give
\[
 \rho_a\circ\mathrm{Vit}\circ\Phi
       \simeq R_a\otimes_A-.
\]
The columns of the frozen extended matrix
$\binom{(B^\vee)_{*,I}}{(I_m)_{*,I}}$ are linearly independent.
Equations~\eqref{eq:finite-star} and~\eqref{eq:upper-freezing} imply that
$
 R\subset R_a\subset K$ and $
 R_0\cap\bigcap_{i\in I}R_i=R$. 
This proves \ref{hyp:restriction-rings}.

To verify \ref{hyp:denominator-units}, write
$
 C=\End_{H^0\Fuk(U)}(E_0)$ and $
 f:A\longrightarrow C
$
for the map induced by $\mathrm{Vit}\circ\Phi$.
Corollary~\ref{cor: frac field embedding} gives
$A\subset C\subset K$.
The maps $h_a:C\to R_a$ induced by $\rho_a$ are $A$-linear,
hence are inclusions in $K$.
It follows that $C\subset R$.

On the cotangent factors, \eqref{eq:marking} reads
\[
 (\varphi_{x_j},\varphi_{s_j})=(-\psi_j,\theta_j),\qquad
 (p_{x_j},p_{s_j})=(d_ju_j,d_jv_j).
\]
Let $\tau_j,\sigma_j$ be the positively oriented coordinate loops
in $\varphi_{x_j},\varphi_{s_j}$, based at zero.
By \cite[Theorem~1.1]{Ab}, we have 
$
 \End(F)\simeq
 \C[\tau_j^{\pm1},\sigma_j^{\pm1}:j\in J]$ where $\tau_j, \sigma_j$ have degree zero. 
The product functor of \cite[Theorem~1.5]{GPS} gives units
$\id_{L_Z}\otimes\tau_j$ in $C$.
Lemma~\ref{lem:dual-product} computes their images in the initial
mirror chart:
$h_0(\id_{L_Z}\otimes\tau_j)=x_j^\vee=h_0(f(x_j^\vee)).
$
Injectivity of $h_0$ gives
$
 f(x_j^\vee)=\id_{L_Z}\otimes\tau_j$ and therefore $f(x_j^\vee)^{-1}=\id_{L_Z}\otimes\tau_j^{-1}$.
This proves \ref{hyp:denominator-units}.

For \ref{hyp:restriction-finite}, the algebraic variety $U$ has
the nowhere-vanishing volume form
\[
 \eta_U=
 \frac{\Omega_{B_{I,I},D_I}^{|I|}}{|I|!}
 \wedge\bigwedge_{j\in J}
       d_j(\dd\log s_j\wedge\dd\log x_j).
\]
It is logarithmic on the dense seed torus.
Lemmas~\ref{lem:log-volume} and~\ref{lem:degreewise} imply
degreewise finiteness of symplectic cohomology for an affine
K\"ahler form on $U$, graded by $\eta_U$.
Lemma~\ref{lem:principal-geometry} identifies this graded symplectic
manifold with $(U,\dd\lambda_U)$.
The symplectomorphism invariance in \cite[Theorem~1.1]{HSinv},
together with the closed--open identification
\cite[Theorem~1.1]{Ganatra}, identifies their graded symplectic
cohomology rings. Hence
$
 \SH^k(U,\lambda_U)$ is finite over $
 \SH^0(U,\lambda_U)$.

Theorem~\ref{thm:affine-freezing} now gives an equivalence
$\Phi_U:\Perf(R)\xrightarrow{\sim}\Fuk(U)$ along with
$\Phi_U(R)=E_0$ and $\mathrm{Vit}\circ\Phi\simeq\Phi_U\circ(R\otimes_A-).$
Here, the fact that we end up with a fixed identification $\Phi_U(R)=E_0$ is ultimately because we have passed the fixed cluster HMS isomorphism $\Phi(A) = L_+$ through the quotient.
The equivalence also intertwines Viterbo restriction to the chosen seed tori with the corresponding localization to Laurent polynomial rings.

Now we pass from \(U\) to \(Z\). 
Let $H_J=\Hom(\pi_1((S^1)^{2|J|},0),\C^*)$ act on $\Fuk(U)$ by
\[
 h\cdot K=K\otimes q^*\mathcal L_h.
\]
Here $q:U\cong Z\times T^*(S^1)^{2|J|}\to(S^1)^{2|J|}$
is the base projection, and $\mathcal L_h$ has monodromy $h$,
trivialized at zero. These trivializations make $E_0=L_Z\boxtimes F$
coherently equivariant.
Lemma~\ref{lem:dual-product} identifies the action on the initial mirror chart as
\[
 \begin{aligned}
 x_j^\vee&\longmapsto h(\tau_j)x_j^\vee,&
 s_j^\vee&\longmapsto h(\sigma_j)s_j^\vee &&(j\in J),\\
 x_i^\vee&\longmapsto x_i^\vee,&
 Y_i^\vee&\longmapsto Y_i^\vee &&(i\in I).
 \end{aligned}
\]
This is translation on the second factor of $U^\vee=Z^\vee\times H_J$.
Since $h_0\otimes\id_{\OO(H_J)}$ is injective, the coactions on
$\End(E_0)=R$ agree.
Lemma~\ref{lem:marked-functors} therefore makes $\Phi_U$ equivariant.
By K\"unneth and \cite[Lemma~9]{ShendeToric}, taking equivariant objects gives
\[
 \Ind\Fuk(Z)
 \simeq(\Ind\Fuk(U))^{H_J}
 \simeq\QCoh(U^\vee)^{H_J}
 \simeq\QCoh(Z^\vee).
\]
Taking compact objects yields
$\Phi_Z:\Perf(Z^\vee) \xrightarrow{\sim} \Fuk(Z)  $, and we descend the  isomorphism 
$\Phi_U(R) \cong E_0$ to
obtain a chosen isomorphism $\Phi_Z(\OO_{Z^\vee}) \cong L_Z$.

Any finite collection of seeds of $Z$ can be added to the chosen charts
using Lemma~\ref{lem:actual-transfer}.
Their restriction functors are $H_J$-equivariant by
Lemma~\ref{lem:equivariant-localization}, so the seed comparisons
also descend to $Z$.

(We do not appeal to the present result until Section \ref{sec:HMS}, so the references to Section \ref{sec: coefficients} create no circularity.) 
\end{proof}

\begin{remark}
Reddening is used to obtain the quadratic embedding for $B_{I,I}$.
We have not proved that the conclusions of Lemma~\ref{lem:principal-geometry}
pass to $B_{I,I}$ without this assumption.
\end{remark}

\begin{remark}\label{rem:rebuilding}
When an auxiliary primitive on $X$ is replaced by the form constructed
in Lemma~\ref{lem:principal-geometry}, both have differential
$\Im\Omega_B$. Assume the auxiliary primitive is also of finite type
and weakly Weinstein. If both vanish on $L_+$ and are exact on the seed tori,
the period argument in Lemma~\ref{lem:change-primitive} gives
$\lambda_1-\lambda_0=\dd f$ with $f|_{L_+}=0$.
That lemma transports HMS and its seed comparisons directly.
\end{remark}

\section{Coefficients} \label{sec: coefficients}

Here we give the proof of Theorem \ref{thm:coefficients}, and discuss other aspects of the various torus actions.  The many lengthy arguments amount to checking that various structures are compatible with various torus actions, and determining the characters of said actions in coordinates which will be useful later.  

\subsection{Equivariance}

\begin{lemma}\label{lem:marked-functors} 
Let $A$ and $B$ be ordinary (not derived) commutative $k$-algebras over any ring $k$.
A functor $F:\Perf(A)\to\Perf(B)$,
equipped with $F(A)\simeq B$, is extension of scalars for
its induced algebra map $f:A\to B$.
If $B\hookrightarrow L$ is injective, $f$ is determined by
its composite $A\to L$.
Comparisons between such functors that preserve the chosen
identifications of the generators are unique and compose coherently.

An algebraic action of an affine algebraic group $H$
flat over $k$ on
$\Mod_A$, with a coherent equivariant structure on $A$,
is similarly determined by its coaction
$A\to A\otimes\OO(H)$.
This coaction can be checked after an injection $A\hookrightarrow L$.  
\end{lemma}
\begin{proof}
By derived Morita theory \cite[Corollary~7.6]{Toen}, the functor
is determined by the right $A$-action on $F(A)\simeq B$.
Since $\End_B(B)=B$ is discrete, this is exactly $f:A\to B$;
comparisons fixing $B$ are unique. For the group action, apply
the same argument over $\OO(H)$ with the specified coherent
trivialization. The coherences are translated into a coassociativity condition, and the injection
$A\otimes\OO(H)\hookrightarrow L\otimes\OO(H)$ detects the coaction.  
\end{proof}

\begin{lemma}\label{lem:normalized-hms-unique}
Let $A$ and $A_{\fs}$ be ordinary commutative $k$-algebras, and let
$\mathcal C$ be an idempotent-complete $k$-linear stable category.
Fix $E\in\mathcal C$,  functors
$R_{\fs}:\mathcal C\to\Perf(A_{\fs})$, isomorphisms
$\kappa_{\fs}:R_{\fs}(E)\xrightarrow{\sim}A_{\fs}$, and algebra maps $f_{\fs}:A\to A_{\fs}$.
Assume that at least one $f_{\fs}$ is injective.
There is at most one tuple
\[
 \Phi:\Perf(A)\xrightarrow{\sim}\mathcal C,\qquad
 \eta:\Phi(A)\xrightarrow{\sim}E,\qquad
 \alpha_{\fs}:R_{\fs}\Phi\xRightarrow{\sim}A_{\fs}\otimes_A^{\mathbb L}-,
\]
up to unique isomorphism preserving all these data. Here  $\alpha_{\fs,A}=\kappa_{\fs}R_{\fs}(\eta)$, and the tensor
product uses $f_{\fs}$.
For a fixed $(\Phi,\eta)$, each normalized $\alpha_{\fs}$ is unique if it exists.
An unnormalized comparison can always be normalized uniquely.
\end{lemma}
\begin{proof}
Suppose two tuples exist. They give algebra isomorphisms
\[
 \rho_i:A\xrightarrow{\sim}\End_{\mathcal C}(E),\qquad
 \rho_i(a)=\eta_i\Phi_i(a)\eta_i^{-1}\quad(i=1,2).
\]
The functor $R_{\fs}$ and $\kappa_{\fs}$ give
$r_{\fs}:\End_{\mathcal C}(E)\to A_{\fs}$. The comparisons imply
\[
 r_{\fs}\rho_1=f_{\fs}=r_{\fs}\rho_2.
\]
Choose $\fs$ for which $f_{\fs}$ is injective. Since $\rho_1$ is an
isomorphism, $r_{\fs}$ is injective, so $\rho_1=\rho_2$.
Derived Morita theory now gives a unique natural isomorphism
$\tau:\Phi_1\xRightarrow{\sim}\Phi_2$ with
$\eta_2\tau_A=\eta_1$. Here $E$ is a generator and its endomorphism
algebra is the ordinary algebra $A$.
Lemma~\ref{lem:marked-functors} makes $\tau$ compatible with every
$\alpha_{\fs}$ and proves the uniqueness of each normalized comparison.
For an unnormalized comparison, its value on $A$ differs from
$\kappa_{\fs}R_{\fs}(\eta)$ by multiplication by a unique unit of $A_{\fs}$.
Multiplying the comparison by the inverse unit gives the required normalization.
\end{proof}

\begin{corollary}
    In the setting of Lemma \ref{lem:change-primitive}, if  $(X, \lambda_0)$ satisfies cluster HMS with an affine mirror, then the equivalence $E_{01}$ and related coherences are uniquely characterized by asking that their composition gives cluster HMS for $(X, \lambda_1)$. 
\end{corollary}

\begin{remark}
    We will not need it here, but the reader may be assured to know that when $H$ is an algebraic torus and $A$ is smooth, actions on $\Mod_A$ are always induced by geometric actions on $\Spec A$, even without assuming given the coherent equivariant structure \cite[Lemma~C.1]{HS}.
\end{remark}

\begin{lemma}\label{lem:equivariant-localization}
Let $j:X_0\to X$ be a Liouville subdomain inclusion for which
Viterbo restriction is defined. Let an algebraic torus $T$
parameterize a multiplicative family of rank-one local systems
$\mathcal L_t$ on $X$. Then
$
 \mathrm{Vit}(E\otimes\mathcal L_t)
 \simeq\mathrm{Vit}(E)\otimes j^*\mathcal L_t$
is compatible with the unit and multiplication in $T$.
The same holds for sectorial restriction, Orlov, product, and
continuation functors when their local systems are pulled back
from a common ambient space.

In Theorem~\ref{thm:affine-freezing}, suppose $T$ acts on $A$
by algebra automorphisms, $\Phi$ is $T$-equivariant,
$E=\Phi(A)$ is coherently equivariant, and $S$ is generated
by $T$-weight vectors.
Then $\Phi_0$ and the square \eqref{eq:affine-freezing-square}
are $T$-equivariant for the induced action on $S^{-1}A$.
\end{lemma}
\begin{proof}
The local-system factors in Floer operations and continuation
maps are parallel transports along their boundary arcs.
Tensoring and restricting commute with these transports.
Apply this construction to the universal local system over
$\OO(T)$; its tensor product gives the unit and multiplication
coherences. This is also the compatibility used in the
descent construction of \cite[Section~5.5]{HS}.

For the last assertion, $E_0=\mathrm{Vit}(E)$ is coherently
equivariant, and hence so is
\[
 A\longrightarrow\operatorname{REnd}(E_0)\simeq S^{-1}A.
\]
The localized action is determined by $t\cdot s^{-1}=(t\cdot s)^{-1}$.
Equivariant Morita theory with generator $E_0$ gives the asserted
equivalence and square. 
\end{proof}

\begin{remark}
    The assertions of Lemma \ref{lem:equivariant-localization} are compatible with the change of Liouville structure of Remark \ref{rem:rebuilding}. 
\end{remark}

\begin{lemma}
\label{lem:frozen-equivariance}
Suppose  cluster HMS holds for $\cU_{\pr}(B)$.  Let
\[
 P=\cU_{\pr}(B)\times(\C^*)^f,\qquad
 P^\vee=\mathcal{U}_{\pr}(-B^{\mathsf T})\times(\C^*)^f.
\]
Let
$\widehat H=(\C^*)^m_t\times(\C^*)^f_u$ act on $\Fuk(P)$ by local system twists, parameterized to have monodromies
$(t^{-1},u)$ in the frozen coordinates $(Z,U)$ on $P$.
The product HMS equivalence is $\widehat H$-equivariant for
the action on $P^\vee$ with weight matrix
\begin{equation}\label{eq:all-frozen-weights}
 \begin{pmatrix}
 I_m&B^{\mathsf T}&0\\
 0&0&I_f
 \end{pmatrix}
\end{equation}
on its coordinates $(X;Z;U)$. Its seed comparisons are
equivariant with the same coherences.
\end{lemma}
\begin{proof}
On the frozen torus use the universal rank-one local system
over its Laurent character ring. Parallel transport around an
integral loop is a Laurent monomial in the parameters. Tensoring
brane local systems with its pullback defines the family of
autoequivalences over $\mathcal O(\widehat H)$. Tensor products
of these universal local systems give the unit and
multiplication isomorphisms and their coherences.

The frozen coordinates take positive real values on $L_+$.
Their real logarithms give a lift to the universal cover of the
frozen torus, hence a trivialization of the pulled-back universal
local system on $L_+$. It is compatible with multiplication of
parameters. Thus the generator corresponding to
$\mathcal O_{P^\vee}$ is coherently fixed. Applying
Lemma~\ref{lem:marked-functors} over $\mathcal O(\widehat H)$
gives a coaction on $\mathcal O(P^\vee)$.

It suffices to compute this coaction on the initial seed. By
\eqref{eq:marking}, with the usual cotangent coordinates on
the extra factors, we have 
$
 \theta_Z=-\varphi_X-B^{\mathsf T}\varphi_Z$ and $\theta_U=\varphi_U$.
For a loop $\gamma$ in the seed torus, set
\[
 (a_X,a_Z,a_U)=\frac1{2\pi}\int_\gamma
       (\dd\varphi_X,\dd\varphi_Z,\dd\varphi_U)
 \in\Z^m\oplus\Z^m\oplus\Z^f.
\]
Its degrees in the frozen coordinates $(Z,U)$ on $P$ are
$(-a_X-B^{\mathsf T}a_Z,a_U)$.
Let $\mathcal L_{(t^{-1},u)}$ be the rank-one local system on
the frozen torus with monodromy parameters $(t^{-1},u)$. Then
\[
 \operatorname{mon}_{(Z,U)^*\mathcal L_{(t^{-1},u)}}([\gamma])
 =(t^{-1})^{-a_X-B^{\mathsf T}a_Z}u^{a_U}
 =t^{a_X+B^{\mathsf T}a_Z}u^{a_U}.
\]
Here $t^a=\prod_i t_i^{a_i}$, and likewise for $u$.

This is precisely
\eqref{eq:all-frozen-weights}. The principal exchange relations
are homogeneous for these weights; hence the action extends
to every seed and to their upper cluster algebra. Since
$\mathcal O(P^\vee)$ injects into the seed Laurent ring, the
computed and geometric coactions agree globally. The same
injection after tensoring with
$\mathcal O(\widehat H)^{\otimes 2}$ verifies compatibility
with multiplication. Seed restriction commutes with the
universal local-system twists by
Lemma~\ref{lem:equivariant-localization}, giving the asserted
equivariant comparisons. 
\end{proof}

\subsection{Algebra}

Here we give the elementary algebraic geometric relationships between the principal coefficient cluster algebras and general ones.  They are presumably well known.

\begin{lemma}\label{lem:coefficient-isogeny}
Let $\widetilde B=\binom BC$ have full rank, and choose
a positive rational diagonal matrix $D$ with
$BD=-DB^{\mathsf T}$ and $D^{-1}$ integral. Put
\[
 S=D^{-1}B+C^{\mathsf T}C,\qquad
 Q=\begin{pmatrix}D^{-1}B&C^{\mathsf T}\\-C&I_f\end{pmatrix}.
\]
Then $\det Q=\det S\ne0$, and
\begin{equation}\label{eq:coefficient-isogeny}
 \begin{gathered}
 \chi:(\C^*)^m_h\times(\C^*)^f_y
          \longrightarrow(\C^*)^m_Z\times(\C^*)^f_U,\\
 Z_i=\prod_jh_j^{(D^{-1}B)_{ij}}\prod_a y_a^{C_{ai}},
 \qquad U_a=\prod_jh_j^{-C_{aj}}y_a
 \end{gathered}
\end{equation}
is a finite \'etale isogeny with kernel
\[
 \ker\chi=\{(\eta,\eta^C):\eta^S=1\},\qquad
 (\eta^C)_a=\prod_i\eta_i^{C_{ai}},\quad
 (\eta^S)_i=\prod_j\eta_j^{S_{ij}}.
\]
\end{lemma}
\begin{proof}
If $Sv=0$ for $v\in\R^m$, then
$0=v^{\mathsf T}Sv=\lVert Cv\rVert^2$, since $D^{-1}B$ is
skew-symmetric. Thus $Cv=Bv=0$, and full rank gives $v=0$.
We have:
\[
 \begin{pmatrix}I_m&-C^{\mathsf T}\\0&I_f\end{pmatrix}
 \begin{pmatrix}D^{-1}B&C^{\mathsf T}\\-C&I_f\end{pmatrix}
 =\begin{pmatrix}S&0\\-C&I_f\end{pmatrix}.
\]
The leftmost matrix has determinant $1$, so $\det Q=\det S$.
The exponent matrix of $\chi$ is $Q$.
Setting $U=1$ gives $y=h^C$, after which $Z=h^S$.
\end{proof}

\begin{lemma}\label{lem:coefficient-fiber}
Let $\widetilde B=\binom BC$ have full rank, and put
\[
 P=\cU_{\pr}(B)\times(\C^*)^f,\qquad
 F:P\longrightarrow(\C^*)^m,\qquad
 F_i=Z_i\prod_a U_a^{-C_{ai}}.
\]
Here $(X;Z)$ are principal seed coordinates and $U$ are the
coordinates on the extra torus. There is an isomorphism
\begin{equation} \label{eq:ys}
 \iota_C:\cU(B;C)\xrightarrow{\sim}F^{-1}(1),\qquad
 (x;y)\longmapsto(x;y^{C^{\mathsf T}};y),\qquad
 (y^{C^{\mathsf T}})_i=\prod_a y_a^{C_{ai}}.
\end{equation}
If $\cU_{\pr}(B)$ is smooth affine of finite type, so is
$\cU(B;C)$.
\end{lemma}
\begin{proof}
Choose $D$ as in Lemma~\ref{lem:coefficient-isogeny}.
Pull back $P$ along \eqref{eq:coefficient-isogeny} and rescale
$X_i=h_i^{d_i^{-1}}x_i$.
Lemma~\ref{lem:coefficients} applies to $(B;I_m;0)$, with
coefficient-change matrices $\binom{D^{-1}}0$ and $Q^{\mathsf T}$.
The new frozen rows are
\[
 \binom{D^{-1}}0B+Q^{\mathsf T}\binom{I_m}0=\binom0C.
\]
Thus the pullback is $(\C^*)^m_h\times\cU(B;C)$, with
$F_i=\prod_jh_j^{S_{ij}}$.
The kernel in Lemma~\ref{lem:coefficient-isogeny} acts freely
and transitively on the possible $h$ coordinates with $h^S=1$.
Each orbit has a unique representative with $h=1$, giving $\iota_C$.
Smoothness and finite type pass to this finite \'etale pullback
and descend along its projection to $\cU(B;C)$.
\end{proof}

\begin{lemma}\label{lem:frozen-parameters}
Let $C\in\operatorname{Mat}_{f\times m}(\Z)$, and put
$\widehat H=(\C^*)^m_t\times(\C^*)^f_u$ and $H=(\C^*)^m$.
There is a split exact sequence
\begin{equation}\label{eq:frozen-torus-extension}
 1\longrightarrow H\xrightarrow{\ h\mapsto(h,h^C)\ }\widehat H
 \xrightarrow{\ (t,u)\mapsto u\,t^{-C}\ }
 T_{\fr}(B;C)\longrightarrow1,
 \qquad \kappa\longmapsto(1,\kappa).
\end{equation}
Restricting the weights \eqref{eq:all-frozen-weights} along
$H\to\widehat H$ gives $[I_m\mid B^{\mathsf T}\mid C^{\mathsf T}]$.
\end{lemma}
\begin{proof}
The map $(t,u)\mapsto u\,t^{-C}$ has kernel
$\{(h,h^C):h\in H\}$ and the displayed section.
The weight calculation is
\[
 \begin{pmatrix}I_m&C^{\mathsf T}\end{pmatrix}
 \begin{pmatrix}I_m&B^{\mathsf T}&0\\0&0&I_f\end{pmatrix}
 =\begin{pmatrix}I_m&B^{\mathsf T}&C^{\mathsf T}\end{pmatrix}.
\]
\end{proof}

\begin{proposition}\label{prop:upper-finite-atlas}
Let $\widetilde B=\binom BC$ have full rank, with $B$
skew-symmetrizable. Choose
\begin{equation}\label{eq:finite-atlas-matrix}
 T=\begin{pmatrix}-B^{\mathsf T}&-C^{\mathsf T}\\ C'&T_{ff}\end{pmatrix},
 \qquad \det T\ne0,\qquad G_T=\Hom(\coker T,\C^*).
\end{equation}
Then there is a finite \'etale atlas and a quotient presentation
\begin{equation}\label{eq:finite-atlas}
 \cU(-B^{\mathsf T};C')\longrightarrow\mathcal X(\widetilde B^\vee),
 \qquad
 \mathcal X(\widetilde B^\vee)
 \cong[\cU(-B^{\mathsf T};C')/G_T].
\end{equation}
For the standard basis $e_1,\ldots,e_{m+f}$, the action on seed
coordinates $(u;v)$ is
\[
 g\cdot u_i=g([e_i])u_i,\qquad
 g\cdot v_a=g([e_{m+a}])v_a.
\]
If $\cU_{\pr}(-B^{\mathsf T})$ is smooth affine of finite type,
so is the atlas, and the quotient is a smooth separated
Deligne--Mumford stack with affine coarse space.
If, in addition, the column lattice of $\widetilde B$ is saturated
in $\Z^{m+f}$, then $T$ can be chosen unimodular and the quotient
is a smooth affine scheme.
\end{proposition}
\begin{proof}
The independent rows $[-B^{\mathsf T}\mid-C^{\mathsf T}]$ extend
to an integral square matrix of nonzero determinant.
The quotient presentation and the displayed group action are
\cite[Proposition~3.10]{HS}.
The same proof applies verbatim to upper cluster varieties.
Its coefficient change uses $Q=T$ and $P=-\binom{I_m}{0}$,
giving
\[
 -\binom{I_m}{0}(-B^{\mathsf T})+T\binom{I_m}{0}
 =\binom0{C'}.
\]
The finite cover in that proof is therefore
$(\C^*)^m\times\cU(-B^{\mathsf T};C')$.
Its commuting quotients by $(\C^*)^m$ and $G_T$ give
\eqref{eq:finite-atlas}.

The matrix $\binom{-B^{\mathsf T}}{C'}$ consists of the first
$m$ columns of $T$, so it has full rank.
Lemma~\ref{lem:coefficient-fiber} gives smoothness and finite
type of the atlas.
The quotient by the finite group $G_T$ is separated
Deligne--Mumford, with coarse space
$\Spec\U(-B^{\mathsf T};C')^{G_T}$.
If the column lattice is saturated, the top rows of $T$ span
a primitive sublattice and extend to an integral basis.
Then $\det T=\pm1$ and $G_T=1$.
For $C=I_m$, the choice $T=E_B^\vee$ gives
$\mathcal X\bigl(\binom B{I_m}^{\!\vee}\bigr)
\cong\cU_{\pr}(-B^{\mathsf T})$.
\end{proof}

\subsection{Geometry}

\begin{lemma}\label{lem:coefficient-geometry}
In Lemma~\ref{lem:coefficient-fiber}, suppose $\cU_{\pr}(B)$ is
smooth affine of finite type and carries a Liouville form
$\lambda_{\pr}$ satisfying the conclusions of
Lemma~\ref{lem:principal-geometry}, with $D^{-1}$ integral. Put
\[
 \lambda_P=\lambda_{\pr}+\sum_a\log|U_a|\,\dd\arg U_a,
 \qquad \alpha=\iota_C^*\lambda_P.
\]
Then $\alpha$ is finite-type Morse--Bott Weinstein, with proper
strictly exact positive locus.
With $S=D^{-1}B+C^{\mathsf T}C$, define
\[
 \begin{gathered}
 \rho_{s,q}(X,Z,U)=
 \left(e^{D^{-1}(s+\mathrm iq)}X,
       e^{D^{-1}B(s+\mathrm iq)}Z,e^{C(s-\mathrm iq)}U\right),\\
 \rho:\cU(B;C)\times T^*\R^m\longrightarrow P,\qquad
 (x,q,p)\longmapsto\rho_{-(S^{\mathsf T})^{-1}p,q}(\iota_C(x)).
 \end{gathered}
\]
Here, $\rho$ is a projection of an isomorphism
\[
 \cU(B;C)\times T^*\R^m
 \xrightarrow{\sim} P\times_{(\C^*)^m}\C^m,\qquad
 (x,q,p)\longmapsto\bigl(\rho(x,q,p),p+\mathrm iSq\bigr),
\]
where $\C^m\to(\C^*)^m$ is the exponential map.  The deck group $\Gamma= \pi_1((\C^*)^m)=\Z^m$ acts by 
\[
 g\cdot(x,q,p)=
 \left(\rho_{0,-2\pi S^{-1}g}(x),q+2\pi S^{-1}g,p\right),
\]
where $x$ is identified with $\iota_C(x)$.

There is a smooth function $G(x,p)$ with $G(x,0)=0$ such that
\[
 \rho^*\lambda_P=\alpha+\sum_i p_i\,\dd q_i+\dd G(x,p),
\]
and hence
\begin{equation}\label{eq:split-real-form}
 \rho^*\dd\lambda_P=\dd\alpha+\sum_i\dd p_i\wedge\dd q_i.
\end{equation}

There is a smooth family $H_\tau\in C^\infty_c(P)$, $0\le\tau\le1$,
with uniformly compact support, $H_0=0$, and $H_\tau|_{L_{P,+}}=0$,
and a $\Gamma$-covering
\[
 \Pi:\cU(B;C)\times T^*\R^m\longrightarrow P,
 \qquad
 \Pi^*(\lambda_P+\dd H_1)=\alpha+\sum_i p_i\,\dd q_i.
\]
It has the same deck action as $\rho$ and is graded for the product of the
induced grading on $\cU(B;C)$ and the standard cotangent grading on the other factor.
The preimage of the positive locus is
\[
 \begin{gathered}
 \Pi^{-1}(L_{P,+})=
 \coprod_{g\in\Z^m}
 \rho_{0,-q_g}(L_{\cU(B;C),+})\times T^*_{q_g}\R^m,\\
 q_g=2\pi S^{-1}g.
 \end{gathered}
\]
Transport along this deformation preserves the positive branes with their induced gradings,
spin structures, and trivial local systems.
\end{lemma}

\begin{proof}
Lemma~\ref{lem:principal-geometry}(iii) and the standard Liouville
form on the extra torus give
\[
 Z_{\lambda_P}\log|Z_i|=\log|Z_i|,
 \qquad Z_{\lambda_P}\log|U_a|=\log|U_a|.
\]
Since $\log|F_i|=\log|Z_i|-\sum_a C_{ai}\log|U_a|$, we obtain
\[
 Z_{\lambda_P}\log|F_i|
 =\log|Z_i|-\sum_a C_{ai}\log|U_a|
 =\log|F_i|.
\]
The compact action $\rho_{0,q}$ preserves $\lambda_P$ and its
product Lyapunov function. Together with
Lemma~\ref{lem:coefficient-fiber}, these are the inputs of
\cite[Lemma~4.14(ii)--(iv)]{HS}, which gives $\alpha$ and the
Liouville-form identity above.
The grading is the one induced from the principal grading and
the standard cotangent grading in \eqref{eq:split-real-form}.
The formula for the deck transformations is a direct calculation from the definitions.

\cite[Lemma~4.15]{HS} compares completions after a compact exact
deformation. To specify the maps, take the finite cover
\[
 \widehat P=\{(z,h)\in P\times(\C^*)^m:F(z)=h^S\},
 \qquad \pi(z,h)=z.
\]
The map $\rho$ factors through the diffeomorphism
\[
 \Phi:\cU(B;C)\times T^*(\R^m/2\pi\Z^m)\longrightarrow\widehat P,
 \qquad
 \Phi(x,q,p)=\bigl(\rho(x,q,p),e^{S^{-1}p+\mathrm iq}\bigr).
\]
Its deck group is $\Delta=\Z^m/S\Z^m$.
Choose a $\Delta$-invariant compactly supported cutoff $\chi(x,p)$,
equal to $1$ near $\Skel(Z_\alpha)\times\{0\}$, and set
$
 H_\tau\circ\pi\circ\Phi=-\tau\chi G.
$
The invariance of $G$ and its vanishing on the lifted positive
loci follow from \cite[Lemma~4.14(iv)]{HS}. Thus $H_\tau$ descends
to $P$, has uniformly compact support, and vanishes on $L_{P,+}$.
Moreover,
\[
 \Phi^*\pi^*(\lambda_P+\dd H_\tau)
 =\alpha+\sum_i p_i\,\dd q_i+\dd\bigl((1-\tau\chi)G\bigr).
\]

By \cite[Lemma~4.15]{HS}, at $\tau=1$ this form and
$\alpha+\sum_i p_i\,\dd q_i$ complete a common compact Liouville
domain $D_0\subset\{\chi=1\}$.
Comparison of their Liouville flows gives a $\Delta$-equivariant
diffeomorphism $\Theta$ to $\widehat P$, equal to $\Phi$ on $D_0$,
with
$
 \Theta^*\pi^*(\lambda_P+\dd H_1)
 =\alpha+\sum_i p_i\,\dd q_i.
$

Lifting this comparison to $q\in\R^m$ gives $\Pi$.
Since $D_0$ is a deformation retract of either completion,
the comparison induces the identity on fundamental groups.
Hence the covering still corresponds to
$
 \ker\bigl(F_*:\pi_1(P)\longrightarrow
                  \pi_1((\C^*)^m)=\Z^m\bigr)
$
and has the stated $\Gamma$-action.

The comparison is the identity on $D_0$, so it preserves the
induced grading there and therefore on the completion.
Both Liouville forms vanish on the lifted positive loci, and
their flows preserve those loci. Thus the comparison preserves
the positive branes, with their transported gradings and spin
structures. Lemma~4.14(iv) of \cite{HS} gives the displayed
preimage of $L_{P,+}$, and pullback preserves its trivial local system.
\end{proof}

\begin{proposition}\label{prop:coefficient-affine}
Under the hypotheses of Lemma~\ref{lem:coefficient-geometry},
assume a closed algebraic embedding
$\iota_0:\cU_{\pr}(B)\hookrightarrow\C^N$ whose coordinates
are weight vectors for the complexification of
\eqref{eq:principal-compact-action} and satisfy
$
 \{\iota_{0,\mu},\iota_{0,\nu}\}_{\Omega_B}
 \in\operatorname{span}_{\C}
       \{1,\iota_{0,\rho},\iota_{0,\rho}\iota_{0,\sigma}\}.$

There is an affine Liouville form $\laff$, a diffeomorphism
$\phi$ isotopic to the identity, and a smooth function $h$ such that
$
 \phi^*\alpha=\laff+\dd h.
$

The symplectomorphism intertwines the grading of
Lemma~\ref{lem:coefficient-geometry} with the grading defined
by the holomorphic volume form with seed expression
$
 \bigwedge_{i=1}^m\dd\log x_i\wedge
 \bigwedge_{a=1}^f\dd\log y_a.
$ 
\end{proposition}
\begin{proof}
If $B$ has a reddening sequence, Lemma~\ref{lem:quadratic}
supplies the embedding $\iota_0$ assumed above.
Use the embedding $\iota=(\iota_0,U,U^{-1})|_{F^{-1}(1)}$.
For the algebraic complex structure $I$,
\[
 \dd\alpha(v,Iv)
 =\sum_a\left(\dd\log|y_a|(v)^2+\dd\arg y_a(v)^2\right)\ge0.
\]
Choose a strictly plurisubharmonic Morse function $R$ agreeing
with $1+|\iota|^2/4$ outside a compact set, and put
$\laff=-\dd R\circ I$ and $g=\dd\laff(\,\cdot\,,I\,\cdot\,)$.
The calculation in \cite[proof of Proposition~B.6]{HS} uses
the quadratic coordinate brackets and the splitting
\eqref{eq:split-real-form}. Its correction bivectors are products
of infinitesimal torus actions, which are linear in the chosen
weight-vector coordinates. It therefore gives
\[
 \| (\dd\alpha)^{-1}\|_g\le c_0(1+|\iota|^2)\le c_1R
 \quad\text{outside a compact set}.
\]
Apply \cite[Theorem~1.4]{HSinv} with growth function $r\mapsto r$.

The volume form on the fiber is the nowhere-zero residue
\[
 \nu=\operatorname{Res}_{F_1=\cdots=F_m=1}
 \frac{(\Omega_B^m/m!)\wedge\bigwedge_a\dd\log U_a}
      {\prod_i(F_i-1)}.
\]
Its seed expression is the stated logarithmic volume, up to a
nonzero constant. On that seed, $\dd\alpha$ and
\eqref{eq:split-real-form} have constant matrices in logarithmic
coordinates. The comparison path of \cite[Section~4]{HSinv}
is $I$-taming away from its initial endpoint.
The convexity argument in the proof of
Lemma~\ref{lem:principal-geometry} compares its grading transport
with a constant linear comparison on the seed.
The resulting grading agrees with that of $\nu$, normalized at
the positive locus. The restriction
$H^1(\cU(B;C);\Z)\to H^1((\C^*)^{m+f};\Z)$ is injective,
since every loop can be moved off the algebraic complement of
the seed. Hence the gradings agree globally.
\end{proof}

\subsection{Descent}

\begin{lemma}\label{lem:character-duality}
Let $\Gamma\cong\Z^r$ act on a stable presentable $\C$-linear category
$\mathcal D$, and put $H=\Hom(\Gamma,\C^*)$.
The category $\mathcal D^\Gamma=\lim_{B\Gamma}\mathcal D$ carries
an algebraic $H$-action, and
$(\mathcal D^\Gamma)^H\simeq\mathcal D.$

For $h\in H$, viewed as a character $h:\Gamma\to\C^*$, the action is
\[
 h\cdot(E,(u_g)_{g\in\Gamma})
   =(E,(h(g)^{-1}u_g)_{g\in\Gamma}),
 \qquad u_g:gE\xrightarrow{\sim}E.
\]
Here the $u_g$ are the descent isomorphisms of a $\Gamma$-equivariant object.
The character identity $h(g+g')=h(g)h(g')$ preserves their cocycle condition.

The equivalence is natural for $\Gamma$-equivariant functors.
If an algebraic torus $K$ acts on $\mathcal D$ and commutes
coherently with $\Gamma$, the equivalence is $K$-equivariant.
\end{lemma}
\begin{proof}
This is categorical Pontryagin duality
\cite[Lemma~9]{ShendeToric}, with the character parameter inverted.
Its naturality applies to the functors and coherence isomorphisms
of the commuting $K$-action.
\end{proof}

\begin{theorem}\label{thm:coefficients}
Suppose $\cU_{\pr}(B)$ and $\cU_{\pr}(-B^{\mathsf T})$ are
smooth affine of finite type, and  cluster mirror
symmetry is realized by a Liouville form $\lambda_{\pr}$ satisfying
the conclusions of Lemma~\ref{lem:principal-geometry}, with the
principal completion \eqref{eq:completion}.
Then, for every $C\in\operatorname{Mat}_{f\times m}(\Z)$ such that
$\binom BC$ has full rank, there is a
$T_{\fr}(B;C)$-equivariant cluster mirror equivalence
\[
 \Fuk(\cU(B;C))\simeq
 \Perf\left([P^\vee/H]\right),\qquad
 P^\vee=\cU_{\pr}(-B^{\mathsf T})\times(\C^*)^f,
 \quad H=(\C^*)^m,
\]
where $H$ acts with weights \eqref{eq:weights}.
\end{theorem}
\begin{proof}
Choose an integer $N>0$ such that $ND^{-1}$ is integral,
and replace $(D,\lambda_{\pr})$ by $(D/N,\lambda_{\pr}/N)$.
Let $\varphi_s$ be the complete flow of $Z_{\lambda_{\pr}}$.
Since $\iota_{Z_{\lambda_{\pr}}}\lambda_{\pr}=0$, Cartan's formula gives
\[
 \mathcal L_{Z_{\lambda_{\pr}}}\lambda_{\pr}=\lambda_{\pr},
 \qquad \varphi_s^*\lambda_{\pr}=e^s\lambda_{\pr},
 \qquad \varphi_{\log N}^*(\lambda_{\pr}/N)=\lambda_{\pr}.
\]
Thus $\varphi_{\log N}$ is an exact symplectomorphism from the original
Liouville manifold to the rescaled one.  (The point of this rescaling is to ensure $(D/N)^{-1}=ND^{-1}$ is integral, as required in
Lemmas~\ref{lem:coefficient-isogeny} and~\ref{lem:coefficient-geometry}.)

Let $P$, $F$, and $\alpha$ be as in
Lemmas~\ref{lem:coefficient-fiber} and~\ref{lem:coefficient-geometry}.
The covering in the latter lemma gives
\[
 \Ind\Fuk(P)\simeq
 \bigl(\Ind\Fuk(\cU(B;C))\bigr)^\Gamma,
 \qquad \Gamma=\Z^m,
\]
by descent \cite[Section~5.5]{HS}.
Concretely, one applies sectorial descent to product charts
over a finite cell decomposition of $\R^m/(2\pi S^{-1}\Z^m)$
and uses K\"unneth on their lifts.
Lemma~\ref{lem:character-duality} then gives
$
 \Ind\Fuk(\cU(B;C))\simeq(\Ind\Fuk(P))^H.
$

For $h\in H=\Hom(\Z^m,\C^*)$, let $\mathcal L_{h^{-1}}$ be
the rank-one local system on $(\C^*)^m$ with character
$a\mapsto h^{-a}=\prod_i h_i^{-a_i}$. The action on
$\Ind\Fuk(P)$ is
\[
 E\longmapsto E\otimes F^*\mathcal L_{h^{-1}}.
\]

Principal HMS and K\"unneth give
$\Ind\Fuk(P)\simeq\QCoh(P^\vee)$.
For $(t,u)\in\widehat H=(\C^*)^m\times(\C^*)^f$, let
$\mathcal L_{(t^{-1},u)}$ be the rank-one local system on
$(\C^*)^m_Z\times(\C^*)^f_U$ with character
\[
 \Z^m\oplus\Z^f\longrightarrow\C^*,
 \qquad (a,b)\longmapsto t^{-a}u^b.
\]
Lemma~\ref{lem:frozen-equivariance} identifies the action
$E\mapsto E\otimes(Z,U)^*\mathcal L_{(t^{-1},u)}$ with the
action on $P^\vee$ of weights \eqref{eq:all-frozen-weights}.
Here $(Z,U):P\to(\C^*)^{m+f}$ is the map given by the frozen
coordinates of the principal factor and the extra torus.
For a loop $\gamma$ in $P$ whose image under $(Z,U)$ has
winding numbers $(a,b)$, the monodromy is given by 
$
 \operatorname{mon}_{F^*\mathcal L_{h^{-1}}}([\gamma])
 =h^{-a+C^{\mathsf T}b}=h^{-a}(h^C)^b.$
Equivalently,
\[
 F^*\mathcal L_{h^{-1}}
 \simeq (Z,U)^*\mathcal L_{(h^{-1},h^C)},
 \qquad
 (h^{-1},h^C)\in\Hom(\Z^{m+f},\C^*).
\]
Thus the subgroup used in descent is $h\mapsto(h,h^C)$.
Its weights are $[I_m\mid B^{\mathsf T}\mid C^{\mathsf T}]$
by Lemma~\ref{lem:frozen-parameters}.

Taking $H$-equivariant objects gives
\begin{equation}\label{eq:residual-equivariant-hms}
 \Ind\Fuk(\cU(B;C))\simeq\QCoh([P^\vee/H]).
\end{equation}

The section $\kappa\mapsto(1,\kappa)$ of
\eqref{eq:frozen-torus-extension} acts on $\Ind\Fuk(P)$ by
$E\mapsto E\otimes U^*\mathcal L_\kappa$.
Its restriction to $F^{-1}(1)$ is
\[
 E\longmapsto E\otimes y^*\mathcal L_\kappa,
 \qquad y=U\circ\iota_C:\cU(B;C)\longrightarrow(\C^*)^f.
\]
This is the $T_{\fr}(B;C)$-action specified in the introduction.

Lemma~\ref{lem:equivariant-localization} makes the descent
comparisons on the product charts and their overlaps equivariant
for these twists. Naturality in Lemma~\ref{lem:character-duality}
then makes \eqref{eq:residual-equivariant-hms}
$T_{\fr}(B;C)$-equivariant.

The same constructions on seed charts give the equivariant
seed comparisons. Positive real logarithms give canonical trivializations for the restrictions of the local system to $L_+$, compatibly with multiplication, so
$L_+$ descends to $\OO$, as in \cite[Section~5.5]{HS}.

We now compare with the parameterized seed tori of
Definition~\ref{def: cms}.  
Fix finitely many seeds. In the finite cover used in
Lemma~\ref{lem:coefficient-geometry}, choose a compact
$\Delta$-invariant Liouville domain $D_0$ for
$\alpha+\sum_i p_i\,\dd q_i$ containing neighborhoods of their lifted
unitary tori. Such a domain is a sufficiently large regular sublevel
of the product Lyapunov function.
Choose the cutoff $\chi(x,p)$ equal to $1$ near $D_0$, and set
$H_\tau\circ\pi\circ\Phi=-\tau\chi G$ as before. Then
\[
 \Phi^*\pi^*(\lambda_P+\dd H_1)
   =\alpha+\sum_i p_i\,\dd q_i
 \quad\text{near }D_0.
\]
The skeleton computation in \cite[proof of Lemma~4.15]{HS} is unchanged
by this choice of cutoff. Both Liouville forms therefore complete
$D_0$. Their flow comparison satisfies $\Theta|_{D_0}=\Phi|_{D_0}$.  Thus the $\Pi$ of 
Lemma~\ref{lem:coefficient-geometry} agrees with the explicit $\rho$ here on the lift of $D_0$.
The potentials vanish on the positive loci.
Lemma~\ref{lem:change-primitive} identifies the seed restrictions before
and after the compact exact deformation, with the same positive fibers.

For a seed $\fs$, let $L_{P,\fs}$ be the unitary product seed torus in $P$,
and let $P_{\fs}^\vee\subset P^\vee$ be its mirror seed torus.
On these neighborhoods, the covering map is $\rho$.
Sectorial descent and Lemma~\ref{lem:character-duality}, applied to
this restriction, give
\[
 \bigl(\Ind\Fuk(T^*L_{P,\fs})\bigr)^H
 \simeq\Ind\Fuk(T^*L_{\fs}),\qquad
 \QCoh(P_{\fs}^\vee)^H\simeq\QCoh([P_{\fs}^\vee/H]).
\]
These equivalences commute with the restrictions from $P$ and $P^\vee$.

The quotient coordinates can be checked on integral loops.
In the initial seed, \eqref{eq:ys} gives
\[
 (\iota_C)_*:\Z^m\oplus\Z^f\longrightarrow
       \Z^m\oplus\Z^m\oplus\Z^f,\qquad
 (a,c)\longmapsto(a,C^{\mathsf T}c,c),
\]
where the target records winding in $(X,Z,U)$ on $P$.
Under principal HMS, the coordinate loops in $x_i$ and $y_a$ thus give
the following characters on $P^\vee$, in the coordinates of
\eqref{eq:weights},
\[
 Z_i\prod_jX_j^{-B_{ij}}=\chi_i,\qquad
 U_a\prod_jX_j^{-C_{aj}}=\chi_{m+a}.
\]
They form a basis of the invariant character lattice.
The monomial coefficient changes in Lemma~\ref{lem:coefficients}
give the corresponding maps on loops in every mutated seed.
The invariant characters are the $\chi_{\fs,i}$ of
\cite[Section~3.6, equations~(25) and~(28)]{HS}.
Hence $[P_{\fs}^\vee/H]=\Spec A_{\fs}$ with
$j_{\fs}^*\chi_{\fs,i}=t_i$, as required in the introduction.

Enlarging $D_0$ treats any finite collection of seeds simultaneously.
The comparisons under enlargement and changes of cutoff are those of
Lemma~\ref{lem:change-primitive}.
They commute with local-system twists, so taking $H$-equivariant
objects gives comparisons for one global descended equivalence.
Let $\Phi_C:\Perf([P^\vee/H])\to\Fuk(\cU(B;C))$ be its inverse.
The component $q=0$ of the lifted positive locus in
Lemma~\ref{lem:coefficient-geometry} is the product of $L_+$ with the
cotangent fiber over $0$.
Positive real logarithms give its trivial equivariant structure.
Thus descent supplies $\eta_C:\Phi_C(\OO)\xrightarrow{\sim}L_+$, along with 
$\alpha_{\fs}:\chi_{\fs}\Vit_{\fs}\Phi_C\xRightarrow{\sim}j_{\fs}^*$ such that $ \alpha_{\fs,\OO}=\kappa_{\fs}\circ\chi_{\fs}\Vit_{\fs}(\eta_C)$.
The comparisons are $T_{\fr}(B;C)$-equivariant by
Lemma~\ref{lem:equivariant-localization}.

Taking compact objects in \eqref{eq:residual-equivariant-hms}
gives the stated result. 
\end{proof}

\subsection{Seed compatibility}

We write
\[
 \Nn := \cU_{\pr}([0]) = \{(a,b,z)\in\C^2\times\C^*:ab=1+z\},
\]
and write $\Nn^\vee$ for its dual principal model.
Its seed charts are $D(a)$ and $D(b)$.
Recall we write  $\mathrm{Vit}_{\fs}:\Fuk(\cU(\widetilde B))\to\Fuk(T^*L_{\fs})$
for Viterbo restriction.
For principal coefficients we use
$\beta_{\fs}^{\cA}\mathrm{Vit}_{\fs}^{\cA}$ and the charts $j_{\fs}^{\cA}$
from Section~\ref{sec:mirror-coordinates}.

Homological mirror symmetry for $\cN$ is well known, and the fact that it respects cluster structures has long been in the literature \cite{seidel-categorical-dynamics,  GLe, PT}.  A version compatible with our present definition of cluster mirror symmetry is available in \cite{HS}; specifically the form is constructed in \cite[Section~4.2, Lemma~4.7]{HS}, and cluster HMS is the
rank-one case of \cite[Theorem~1.3 and Sections~5.3--5.4]{HS}. 
Lemma~\ref{lem:change-primitive} lets us
use the form of Lemma~\ref{lem:principal-geometry}.

On principal seed tori we use the parameterization in \eqref{eq:marking}, corresponding to $\cA$ coordinates on $\cU_{\pr}(B^\vee)$, in order to appeal to the 
expression for the freezing product  given in Lemma~\ref{lem:principal-product}.
The cotangent change of basis in Section~\ref{sec:mirror-coordinates}
converts these comparisons into \eqref{eq:cms}.

\begin{proposition}\label{prop:node-rank-one}
Every full-rank seed with one mutable variable satisfies
equivariant cluster mirror symmetry.
\end{proposition}
\begin{proof}
Write the frozen column as $c\in\Z^f\setminus\{0\}$ and put
$c_\pm=(\max(\pm c_i,0))_i$. The exchange relation is
\[
 xx'=y^{c_+}+y^{c_-},\qquad b=x'y^{-c_-},\qquad xb=1+y^c.
\]
Let $R=\C[y_1^{\pm1},\ldots,y_f^{\pm1}]$ and $p=1+y^c$.
For a Laurent polynomial $\sum_n r_nx^n\in R[x^{\pm1}]$, membership in
$R[(p/x)^{\pm1}]$ is equivalent to $r_{-n}\in p^nR$ for every $n>0$.
Consequently
\[
 \U\binom0c=R[x^{\pm1}]\cap R[(p/x)^{\pm1}]
 =R[x,p/x]\cong R[x,b]/(xb-1-y^c).
\]
Thus $X_c:=\cU\binom0c$ is the pullback of $\Nn$
along $(\C^*)^f\to\C^*$, $y\mapsto y^c$.
Its two seed charts are $D(x)$ and $D(b)$.

Both principal models are $\Nn\cong D(ab-1)\subset\C^2$,
so they are smooth affine of finite type.
Apply Theorem~\ref{thm:coefficients} with $B=(0)$ and $C=c$ to obtain
\[
 \Phi_c:\Fuk(X_c)\xrightarrow{\sim}
 \Perf\bigl([(\Nn^\vee\times(\C^*)^f)/\C^*]\bigr).
\]
For $h\in\C^*$, the quotient action is
\[
 h\cdot(a^\vee,b^\vee,z^\vee,u_i)
 =(ha^\vee,h^{-1}b^\vee,z^\vee,h^{c_i}u_i).
\]
The residual torus $(\C^*)^f$ scales the $u_i$ and corresponds to
local-system twists in the frozen variables on $X_c$.
For either seed $\fs$, let $j_{\fs}$ be its dual chart inclusion.
The seed comparisons in Theorem~\ref{thm:coefficients} give
$j_{\fs}^*\Phi_c\simeq\chi_{\fs}\mathrm{Vit}_{\fs}$ and $
 \Phi_c(L_{c,+})\simeq\OO$.
\end{proof}

\begin{lemma}\label{lem:one-seed-hms}
Let $X=\cU_{\pr}(B)$ and $X^\vee=\cU_{\pr}(B^\vee)$ be smooth affine
of finite type. Let $\lambda_X$ satisfy the conclusions of
Lemma~\ref{lem:principal-geometry}.
Let $\Phi:\Fuk(X)\xrightarrow{\sim}\Perf(X^\vee)$ carry $L_+$ to $\OO_{X^\vee}$.
For each seed $\fs$, write $j_{\fs}^{\cA}:T_{\fs}^{\cA}\hookrightarrow X^\vee$
for its dual $\cA$ chart inclusion.
If $(j_{\fs}^{\cA})^*\Phi\simeq\beta_{\fs}^{\cA}\mathrm{Vit}_{\fs}^{\cA}$ for one seed $\fs$,
then this holds for every seed.
\end{lemma}
\begin{proof}
Suppose the comparison holds at $\fs$, and mutate at $k$ to obtain $\fs'$.
By Lemma~\ref{lem:reseeding}, we may take $\fs$ as the original seed.
Set $J=\{1,\ldots,m\}\setminus\{k\}$.
Lemmas~\ref{lem:upper-freezing} and~\ref{lem:principal-product} give
\[
 U:=X_{x_J}\cong\Nn\times(\C^*)^{2|J|},\qquad
 U^\vee:=X^\vee_{x_J^\vee}\cong\Nn^\vee\times(\C^*)^{2|J|}.
\]
The first product uses
\[
 a=x_k,\qquad Y=y_k\prod_{j\in J}x_j^{B_{jk}},\qquad
 M_-:=\prod_{j\in J}x_j^{[-B_{jk}]_+},\qquad b=x'_k/M_-.
\]
Here $[r]_+=\max(r,0)$ and $ab=1+Y$.
The other coordinates are $(x_j,s_j)_{j\in J}$ from
\eqref{eq:principal-freezing-coordinates}.
The maps are monomial on each seed chart, so they identify both
unitary seed tori and their positive loci with the corresponding product loci.

Lemma~\ref{lem:dual-product} computes the dual coordinates on the first chart.
On the second chart, $\dd\log b=\dd\log x'_k-\sum_{j\in J}[-B_{jk}]_+\dd\log x_j$.
The cotangent calculation \eqref{eq:angle-map}, applied to this monomial map
with the mutated completions, gives
\[
 b^\vee=(x_k^\vee)'\prod_{j\in J}
 (x_j^\vee)^{-(d_k/d_j)[-B_{jk}]_+}
 =\frac{(x_k^\vee)'}{\prod_{j\in J}(x_j^\vee)^{[B_{kj}]_+}}.
\]
The second equality uses $B_{jk}d_k=-d_jB_{kj}$.
Together with
\[
 Y^\vee=y_k^\vee\prod_{j\in J}(x_j^\vee)^{-B_{kj}},\qquad
 x_k^\vee b^\vee=1+Y^\vee,
\]
this is exactly the dual exchange relation after the frozen-unit rescaling.
The remaining monodromies are those of Lemma~\ref{lem:dual-product}.
The cotangent lifts preserve the Lie group spin structures and the gradings
normalized at the positive fibers.

Give the $\cN$ factor the Liouville form of \cite[Lemma~4.7]{HS}, scaled by $d_k$, and the other factors their
cotangent forms. Denote the product form by $\lambda_U$.
The factor $d_k$ matches the term $d_k\,\dd\log Y\wedge\dd\log a$
in \eqref{eq:product-omega}.
On a cotangent chart this rescaling sends $(\varphi,p)$ to $(\varphi,d_kp)$,
preserving the circle basis and the positive fiber with its grading and Lie group spin structure.

By \eqref{eq:product-omega}, $\lambda_U-\lambda_X|_U$ is closed.
By Lemma~\ref{lem:dual-product} and the seed-primitive condition in
Lemma~\ref{lem:principal-geometry}, it is exact on the initial seed torus.
The complement of this torus has real codimension at least two.
Transversality to an algebraic stratification moves every loop in $U$ into the torus.
All periods therefore vanish, and therefore we may write
$\lambda_U-\lambda_X|_U=\dd f$ where $f|_{L_+\cap U}=0$.
Now choose by
Lemma~\ref{lem:actual-transfer}  nested Liouville domains containing
neighborhoods of both compact seed tori.
Strict exactness makes Viterbo restriction carry $L_+$ to the product positive brane,
as in \cite[Lemma~5.8]{HS}.
Transitivity identifies the two seed restrictions with their composites through
$\mathrm{Vit}_U$.
The equivalence recalled in Proposition~\ref{prop:node-rank-one} and K\"unneth \cite[Theorem~1.5]{GPS} now give
$
 \Psi:\Fuk(U)\xrightarrow{\sim}\Perf(U^\vee),
$
carrying the product positive brane to $\OO_{U^\vee}$ and compatible with both charts.

Write $i:U^\vee\hookrightarrow X^\vee$, $A=\OO(X^\vee)$ and $S=\OO(U^\vee)$.
The functor $G=\Psi\circ\mathrm{Vit}_U\circ\Phi^{-1}$ sends $A$ to $S$.
By Lemma~\ref{lem:marked-functors}, it is extension of scalars along
its induced homomorphism $g:A\to S$.
Let $r:S\to\OO(T_{\fs}^{\cA})$ be restriction to the first dual seed chart.
The comparisons for $\Psi$ and $\Phi$ at $\fs$ give
$
 r\circ g=r\circ i^*:A\longrightarrow\OO(T_{\fs}^{\cA}).
$
Changing a trivialization of the free rank-one module conjugates by a unit
and hence does not change this equality.

Since $T_{\fs}^{\cA}$ is dense in the integral affine variety $U^\vee$,
the map $r$ is injective. Thus $g=i^*$ and
$\Psi\circ\mathrm{Vit}_U\simeq i^*\circ\Phi$.
Restricting to the second chart gives
$(j_{\fs'}^{\cA})^*\Phi\simeq\beta_{\fs'}^{\cA}\mathrm{Vit}_{\fs'}^{\cA}$.  This establishes the result for a single mutation; in general, it follows by induction.
\end{proof}

\section{Mirror symmetry for sub-peelable matrices}
\label{sec:HMS}

\begin{theorem}\label{thm:sub-peelable}
Let $\widetilde B=(B;C)$ have full rank, where $B$ is
skew-symmetrizable and sub-peelable.
Then $\cU(\widetilde B)$ is smooth affine of finite type and
carries a graded Weinstein structure symplectomorphic to its
usual affine K\"ahler structure. $\cX(\widetilde B^\vee)$ is a smooth Deligne-Mumford stack, and is a smooth affine scheme if $\Z^{m+f}/\widetilde B \Z^m$ is torsion free.   These spaces satisfy equivariant cluster mirror symmetry; in particular, $\Fuk(\cU(\widetilde B)) \cong \Perf(\cX(\widetilde B^\vee))$. 
\end{theorem}

\subsection{Recollections from \cite{HS}}

First, we summarize the basic intertwining of A- and B-model gluings in the inductive step of \cite{HS}.

\begin{proposition}\label{prop:mirror-gluing}
Let $A$ and $S$ be smooth finitely generated $\C$-algebras.
Fix an open immersion
$\Spec S[u^{\pm1},y^{\pm1}]\hookrightarrow\Spec A$
such that $y$ extends to a unit of $A$, and set
\[
 Y=\Spec A\cup_{\Spec S[u^{\pm1},y^{\pm1}]}
                      \Spec S[u^{\pm1},y],\qquad
 D=\{y=0\},\qquad H_0=\{y=0,u=-1\}.
\]
Let $\widetilde D$ be the strict transform of $D$ and put
$Y'=\Bl_{H_0}Y\setminus\widetilde D$.

Suppose given graded Liouville manifolds $V, Z$ and equivalences
$\Fuk(V)\simeq\Perf(A)$ and $\Fuk(Z)\simeq\Perf(S)$.  Suppose in addition given the following geometric data: 
\begin{enumerate}[label=\textup{(G\arabic*)}]
\item\label{hyp:gluing-geometry}
There are Weinstein domains $D_O\subset D_V$ whose boundaries
agree near a common Weinstein stop $\mathcal F$.
Their completions are $O$ and $V$.
There are compatible identifications
\[
 \mathcal F=T^*S^1\times Z,\qquad
 C:=T^*D^2\times Z,\qquad
 M:=\widehat{(V,\mathcal F)\cup_{\mathcal F}C},
\]
and this gluing gives a Weinstein sectorial cover.
\end{enumerate}
This data gives a commutative diagram of Orlov and Viterbo functors: 
\[
 \begin{tikzcd}[column sep=large]
 \Fuk(\mathcal F) \arrow[r, "J"] \arrow[rd, "J_0"]
    &\Fuk(V, \mathcal F )\arrow[d,"q"]\\
 & \Fuk(O, \mathcal{F})
 \end{tikzcd}
\]

Suppose in addition
\begin{enumerate}[label=\textup{(G\arabic*)}] 
\setcounter{enumi}{1}
\item\label{hyp:gluing-marking}

The pair $(O,\mathcal F)$ is identified with
\[
 \bigl(T^*S^1_y,\{(\varphi_y,p_y)=(\pi,-\infty)\}\bigr)
      \times T^*S^1_u\times Z,
\]
with the usual cotangent gradings and Lie group spin structures.  
The indicated point of the stop lies at cotangent infinity.
\end{enumerate}
Then K\"unneth and the 
rank-one calculation of \cite[Corollary~4.10 and Example~4.11]{GLe} give a commutative diagram:
\begin{equation} \label{eq: rank one boosted}
 \begin{tikzcd}[column sep=large]
  \Fuk(O)  \arrow[d, equal] & \Fuk(O, \mathcal F) \arrow[d, equal] \arrow[l] & \Fuk(\mathcal F) \arrow[d, equal] \arrow[l, "J_0"] \arrow[r] & \Fuk(C) \arrow[d, equal] \\
  \Perf(S[u^{\pm1},y^{\pm1}]) & \Perf(S[u^{\pm1},y]) \arrow[l] & \Perf(S[u^{\pm1}]) \arrow[l, "y=0"] \arrow[r, "u=-1"]& \Perf(S) 
 \end{tikzcd}
\end{equation}
Finally we ask: 
\begin{enumerate}[label=\textup{(G\arabic*)}] 
\setcounter{enumi}{2}
\item\label{hyp:gluing-twist}
The stop $\cF$  is swappable in $V$, and its monodromy autoequivalence preserves
$\ker(\Ind q)$.
\item \label{hyp:gluing-coherence} The given isomorphisms intertwine Viterbo restriction $\Fuk(V)\to\Fuk(O)$ with extension
of scalars along $A\to S[u^{\pm1},y^{\pm1}]$.
\end{enumerate}

Then the specified comparisons extend to equivalences
\[
 \Fuk(V,\mathcal F)\simeq\Perf(Y),\qquad
 \Fuk(M)\simeq\Perf(Y').
\]
For $i:D\hookrightarrow Y$ and $j:H_0\hookrightarrow D$, they
give the commutative diagram
\[
\begin{tikzcd}[column sep=large]
 \Fuk(V,\mathcal F)\arrow[d,"\sim"']
 &\Fuk(\mathcal F)\arrow[l,"J"']\arrow[r,"\mathrm{cap}"]
       \arrow[d,"\sim"]
 &\Fuk(C)\arrow[d,"\sim"]\\
 \Perf(Y)&\Perf(D)\arrow[l,"i_*"']\arrow[r,"j^*"]&\Perf(H_0).
\end{tikzcd}
\]
Here $J$ is the Orlov functor and $\mathrm{cap}$ is induced by
the attachment $\mathcal F\subset\partial C$.
For the open immersion $a:\Spec A\hookrightarrow Y$, the
stop-removal functor gives a second commutative square
\[
\begin{tikzcd}[column sep=large]
 \Fuk(V,\mathcal F)\arrow[r,"\mathrm{stop\ removal}"]
       \arrow[d,"\sim"']
 &\Fuk(V)\arrow[d,"\sim"]\\
 \Perf(Y)\arrow[r,"a^*"']&\Perf(A).
\end{tikzcd}
\]

If an object of $\Fuk(V,\mathcal F)$ restricts to the free rank-one
objects on both affine charts, its image on $Y$ is a line bundle.
Tensoring the final equivalence by the inverse of its pullback
to $Y'$ sends the corresponding object of $\Fuk(M)$ to $\OO_{Y'}$.
\end{proposition}
\begin{proof}
We briefly recall some ideas, see \cite[Sec. 5]{HS} for details. 
The map $A\otimes S[u^{\pm1},y]\to S[u^{\pm1},y^{\pm1}]$
is surjective, so $Y$ is separated. It is smooth of finite type,
and $H_0\subset D\subset Y$ are smooth divisors.
Thus the proof of \cite[Proposition~5.3]{HS} applies to this
open immersion.

The functor $J_0$ is conservative by \eqref{eq: rank one boosted}, which follows from \ref{hyp:gluing-marking}.
Together with 
\ref{hyp:gluing-twist}, this verifies the hypotheses of \cite[Theorem~2.1]{HS}; using also \ref{hyp:gluing-coherence}, we deduce: 
\begin{equation}
    \label{eq:B-pull}
    \Ind\Fuk(V,\mathcal F)\simeq
 \Mod_A\mathop{\times}_{\Mod_{S[u^{\pm1},y^{\pm1}]}}
                    \Mod_{S[u^{\pm1},y]}
 \simeq\QCoh(Y)
\end{equation}
Then the geometry assumed in \ref{hyp:gluing-geometry} gives, by sectorial descent \cite{GPS}: 
\begin{equation}\label{eq:B-push}
 \Ind\Fuk(M)\simeq
 \QCoh(Y)\mathop{\coprod}_{\QCoh(D)}\QCoh(H_0)
 \simeq\QCoh(Y').
\end{equation}
Taking compact objects gives the stated equivalences.
The asserted commutativities with $i_*$ and $j^*$ are deduced by pullback and pushout from \eqref{eq: rank one boosted} and 
Hypothesis~\ref{hyp:gluing-coherence}.  The remaining assertions are obvious.
\end{proof}

\subsection{Proof of Theorem \ref{thm:sub-peelable}}

\begin{proposition}\label{prop:sink-step}
Suppose $s$ is a sink of $B$. Put $R=B_{J,J}$, where
$J=\{1,\ldots,m\}\setminus\{s\}$. 
Suppose $R$ (or equivalently $R^\vee$) has a reddening sequence, and that  $\cU_{\pr}(R)$ and $\cU_{\pr}(R^\vee)$ satisfy cluster mirror symmetry when equipped with the geometry of Lemma \ref{lem:principal-geometry}.
Then the same is true for 
$B$ and $B^\vee$. 
\end{proposition}
\begin{proof}
Let
$N=\{j:B_{js}>0\}$ and $I=J\setminus N$.
The remainder of this proof largely amounts to the observation that, while in \cite{HS} we had to inductively assume that $\cU_{\pr}(B_{I, I})$ also satisfied (the geometric and) mirror symmetric hypotheses -- hence we had to assume sink recurrence -- now we can deduce these hypotheses instead from 
Theorem~\ref{thm:freezing}.  We proceed with the details, which are in large part a recapitulation of  \cite{HS}.

Adding a source or sink preserves existence of a reddening
sequence \cite[Theorem~4.5 and Remark~4.6]{CL}.
  Mutation changes sources into sinks, so we discuss the sink case below.

We abbreviate $X:=\cU_{\pr}(B)$ and $X^\vee:=\cU_{\pr}(B^\vee)$. 
The exchange relation at the sink is: 
\begin{equation}\label{eq:sink-relation}
 xb=1+t,\qquad x=x_s,\qquad
 t=y_s\prod_{j\in N}x_j^{B_{js}}.
\end{equation}
Since $1=xb-t$, the space $X$ is covered by its localizations $X_x$ and $X_t$. 

Per Lemma \ref{lem:principal-product},
 $X_x\cong \cU_{\pr}(R)\times(\C^*)^2$.
After freezing $N$, the vertex $s$ has no mutable neighbors. 
The argument of Lemma~\ref{lem:principal-product} gives:
\begin{equation}\label{eq:sink-product}
 X_t\cong\Nn\times Z,\qquad
 X_{xt}\cong(\C^*)^2\times Z,
 \quad Z\cong \cU_{\pr}(B_{I,I})\times(\C^*)^{2|N|}.
\end{equation}
In Lemma~\ref{lem:principal-product} with frozen set $N$, the
retained mutable matrix is $0\oplus B_{I,I}$ and its coefficient
at $s$ is $Y_s=t$. Hence the same integral coordinates give
\[
 \Omega_B=d_s\,\dd\log t\wedge\dd\log x+\Omega_Z.
\]
Applying Lemma~\ref{lem:dual-product} to this splitting gives
\[
 t^\vee=y_s^\vee\prod_{j\in N}(x_j^\vee)^{-B_{sj}},\qquad
 Z^\vee\cong \cU_{\pr}((B^\vee)_{I,I})\times(\C^*)^{2|N|}.
\]
The corresponding dual cover has equation $x^\vee b^\vee=1+t^\vee$.

On the $(x,t)$ factor, the cotangent coordinates of
\eqref{eq:marking} are
\[
 \varphi_y=-\arg t,\qquad \varphi_u=\arg x,\qquad
 p_y=d_s\log|x|,\qquad p_u=d_s\log|t|.
\]
For a rank-one local system $\mathcal L$ on the base torus,
Lemma~\ref{lem:dual-product} identifies
\[
 x^\vee=\operatorname{mon}_{\mathcal L}([S^1_{\varphi_y}]),
 \qquad
 t^\vee=\operatorname{mon}_{\mathcal L}([S^1_{\varphi_u}]).
\]
Here $S^1_{\varphi_y}$ and $S^1_{\varphi_u}$ are the positively
oriented coordinate circles. The seed coordinates of $Z^\vee$
are the monodromies along the remaining base circles, with the
integral coordinate change specified in Lemma~\ref{lem:dual-product}.

Both opens are smooth of finite type. They cover an affine scheme,
so $X= \cU_{\pr}(B)$ is smooth of finite type. The product form on the
second open agrees with the form on $X_x$ by
\eqref{eq:principal-omega}; it is nondegenerate everywhere. 
Indeed, on $X_x$ it is the product form
\[
 \Omega_R+d_s\,\dd\log y_s\wedge\dd\log x.
\]
On $X_t$, the equation $t=xb-1$ gives
\[
 \Omega_B=d_s\frac{\dd b\wedge\dd x}{xb-1}+\Omega_Z.
\]
The function $xb-1$ is invertible on this chart. Thus the first
summand is nondegenerate on the two-dimensional factor, including
at $x=b=0$, and $\Omega_Z$ is nondegenerate on $Z$.

Thus $\Omega_B$ extends globally. The same argument for $B^\vee$ proves
that $X^\vee$ is smooth affine of finite type
and that $\Omega_{B^\vee}$ is a nondegenerate algebraic two-form without poles.

We turn to symplectic geometry. We write 
\[
 \Nn:=\cU_{\pr}([0])
 =\{(x,b,t)\in\C^2\times\C^*:xb=1+t\}.
\]
 We use the Liouville form $\lambda_{\Nn}$ of
\cite[Section~4.2]{HS}.

The hypothesis gives a form $\lambda_R$ on $\cU_{\pr}(R)$.
Using the reddening assumption on $R$, apply Theorem~\ref{thm:freezing}
with frozen set $N$ to obtain a form $\lambda_I$ on
$\cU_{\pr}(B_{I,I})$ with the properties of Lemma~\ref{lem:principal-geometry}.
In the product coordinates of Lemma~\ref{lem:principal-product}, put
\[
 \begin{aligned}
 \lambda_{X_x}
 &=\lambda_R+d_s\bigl(\log|y_s|\,\dd\arg x
                         -\log|x|\,\dd\arg y_s\bigr),\\
 \lambda_Z
 &=\lambda_I+\sum_{j\in N}d_j
       \bigl(\log|s_j|\,\dd\arg x_j
                         -\log|x_j|\,\dd\arg s_j\bigr),\\
 \lambda_{X_t}&=d_s\lambda_{\Nn}+\lambda_Z.
 \end{aligned}
\]
For $X_x$, the coefficient coordinates on the $R$ factor are
$Y_i=y_i x^{B_{si}}$ for $i\in J$.
For $Z$, the coordinates $(x_j,s_j)_{j\in N}$ are those of
\eqref{eq:principal-freezing-coordinates} with frozen set $N$.
The product assertion of Lemma~\ref{lem:principal-geometry}
gives finite-type Morse--Bott Weinstein forms on $X_x$ and $Z$,
with proper strictly exact positive loci.

On $X_x=\cU_{\pr}(R)\times(\C^*)^2_{x,y_s}$, use the circle action
\[
 e^{i\alpha}\cdot(q,x,y_s)=(q,x,e^{i\alpha}y_s),
 \qquad q\in\cU_{\pr}(R).
\]
Its restriction to $X_{xt}$ is
$e^{i\alpha}\cdot(x,t,z)=(x,e^{i\alpha}t,z)$.
Writing $\partial_\eta$ for its generator and
$\omega=\Im\Omega_B$, we find:
\[
 Z_{\lambda_{X_x}}\log|x|=\log|x|,\qquad
 \iota_{\partial_\eta}\omega=d_s\,\dd\log|x|,\qquad
 \lambda_{X_x}(\partial_\eta)=-d_s\log|x|.
\]
Thus the constant $\kappa$ in \cite[Setup~4.8]{HS} is $d_s>0$.

Let $T_Z$ be the initial seed torus of $Z$, and set
$T=(\C^*)^2_{x,t}\times T_Z\subset X_{xt}$.
Both Liouville forms have differential $\omega$, so
$\delta=\lambda_{X_t}-\lambda_{X_x}$ is closed.  We will compute the periods of $\delta$ to show it is exact.
On $T$, compare them with
\[
 \lambda_{\mathrm{tor},B}
 =d_s\bigl(\log|t|\,\dd\arg x-\log|x|\,\dd\arg t\bigr)
       +\lambda_{\mathrm{tor},Z}.
\] 
The seed-torus assertion of Lemma~\ref{lem:principal-geometry}
and the formula for $\lambda_{\Nn}$ in \cite[Section~4.2]{HS}
give
\[
 (\lambda_{X_x}-\lambda_{\mathrm{tor},B})|_T\in\dd C^\infty(T),
 \qquad
 (\lambda_{X_t}-\lambda_{\mathrm{tor},B})|_T\in\dd C^\infty(T).
\]
Consequently $\int_\gamma\delta=0$ for every closed loop
$\gamma\subset T$, including loops in $T_Z$.
The complement of $T$ in the smooth variety $X_{xt}$ is an
algebraic subset of complex codimension at least one.
Every loop in $X_{xt}$ can be moved off this
subset. Hence
\[
 H_1(T;\R)\longrightarrow H_1(X_{xt};\R)
 \quad\text{is surjective},\qquad
 \delta=\dd F\quad\text{on }X_{xt}.
\]
This checks all periods, including those represented by loops in $Z$.
Both forms vanish on the connected positive locus of $X_{xt}$;
normalize $F$ to vanish there.

The compact torus $\T$ of \eqref{eq:principal-compact-action}
preserves these product forms. On $X_t$ it fixes $t$, rotates
$x$, and acts on $Z$ through its frozen torus actions.
For a periodic infinitesimal generator $\xi$ of $\T$,
Cartan's formula gives $\dd(\delta(\xi))=0$.
Its integral around a $\T$-orbit is zero because $\delta=\dd F$;
thus $\delta(\xi)=0$.
The periodic generators span $\operatorname{Lie}\T$ over $\R$.
Hence, for every $\xi\in\operatorname{Lie}\T$:
\[
 \lambda_{X_t}(\xi)-\lambda_{X_x}(\xi)
 =\delta(\xi)=\dd F(\xi)=\xi F=0.
\]
The product Lyapunov functions are also $\T$-invariant.
Together with the cover \eqref{eq:sink-product} and the proper
strictly exact positive loci, these verify all assumptions of
\cite[Setup~4.8]{HS}.

Then \cite[Proposition~4.9]{HS} gives
 \begin{equation}\label{eq:handle-gluing}
 (X,\lambda_X)\cong
 \widehat{(X_x,\mathcal F)\cup_{\mathcal F}(T^*D^2\times Z)},
 \qquad \mathcal F=T^*S^1\times Z.
\end{equation}
Here $\lambda_X$ is a finite-type Morse--Bott Weinstein Liouville
form, $\dd\lambda_X=\omega$, and $L_+\subset X$ is proper with
$\lambda_X|_{L_+}=0$.
The same proposition supplies Liouville domains
$\overline X_{xt}\subset \overline X_x\subset \overline X$ with completions
$X_{xt}$, $X_x$, and $X$.
The boundaries of $\overline X_{xt}$ and $\overline X_x$ agree near their common
Weinstein hypersurface $\mathcal F$.

If $\Psi_{X_x}$ and $\Psi_Z$ are the input Lyapunov functions,
the critical locus of the resulting function on $X$ is
$
 \Crit\Psi_X=\Crit\Psi_{X_x}\ \sqcup\
(\{x=b=0,\ t=-1\}\times\Crit\Psi_Z).
$
The germs along $\Crit\Psi_{X_x}$ are unchanged.
A new component indexed by $C\subset\Crit\Psi_Z$ has normal
index $2+\operatorname{ind}_C\Psi_Z$.

We apply Proposition~\ref{prop:mirror-gluing} to the algebraic data
\[
 A=\OO(X^\vee_{x^\vee}),\qquad S=\OO(Z^\vee),\qquad
 y=x^\vee,\qquad u=t^\vee,
\]
and Liouville manifolds $V=X_x$ and $Z$.

By K\"unneth for Fukaya categories \cite{GPS}, we may tensor the standard mirror symmetry for cotangent bundles \cite[Theorem~1.1]{Ab} with the assumed equivalence for $R$ and the equivalence
for $B_{I,I}$ from Theorem~\ref{thm:freezing}, to get:
\[
 \begin{aligned}
 \Fuk(X_x)
 &\simeq\Perf(\cU_{\pr}(R^\vee))\otimes\Perf(\C[\Z^2])
   \simeq\Perf(X^\vee_{x^\vee}),\\
 \Fuk(Z)
 &\simeq\Perf(\cU_{\pr}((B^\vee)_{I,I}))
             \otimes\Perf(\C[\Z^{2|N|}])
   \simeq\Perf(Z^\vee).
 \end{aligned}
\]
The positive locus in
a torus factor is a cotangent fiber, so the product
equivalences still send the positive Lagrangians to the
structure sheaves.

The variables $y=x^\vee$ and $u=t^\vee$ record monodromy along
$S^1_{\varphi_y}$ and $S^1_{\varphi_u}$, respectively, as computed
above. In these coordinates,
\[
 d_s\bigl(\log|t|\,\dd\arg x-\log|x|\,\dd\arg t\bigr)
      =p_y\,\dd\varphi_y+p_u\,\dd\varphi_u.
\]
The attaching hypersurface is given, in the rank-one model, by
$|x|=\epsilon$ and $\arg t=\pi$ for sufficiently small
$\epsilon>0$. Its $y$ coordinates are
$\varphi_y=\pi$ and $p_y=d_s\log\epsilon<0$.
After completion, the stopped product identification of
\cite[Sections~4.2 and~5.3]{HS} is therefore
\[
 (X_{xt},\mathcal F)\cong
 \bigl(T^*S^1_y,\{\varphi_y=\pi,\ p_y=-\infty\}\bigr)
       \times T^*S^1_u\times Z.
\]
The cotangent gradings and Lie group spin structures are the
ones in \eqref{eq:marking}. The weight $d_s>0$ only rescales
the momenta and does not change the coordinate circles.

In the notation of Proposition ~\ref{prop:mirror-gluing}, we have
\[
 Y=\Spec A\cup_{\Spec S[u^{\pm1},y^{\pm1}]}
                   \Spec S[u^{\pm1},y],\qquad
 D=\{y=0\},\qquad H_0=\{y=0,u=-1\}.
\]
Both $A$ and $S$ are smooth and finitely generated, and $y\in A^*$.
Note that the overlap in the gluing of $Y$ is:
\[
 X^\vee_{x^\vee t^\vee}
   =D\Bigl(\prod_{j\in N}x_j^\vee\Bigr)
          \subset X^\vee_{x^\vee},\qquad
 A\Bigl[\Bigl(\prod_{j\in N}x_j^\vee\Bigr)^{-1}\Bigr]
   =S[u^{\pm1},y^{\pm1}].
\]
Indeed, $y_s^\vee$ is invertible and $-B_{sj}>0$ for $j\in N$,
so inverting $t^\vee$ is equivalent to inverting every
$x_j^\vee$ with $j\in N$.
This is the localization obtained by freezing $N$ in the $R^\vee$
factor of $X^\vee_{x^\vee}$.

By \cite[Lemma~5.2]{HS}, the affine modification is obtained by
gluing $\Spec A$ to the chart $yb=1+u$. This is the dual sink
cover above, so $\Bl_{H_0}Y\setminus\widetilde D\cong \cU_{\pr}(B^\vee)$.

We now verify the enumerated hypotheses of Proposition~\ref{prop:mirror-gluing}.
Hypothesis~\ref{hyp:gluing-geometry} is the domain inclusion and
sectorial gluing of \cite[Proposition~4.9(vi)]{HS}, as written in
\eqref{eq:handle-gluing}.
Hypothesis~\ref{hyp:gluing-marking} is the stopped cotangent
product identification just computed.

For \ref{hyp:gluing-twist}, the circle action
$e^{i\alpha}\cdot(q,x,y_s)=(q,x,e^{i\alpha}y_s)$ is defined on
all of $X_x$. It gives the positive wrapping of $\mathcal F$
in \cite[Proposition~4.9(vi)]{HS}.
In cotangent coordinates it is translation of
$\varphi_y$ by $-\alpha$. Its derivative preserves the cotangent
polarization, and the Lie group spin structure is translation
invariant. Thus its lift to the grading and spin data is periodic.
By \cite[Lemma~2.2]{HS}, its return functor is the identity;
in particular it preserves the kernel required in
\ref{hyp:gluing-twist}.

Finally, \ref{hyp:gluing-coherence}
comes from the compatibility with Viterbo restriction in the proof of  Theorem~\ref{thm:freezing}, which identifies
$
 \Fuk(X_x)\to\Fuk(X_{xt})$ with 
 $S[u^{\pm1},y^{\pm1}]\otL_A-$.
 
These verify all hypotheses of Proposition~\ref{prop:mirror-gluing}, so we get an equivalence
$\Fuk(\cU_{\pr}(B))\simeq\Perf(\cU_{\pr}(B^\vee))$, sending $L_+ \mapsto \OO$.  By construction, the equivalence is compatible with restriction to the seed torus.

Finally, the construction of \cite[Proposition~4.9]{HS} has
$\lambda_X|_{X_x}=\lambda_{X_x}+\dd(\chi F)$ for its cutoff
$\chi$. Hence $\lambda_X-\lambda_{\mathrm{tor},B}$ is exact on
the initial seed torus. By \cite[Lemma~4.13]{HS}, exactness
propagates to every seed torus.
We have proved smoothness, finite type, nondegeneracy of
$\Omega_B$, and properness of $L_+$.
The reddening sequence for $B$ gives the equivariant quadratic
Poisson embedding by Lemma~\ref{lem:quadratic}. (Reddening is also used  to apply Theorem~\ref{thm:freezing} to $R$.)

Thus Lemma~\ref{lem:principal-geometry} applies with
$\alpha_0=\lambda_X$ and compares $\Im\Omega_B$ with an affine
K\"ahler form, with the specified grading. We obtain a cluster-admissible form.

We keep the symplectic form $\Im\Omega_B$ throughout the induction.
Write $\widehat\lambda_X$ for the primitive produced by
Lemma~\ref{lem:principal-geometry}. Then
$
 \dd\widehat\lambda_X=\dd\lambda_X=\Im\Omega_B$ and 
$
 \widehat\lambda_X|_{L_+}=\lambda_X|_{L_+}=0.$
Both primitives are of finite type and weakly Weinstein.
Their difference is exact on the initial seed torus, by the construction
of $\lambda_X$ above and the first conclusion of
Lemma~\ref{lem:principal-geometry}. The period argument of
Lemma~\ref{lem:change-primitive} therefore gives
$\widehat\lambda_X-\lambda_X=\dd f$ with $f|_{L_+}=0$.
Apply that lemma to obtain
\[
 E:\Fuk(X,\lambda_X)\xrightarrow{\sim}\Fuk(X,\widehat\lambda_X),
 \qquad E(L_+)\simeq L_+,
 \qquad \Vit_{\fs,\widehat\lambda_X}E\simeq\Vit_{\fs,\lambda_X}
\]
for the initial seed $\fs$.
If $\Phi_0$ is the equivalence just constructed by gluing, put
$\widehat\Phi=\Phi_0E^{-1}$.
It sends $L_+$ to $\OO$ and satisfies the initial seed comparison.
Lemma~\ref{lem:one-seed-hms} gives the comparisons for all seeds.

The same construction applies to $B^\vee$, with weight $d_s^{-1}$,
using the hypotheses for $R^\vee$.
\end{proof}

\begin{proof}[Proof of Theorem~\ref{thm:sub-peelable}]
First we treat the principal coefficient case.
For the empty mutable matrix, the principal model is a
point, so the geometric assertions are immediate, as is mirror symmetry. Isolated vertices are
then included by Proposition~\ref{prop:sink-step} with $N=\varnothing$.
The result for peelable matrices now follows by induction using Proposition~\ref{prop:sink-step} (using Lemma~\ref{lem:reseeding} to propagate the cluster-admissible form across mutations).  The sub-peelable case follows from Theorem~\ref{thm:freezing}, whose reddening assumption holds by Lemma~\ref{lem:class}.
This establishes the non-equivariant principal case.

As in the proof of Proposition \ref{prop:sink-step}, 
the reddening sequence for $B$, supplied by Lemma~\ref{lem:class}, allows us to appeal to Lemma~\ref{lem:quadratic}, hence Lemma~\ref{lem:principal-geometry},
hence by Proposition~\ref{prop:coefficient-affine} we may compare the $\dd\alpha$ and K\"ahler forms for any full rank $\cU(B;C)$.

The principal mirror $\cU_{\pr}(B^\vee)$ is smooth affine of
finite type, so Proposition~\ref{prop:upper-finite-atlas} gives the stated algebro-geometric properties. 
Theorem~\ref{thm:coefficients} gives equivariant cluster mirror
symmetry for every full-rank $\widetilde B=(B;C)$.
\end{proof}

\newpage

\bibliographystyle{amsplain}
\bibliography{refs}
\end{document}